\documentclass[10pt]{article}
\usepackage{color}
\usepackage{amsmath}
\usepackage{amssymb}
\usepackage{amsfonts}
\usepackage{mathrsfs} 
\usepackage[utf8]{inputenc}
\usepackage[english]{babel}
\usepackage{amsthm}
\usepackage{genyoungtabtikz}
\usepackage{tikz}
\usepackage[normalem]{ulem}

\usepackage[left=3cm,top=2cm,right=3cm,bottom=2cm]{geometry}
\usepackage{graphicx}

\newcommand{\im}{\operatorname{im}}
\newcommand{\Res}{\operatorname{Res}}
\newcommand{\ms}{\scriptstyle}
\newcommand{\LL}{{ \mathbb L \mathbb L }}

\newtheorem{teo}{Theorem}[section]
\newtheorem{lem}[teo]{Lemma}

\newtheorem{pro}[teo]{Proposition}
\newtheorem{defi}[teo]{Definition}
\newtheorem{defi/teo}[teo]{Definition/Theorem}

\newcommand{\Z}{ \mathbb Z}

\newcommand{\F}{ {\mathbb F}}

\newcommand{\e}{{\bf e}}
\newcommand{\f}{{\bf f}}

\newenvironment{dem}{\noindent \textit{Proof:} }{\quad \hfill $\square$}

\newcommand{\tab}{\operatorname{Tab}}

\newcommand{\Hom}{\operatorname{Hom}}
\numberwithin{equation}{section}

\date{} 
\begin{document}
  \Yvcentermath1
  \title{Jucys-Murphy
%    {\color{green}{\sout{{operators}}}}
{\color{black}{elements}}
    for Soergel bimodules}
\author{\sc  steen ryom-hansen\thanks{Supported in part by FONDECYT grant 1171379  } }

%\date{\today}
\maketitle

\pagenumbering{roman}

\begin{abstract}
We produce Jucys-Murphy elements for the diagrammatical category of Soergel bimodules associated with general Coxeter groups, and 
use them to diagonalize the bilinear form 
on the cell modules. This gives rise to an expression for the determinant of the forms and Jantzen type sum formulas.
\end{abstract}
%\newpage

% \tableofcontents
% \newpage
\pagenumbering{arabic}

\section{Introduction}
Soergel bimodules were introduced by Soergel around 1990 in order to prove
the Koszul duality conjectures for category $ \cal O $ 
of a complex simple Lie algebra. At the same time they also gave rise to a new proof
of the Kazhdan-Lusztig conjectures for the composition 
factor multiplicities 
of Verma modules. Although Soergel bimodules are essentially combinatorial objects
defined in terms of the corresponding 
Weyl group, the proof of these conjectures relied on  
intersection homology methods.

\medskip
Some time later Soergel realized that the bimodule theory makes sense for arbitrary Coxeter systems $ (W,S) $ and formulated 
in this general context the {\it Soergel conjecture}, relating the indecomposable bimodules with the Kazhdan-Lusztig basis of the associated 
Hecke algebra. This conjecture implies the positivity conjectures for Kazhdan-Lusztig polynomials.
It was recently proved in a celebrated work by Elias and Williamson. Their proof was entirely algebraic and 
even gave rise to a proof of the Kazhdan-Lusztig conjectures about composition factors of Verma modules that avoids intersection homology. 
%% One of the key point of their work was the statement
%% that a certain explicitly defined element $ \rho$ acts as a so-called {\it Lefschez operator} in the bimodules, with respect to the 
%% natural $ \mathbb Z$-grading on them.

\medskip
In positive characteristic and for $ W $ infinite, the category of Soergel bimodules does not categorify
the Hecke algebra associated with $ (W,S) $ and the relationship to representation theory becomes unclearer. But 
over the last decade, starting with the work of Elias and Khovanov, 
a diagrammatical category $ \cal D $ has been developed that remedies this problem. Apart
from the theoretical advantages of $ \cal D $ over the bimodule category, 
certain calculations are 
more accessible in $ \cal D $ than in the bimodule category. Indeed Williamson's counterexamples
to the explicit bound given in Lusztig's famous conjecture for the representation
theory of algebraic groups in characteristic $p $ used calculations involving $ \cal D$.

\medskip
In this paper we study the diagrammatical category $ \cal D $ starting from the fact, proved by Elias and Williamson, that 
$ \cal D $ is a cellular category. The cellular basis is here a diagrammatical adaption of Libedinsky's light leaves.

\medskip
The cellular algebra approach to the representation theory of the symmetric group and Hecke algebra in type $ A$ 
was pioneered by Murphy.  One of his major insights was the proof that the classical
{\it Jucys-Murphy} elements $ \{ L_i \} $ act lower triangularly on the
cellular basis in these cases. To a large degree it is this triangularity property that explains the importance of the 
Jucys-Murphy elements in modern representation theory. In fact,  
it has been axiomatized by Mathas as the defining property of a  
cellular algebra {\it endowed with a family of JM-elements}. 

\medskip
The objects of $ \cal D $ are expressions over $S$, that is words in the alphabet $ S$, and 
for each such expression $ \underline{w} $ we consider the endomorphism algebra
$ A_{\underline{w}} = { \rm End}_{\cal D} (\underline{w})$. 
It is a cellular algebra in the original sense of Graham and Lehrer.
The first main result of the paper is our Theorem 
{\ref{firstlightleavesJM}} stating 
that $ A_{\underline{w}} $ is a cellular algebra endowed with a family of JM-elements, 
in the sense of Mathas. 
Apart from the intrinsic value of this result, the JM-elements are of interest by themselves. 
They are of the form 
\begin{equation}
\begin{aligned}
  \begin{tikzpicture}[xscale=0.3, yscale=0.4]
\draw[black, very thick](0,0) -- (0,2);
\draw[black, very thick](1,0) -- (1,2);
 \node[circle,draw,black,fill,minimum size=0.5mm, inner sep=0pt]  at (2,1) {};
 \node[circle,draw,black,fill,minimum size=0.5mm, inner sep=0pt]  at (3,1) {};
 \node[circle,draw,black,fill,minimum size=0.5mm, inner sep=0pt]  at (4,1) {};
\draw[blue, very thick](5,0) -- (5,0.5);
\draw[blue, very thick](5,1.5) -- (5,2);
 \node[circle,draw,blue,fill,minimum size=1mm, inner sep=0pt]  at (5,0.5) {};
 \node[circle,draw,blue,fill,minimum size=1mm, inner sep=0pt]  at (5,1.5) {};
\node at (5,3) {\small $ i $'th  position};
 \node[circle,draw,black,fill,minimum size=0.5mm, inner sep=0pt]  at (6,1) {};
 \node[circle,draw,black,fill,minimum size=0.5mm, inner sep=0pt]  at (7,1) {};
 \node[circle,draw,black,fill,minimum size=0.5mm, inner sep=0pt]  at (8,1) {};
\draw[black, very thick](9,0) -- (9,2);
\draw[black, very thick](10,0) -- (10,2);
\node at (-2,1) {$ L_i =  $};
\end{tikzpicture}
\end{aligned}
\end{equation}
and already 
appear in the literature. For example, whenever $ \underline{w} $ corresponds to a reduced expression of 
an element of $W$,   
the Lefschez operator $ \rho $ from Elias and Williamson's proof of Soergel's 
conjecture can be written as a linear combination with coefficients in $ {\mathbb R}^+$ of these $ L_i$'s. 
We hope that this connection between $ \rho $ and the $ L_i$'s can be the inspiration for the construction of 
Lefschez operators for 
other $ \Z$-graded representation theories that arise in algebraic Lie theory.

\medskip
We are interested in the representation theory of the cellular algebras $ A_{\underline{w}}$. 
By the general theory for cellular algebras, for each $ A_{\underline{w}}$ there 
exists a family of cell modules $ \Delta_{\underline{w}}(y ) $,  
each endowed with a
bilinear form $   \langle \cdot, \cdot \rangle_{y } $ such that the simple $ A_{\underline{w}}$-modules are of the form
$ L_{\underline{w}}(y) = \Delta_{\underline{w}}(y )/ {\rm rad}\langle \cdot, \cdot \rangle_{y }$. 
In our particular case $ y $ runs over those elements of $ W $ 
that appear as a subexpression of $ \underline{w} $.

\medskip
In fact, $ \cal D $ is a $ \Z$-graded cellular category and so $ A_{\underline{w}}$ becomes a
$ \Z$-graded cellular algebra in the
sense of Hu and Mathas. 
In the case of the ground field $ \mathbb R  $,
the representation theory of $ A_{\underline{w}}$ was studied by D. Plaza who showed that its 
associated graded decomposition numbers are the Kazhdan-Lusztig polynomials associated with $ (W,S)$.
In characteristic $ p $ the decomposition numbers are the $p$-Kazhdan-Lusztig polynomials.

\medskip
One of the important applications of the triangularity property of the  Jucys-Murphy elements, first due to James and Murphy for the symmetric group, is to 
diagonalize the bilinear form. This gave rise to 
a cancellation free expression for the determinant of the form.  In this paper we show that James and Murphy's idea can be carried out in the setting of $ \cal D $ as 
well, even though the combinatorics of $ \cal D $ is very different from the Young diagram combinatorics of the symmetric group, of course. 
Our formula, given in Theorem {\ref{lightleavesJM}} of the paper, is the following one 
\begin{equation}{\label{detexintrp}}
\det  \langle \cdot, \cdot \rangle_{y} = \pm  \prod_{  \ms \beta > 0, s_{\beta} y > y }   \beta^{ \dim \Delta( s_{\beta} y)} 
\end{equation}
where $ \beta > 0 $ means that $ \beta $ is positive root for $ W $: we must here impose conditions 
for this notion to make sense. Note that the formula 
has much resemblance with the classical Shapovalov formula for Verma modules. 

\medskip
Following Jantzen's original ideas for Verma modules, 
we next go on to construct a filtration of $ A_{\underline{w}}$-submodules 
$ \Delta_{  \underline{w}}(y ) =  \Delta_{  \underline{w}}^0(y) \supseteq \Delta_{  \underline{w}}^1(y) \supseteq \Delta_{  \underline{w}}^2(y) \supseteq \ldots $ 
on $  \Delta_{  \underline{w}}(y ) $ and wish to deduce from ({\ref{detexintrp}})
a sum formula for this filtration valid in the Grothendieck group. 
Here the lack of quasi-heredity of $ A_{\underline{w}}$ turns out to be 
an obstacle for a direct translation of Jantzen's ideas to our setting, 
since the $ \Delta_{  \underline{w}}(y )$'s do not induce a basis of 
the corresponding Grothendieck group. 
We resolve this problem by replacing $ \underline{w} $ by a certain 
subset $ \underline{\pi} $  
of the expressions over $ S $, containing $ \underline{w} $, and satisfying that the corresponding subset $ \pi \subseteq  W $ is an ideal in $ W $ with respect to the 
Bruhat order. This gives rise to another cellular algebra $ A_{\underline{\pi }} $ with cell modules $ \Delta_{  \underline{\pi}}(y )$. 
The algebra $ A_{\underline{\pi }} $ is quasi-hereditary as we show in Theorem {\ref{lightleavesJM}},
and therefore solves the above mentioned problem.

\medskip
The above results are valid for any ground field {\color{black}{$  \Bbbk  $}}
such that $ {\rm char } \,  \Bbbk > 3 $ 
and so for any valuation $ \nu $
on the ground ring $ R $ whose fraction field is {\color{black}{$  \Bbbk  $}} we obtain a filtration 
$ \Delta_{  \underline{\pi}}(y ) =  \Delta_{  \underline{\pi}}^0(y) \supseteq \Delta_{  \underline{\pi}}^1(y) \supseteq \Delta_{  \underline{\pi}}^2(y) \supseteq \ldots $ 
with corresponding sum formula. 
Although the first term of each of these filtrations $ \Delta_{\underline{\pi}}^1(y) $ is always the radical of the form on 
$ \Delta_{\underline{\pi}}(y) $, 
the other terms of the filtration will in general depend on the valuation $ \nu$. 

\medskip
A particularly interesting case of our theory is the case where $ W $ is of type $ A_{n-1} $ 
and $  \Bbbk   = \overline{\F}_p $, since the decomposition numbers for $ A_{ \underline{\pi}}$
in this case are the decomposition numbers for the algebraic group $ Sl_n(\overline{\, \F}_p )$ around the Steinberg weight, by Riche and Williamson's recent 
work. 
In this case, our sum formula looks as follows
\begin{equation}
\sum_{i>0}  [ \Delta_{{\pi},\overline{\F}_p}^i(y)] = \sum_{\beta > 0, s_{\beta} y > y } \nu_{p}(  \beta )  [ \Delta_{{\pi}, \overline{\F}_p}( s_{\beta} y) ] 
\end{equation}
where $ \nu_{p}(  \beta ) $ is the canonical $p $-adic valuation, see Theorem {\ref{Similarly, we have the following Theorem}}.
We give in the end of the paper an example to illustrate how to use this formula to obtain decomposition numbers. 

\medskip
There are several Jantzen type filtrations with associated sum formulas available in the literature, for example Andersen's filtrations for 
tilting modules. At present, we do not know if it makes sense to ask for possible relations between these filtrations and ours.

\medskip
The paper is organized as follows.
In section 2 we describe the basic facts of Soergel bimodules. In section 3
we explain the diagrammatical category $\cal D $ 
and the diagrammatical category $ {\cal D}^{std} $, with an emphasis on the localization methods. Next, 
in section 4, we describe the cellular basis for $ \cal D $, the light leaves basis. In section 5 we introduce the $L_i$'s and 
prove that they are JM-elements. This relies heavily on the previous sections. We moreover show that 
they satisfy a separability condition over the fraction field. In section 6 we obtain via the
JM-elements a first version of the determinant formula and
then, in section 7, we get via Plaza's branching rule the Shapovalov type version of the determinant formula,
mentioned above.
In section 8 we construct the quasi-hereditary algebra $ A_{\pi} $ and use it to deduce the sum formula from the determinant expression.
Finally, in section 9 we give an application of our formula.

\medskip
It is a pleasure to thank David Plaza for useful conversations related to this article.
{\color{black}{It is also a great pleasure to thank the anonymous referee for a 
    long list of comments and suggestions that helped us improve the exposition of the article.}}

\medskip
This paper is published in Journal of Algebra, but this version of the paper differs in some places 
from the published version. After the
publication of the paper H. H. Andersen pointed several inaccuracies out to us, 
especially in the last section 9 where the $ A_2 $-calculations are done, 
and the differences from the published version of the paper are essentially
due to the corresponding corrections. We thank H. H. Andersen for mentioning these
inaccuracies.

\section{The category of Soergel bimodules}
In this section, we briefly explain the basic facts of the category of Soergel bimodules.

\medskip
Let $ (W, S) $ be a Coxeter system. Thus $W$ is the group generated by the finite set $ S $ subject to the relations
$ (s t )^{ m_{st} } = 1$ where the $ m_{st}$'s are fixed numbers $ m_{st} \in \{1,2,3, \ldots, \infty \} $ satisfying 
$ m_{st}= m_{ts}$ and {\color{black}that $ m_{st} = 1 $ iff $ s = t $}. Here 
$m_{st} = \infty $ means that the relation $ (s t )^{ m_{st} } = 1$ is omitted.
We denote by $ < $ the (Chevalley)-Bruhat order on $W$; we use the convention that $1 \in W$ is the smallest
element with respect to $<$.

\medskip
We shall make a careful distinction between elements $ w$ of $W $ and \emph{\color{black}expressions} $ \underline{w}$ in $ S$, 
that is \emph{\color{black}words} of finite length {\color{black} in} the alphabet $ S$,
{\color{black} including the empty word}.
We denote by $ {\bf{exp}}_s$ the set of expressions in $S$ and 
by $ {\bf{rexp}}_s$ the subset of $ {\bf{exp}}_s$
consisting of {\it reduced} expressions in $ S$. 
There is a canonical map $ {\bf{exp}}_s  \rightarrow W, \, \underline{w}
\mapsto w$ which we use tacitly. For example, if $ s \in S$ and $\underline{w} :=ss $ then $\underline{w} \mapsto w= 1 $.

\medskip
Let $ \Bbbk  $ be a field. The setting of the theory of Soergel bimodules is 
a representation $ \mathfrak h $ of $W$ over $ \Bbbk  $ which 
is {\it reflection faithful} in the sense of Soergel, see \cite{So}.
Soergel used in {\it loc. cit.} an extension of a construction of Kac 
to show that a reflection faithful representation of $ (W,S) $ exists 
when $\Bbbk  = \mathbb R $. However, for infinite $W $ and 
for $ \Bbbk $ of positive characteristic, it does not exist in general.

\medskip
% Let $ \Bbbk  $ be a commutative ring and let $ (W, S) $ be a Coxeter system. Let $ \mathfrak h $ be a realization, in the sense of 
% [EW], of $ (W, S) $ over $ \Bbbk $. That is, $ \mathfrak h $ is a finite rang $ \Bbbk$-module
% together with sets $ \{ \alpha_s^{\vee} | s \in S \} \subset {\mathfrak h} $ and $ \{ \alpha_s | s \in S \}  \subset {\mathfrak h}^{\ast}$, 
% subject to certain conditions generalizing the concept of roots systems.

Let $ \mathfrak h $ be a reflection faithful representation of $ (W,S) $ and let 
$ R := \oplus_m S^{m}({\mathfrak h}^{\ast}) $ be the symmetric algebra on $ {\mathfrak h}^{\ast}$. 
The elements $ f $ of $ R $ are the polynomial functions in $ \mathfrak h$. 
The $ W $-action on $ {\mathfrak h} $ induces a $W$-action on $ {\mathfrak h}^{\ast}$, and so  
on $ R$, and for $ s \in S $ the 
corresponding subalgebra of $ s$-invariants is denoted $ R^s $, that is $ R^s := \{ x \in R | sx =x \} $.

\medskip

There is a natural $ \mathbb Z$-grading on $ R$ which is chosen such  
that the elements of $ {\mathfrak h}^{\ast} $ are of degree two.
The $W$-action on $ R $ is degree preserving and so $ R^{s} $ is a $ \mathbb Z$-graded subalgebra of $ R$. 

\medskip
When referring to 'grading', we always mean 
'$ \mathbb Z$-grading'. In general, for a graded commutative ring $ S $ and a graded $ S $-module 
$ M = \oplus_i M^i $ we let $M(k) $ denote the graded $ S$-module obtained by shifting
the degree, that is $ M(k) := \oplus_i M(k)^i  $ where $ M(k)^{i} :=  M^{i+k} $.

\medskip
For $ s \in S $ we define the graded $ R $-bimodule $ B_s := R \otimes_{R^s} R(1) $ 
and for 
$ \underline{w} = s_{i_1}  s_{i_2} \cdots s_{i_k}  \in {\bf{exp}}_s$ 
we define the graded $ R $-bimodule 
\begin{equation}B_{\underline{w} } :=B_{s_{i_1}}   \otimes_R B_{s_{i_2}} \otimes_R  \cdots  \otimes_R B_{s_{i_k}}.
\end{equation}
The action of $ R $ on $ B_{\underline{w} } $ is given by left and right multiplication.
This is the 
\emph{Bott-Samelson bimodule} associated with $ \underline{w}$.
The category of Bott-Samelson bimodules $ {\mathbb B} {\mathbb S}\rm Bim $ is defined as the category whose objects are sums of shifts of 
Bott-Samelson bimodules 
and the category of Soergel bimodules $ {\mathbb S}\rm Bim $ is defined as the category 
whose objects are sums of shifts of \emph{summands} of
Bott-Samelson bimodules, where in both categories
the morphisms are homogeneous $ R  $-bimodule homomorphisms. 
In other words, $ {\mathbb S}\rm Bim $ is the Karoubian envelope of $ {\mathbb B} {\mathbb S}\rm Bim $.

\medskip
Let ${ \cal H}(W,S)$ be the Hecke algebra associated with $(W,S)$. 
It is the 
$ {\mathbb Z}[q, q^{-1} ] $-algebra 
with generators $ \{ H_s |  s \in S \} $ subject to the relations 
$(H_s -q ) (H_s + q^{-1} ) = 0 $
and 
\begin{equation} \underbrace{H_s H_t H_s \ldots}_{m_{st}\, {\rm factors}}   = \underbrace{H_t H_s H_t  \ldots}_{m_{st}\, { \rm factors}} 
\end{equation}
For $ w \in W $ represented by the reduced expression $  w = s_{i_1} s_{i_2} \cdots s_{i_k} $ 
we set $ H_w := H_{i_1} H_{i_2} \cdots H_{i_k} $ which by Matsumoto's Theorem does not depend on the
choice of reduced expression.
Then $ \{ H_w | w \in W \} $ is a $ {\mathbb Z}[q, q^{-1} ] $-basis for $ {\cal H}(W,S)$.  
Let $ \{ \underline{H}_w | w \in W \} $ be the Kazhdan-Lusztig basis for $ {\cal H}(W,S)$ in the normalization used in \cite{So}; for example
$\underline{H}_s  = H_s +q $.
The entries of the base change matrix between the two bases for $ {\cal H}(W,S)$
are the Kazhdan-Lusztig polynomials, by definition. 

\medskip
Let $ \langle  {\mathbb S}{\rm Bim}  \rangle $ be the split Grothendieck group of ${\mathbb S}{\rm Bim} $ and let 
$ [M ] $ denote the class in $ \langle  {\mathbb S}{\rm Bim}  \rangle $ of the object $ M \in {\mathbb S}{\rm Bim} $.
We make $ \langle  {\mathbb S}{\rm Bim}  \rangle $ into a $ {\mathbb Z}[q, q^{-1} ] $-algebra by the rules $[M] [N] := [ M \otimes_R N] $ and 
$ q  [M ] := M[-1] $. 

\medskip
The following important Theorem was shown in \cite{So}. It is known as 
Soergel's categorification Theorem.

\begin{teo}{\label{SoergelCategorificationBimodules}} 
a) For each reduced expression $ \underline{w} \in  {\bf{rexp}}_s  $ of $ w \in W$ 
there is a unique indecomposable bimodule $ B_w $ in $ {\mathbb S}{\rm Bim} $ that occurs as a summand of $ B_{\underline{w}} $ (with multiplicity one) and 
does not occur as a summand of any $ B_{\underline{y}} $ such that $ \underline{y} 
\in  {\bf{rexp}}_s  $ and $ y < w$.
The set $ \{ B_w(m) \mid w \in W, m \in {\mathbb Z}\} $ classifies the indecomposable bimodules in 
$  {\mathbb S}{\rm Bim} $ up to isomorphism. \newline
b) There is a unique algebra homomorphism $F: {\cal H}(W,S) \rightarrow \langle  {\mathbb S}{\rm Bim}  \rangle$ given by 
$\underline{H}_s   \mapsto [ B_s(1) ]$. It is an isomorphism of $ {\mathbb Z}[q, q^{-1}] $-algebras. 
\end{teo}
Soergel's conjecture from \cite{So}
states that if $ \Bbbk  = \mathbb R$, then $ F( \underline{H}_{w}) = [B_w] $.
It implies positivity properties for the Kazhdan-Lusztig polynomials
and if $ W $ is the Weyl group associated with a complex semisimple Lie algebra $ \mathfrak g $, 
it also implies the Kazhdan-Lusztig conjectures for $ \mathfrak g $.
For Weyl groups it was shown by
Soergel in \cite{So1} and for general Coxeter groups it  was recently shown in a seminal work 
by Elias and Williamson, see \cite{EW}.

\section{The diagrammatical categories $ \cal D $ and $ {\cal D}^{std} $}
Let  $ A_{\underline{w}} :=  {\rm End}_{{\mathbb B} {\mathbb S}\rm Bim}( B_{\underline{w} })$.
It is shown in \cite{EW1} that $ {\mathbb B} {\mathbb S}\rm Bim $ is a
cellular category, and in particular 
that $ A_{\underline{w}} $
is a cellular algebra in the sense of Graham and Lehrer. This result is a starting point of
our paper and we need to explain it in some
detail. 

\medskip
The cellularity of $A_{\underline{w}}$ comes from the representation of 
$ {\mathbb B}{ \mathbb S}\rm Bim $ as a diagrammatical
category $ \cal D $.
For general Coxeter systems $ (W,S) $ this diagrammatical category
$ \cal D  $ 
was constructed by Elias and Williamson in \cite{EW1} 
although preliminary versions of $ \cal D $ already appeared in 
\cite{EK} and \cite{E}.  
One advantage of $ \cal D $ over the bimodule category $ {\mathbb B}{ \mathbb S}\rm Bim $
is that it makes sense even for  
infinite $ W $
in positive characteristic.

\medskip
The setting for $ \cal D $ is
a {\it realization} $ \mathfrak h $ of $ (W,S)$, as introduced in section 3 of \cite{EW1}. 
\begin{defi}
A realization of $ (W,S)$ over a commutative ring
$ \Bbbk  $ is a free, finite rank module $ \mathfrak h $ over $ \Bbbk  $ together with sets $ \{ \alpha_s | s \in S \} \subset {\mathfrak h}^{\ast}$ and 
$ \{ \alpha_s^{\vee} | s \in S \} \subset {\mathfrak h}$ such that $ \langle \alpha_s^{\vee}, \alpha_s \rangle = 2 $ and such that the 
rule $ s(v) := v - \langle v, \alpha_s \rangle \alpha_s^{\vee} $ defines a representation of $ W  $ in $ \mathfrak h$. Furthermore, 
a technical 'balancedness' condition should be verified.
\end{defi}

We shall always assume that Demazure surjectivity holds.
\begin{defi}
A realization $ \mathfrak h $ of $ (W,S)$ over 
$ \Bbbk  $ is said to be Demazure surjective if for all $ s \in S $ the evaluation maps $ \alpha_s: \mathfrak h \rightarrow \Bbbk $ and
$ \alpha_s^{\vee}: {\mathfrak h}^{\ast} \rightarrow \Bbbk $ are surjective.
\end{defi}

Given a realization of $ \mathfrak h $ of $ (W,S)$, the graded commutative ring $ R$ is defined as in the bimodule case. 
Let us now explain the various other ingredients of $ \cal D $, as introduced in \cite{EW1}. 

\begin{defi}
A diagram for $ (W,S) $ (or simply a diagram when confusion is not possible) is a finite diagram on a strip ${ \mathbb R  } \times [0,1]$. 
The arcs are decorated with elements of $ S$. The vertices are the points where the arcs end or meet.
The intersection points between the arcs and the upper (lower) border of the strip
${ \mathbb R  } \times \{1\}$ (${ \mathbb R  } \times \{0\}$) are called boundary point and are not vertices.
They define sequences of elements of $ S $, called the top sequence and bottom 
sequence of the diagram.
The arcs all end in vertices or boundary points.
Loops are allowed.
The regions defined by the arcs may be decorated by homogeneous elements of $ R$. All diagrams
are considered up to isotopy.

The diagram is called standard if the only vertices are $ 2 m_{st} $-valent, with incident lines of alternating decorations $ s, t, \ldots, s, t$.

The diagram is called a Soergel diagram if each vertex is either $ 0 $-valent, that is an endpoint of an arc, or 
$ 3$-valent with the three incident arcs of the same colour, or 
$ 2 m_{st} $-valent, with incident arcs of alternating decorations $ s, t, \ldots, s, t$.
The degree of a Soergel diagram is the sum of the degrees of all its vertices and polynomials where a
$ 0 $-valent vertex has degree 1, a $3 $-valent vertex has degree -1 and a $ 2 m_{st} $-valent vertex has degree 0.
\end{defi}

When we draw diagrams, we identify $ S $ with a set of physical colours and 
indicate the decorations of the lines by using those colours.
Below is an example of a Soergel diagram of degree 5 with $ S:=\{ red,blue, green \} $  and
$ m_{ red,blue} = 3$, $ m_{ red,green} = 2$ and $ m_{ blue,green} = 2$.

\medskip
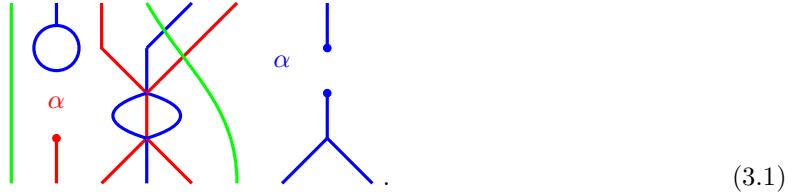
\begin{equation}{\label{picture1}}
\begin{aligned}
  \begin{tikzpicture}[scale=0.6]
\draw [blue,  very thick] (-1,3) circle [radius=0.5];;
\draw [blue,  very thick] (-1,3.5) -- (-1, 4);
\node[circle,red,draw,fill,minimum size=1mm, inner sep=0pt]  at (-1,1) {};
\node at (-1,1.8) {$\color{red} \tiny  \alpha$};
\draw[green,very thick] (-2,0) --(-2,4);
\draw[red,very thick] (-1,0) --(-1,1);
\draw[red,very thick] (0,0) -- (1,1) -- (1,2) -- (0,3) -- (0,4);
\draw[blue,very thick]  (1,0)  -- (1,1);
\draw[blue, very thick] (1,1)  .. controls (0,1.3) and (0,1.7) .. (1,2);
\draw[blue, very thick](1,2) -- (1,3);
\draw[blue,very thick](1,3) -- (2,4);
\draw[red,very thick](2,0) -- (1,1);
\draw[blue, very thick] (1,1)  .. controls (2,1.3) and (2,1.7) .. (1,2);
\draw[red,very thick](1,2) -- (3,4);
\draw[green,very thick](3,0) to[out=90,in=-90+30](1,4);
\draw[blue, very thick] (4,0) -- (5,1) -- (5,2);
\draw[blue, very thick] (6,0) -- (5,1);
\node[circle,blue,draw,fill,minimum size=1mm, inner sep=0pt]  at (5,2) {};
\node[circle,blue,draw,fill,minimum size=1mm, inner sep=0pt]  at (5,3) {};
\node at (4,2.7) {$\color{blue} \tiny \alpha$};
\draw[blue, very thick] (5,3) -- (5,4);
\end{tikzpicture} \, \, .
\end{aligned}
\end{equation}

Standard diagrams are drawn using dashed lines. Below is an example, using the same 
$S $ and $ m$'s as before.
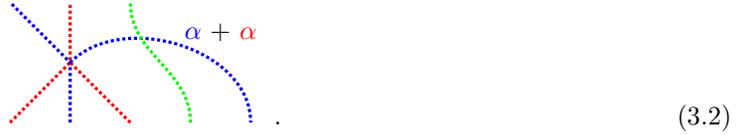
\begin{equation}{\label{picture2}}
\begin{aligned}
  \begin{tikzpicture}[scale=0.8]
\draw[red,very thick, densely dotted] (-2,0) --(-1,1) -- (-1,2);
\draw[blue,very thick, densely dotted] (-1,0) --(-1,1) -- (-2,2);
\draw[red,very thick, densely dotted] (0,0) --(-1,1) ;
\draw[blue,very thick,densely dotted](-1,1) to[out=45,in=90](2,0);
\node at (2-0.5,2-0.5) {$\color{blue} \alpha$ $ +$ $ \color{red} \alpha $};
\draw[green,very thick,densely dotted](1,0) to[out=90,in=-90](0,2);
\end{tikzpicture} \, \, .
\end{aligned}
\end{equation}

We can now give the definition of the category $ \cal D$.

\begin{defi}
Let $ \Bbbk $ be as before. Then 
$ \cal D$ is the $ \Bbbk $-linear monoidal category whose objects are the elements of $ {\bf{exp}}_s$.
{\color{black}{For $ \underline{w}, \underline{z} \in {\bf{exp}}_s $}},
the morphism {\color{black}{space}} $ \Hom_{\cal D} (\underline{w}, \underline{z} )$ is the
free $\Bbbk $-module spanned by 
Soergel diagrams with bottom sequence $ \underline{w} $ and top sequence $ \underline{z}$ modulo a number of 
relations, that we explain below. The composition $ g \circ f  $ of a morphism $ f \in \Hom_{\cal D} (\underline{w}, \underline{z} )$
and a morphism $ g \in \Hom_{\cal D} (\underline{z}, \underline{x} )$ is given by vertical concatenation,
{\color{black}{with $ g $ on top of $ f$}}.

The monoidal structure is given by horizontal concatenation of expressions.
\end{defi}
Let us now explain the relations between the diagrams in $ \Hom_{\cal D} (\underline{w}, \underline{z} )$. 
The isotopy relation has already been mentioned. The other relations 
%% Let $ D$ be a diagram 
%% with upper rightmost boundary point $ P $ of color blue, say, then the multiplication on top of $ D $
%% with the diagram 
%% \begin{equation}{\label{picture3}}
%% \begin{tikzpicture}[ scale=0.3]
%% \draw[very thick, blue](0,0-2) to[out=90,in=180] (1,1-2) to[out=0,in=90] (2,0-2);
%% \node at (0,-0.5-2.5) {$ P$};
%% \end{tikzpicture}$$
%% is an operation called bending downwards the point $ P$.
%% Similarly the operation of bending upwards the lower leftmost point $Q$ is obtained using 
%% $$\begin{tikzpicture}[very thick, scale=0.3]
%% \draw[blue](0,0) to[out=-90,in=180] (1,-1) to[out=0,in=-90] (2,0);
%% \node at (2,0.5+0.5) {$ Q$};
%% \end{tikzpicture}
%% \end{equation}
%% The composite operation of bending downwards all the top points of $ D $ and bending upwards all the bottom points of $ D $, 
%% is called rotation by $180^{\circ}$ of $ D$ and performing it twice is called  rotation by $360^{\circ}$ of $ D$.
%% The cyclicity relation states that rotation by $360^{\circ}$ of $D $ is equal to $ D $ itself for all $ D $.
are divided into the polynomial relations,  
the one-colour, two-colour, three-colour relations and the cyclicity relation.
They are all \emph{local}, involving a small part of the diagram and leaving the rest unchanged.
The 
polynomial and one-colour relations are the following ones (for the colour blue):
\begin{equation}{\label{picture4}}
\begin{aligned}
  \begin{tikzpicture}[scale=0.4]
\node[circle,blue,draw,fill,minimum size=1mm, inner sep=0pt]  at (-3,-2) {};
\node[circle,blue,draw,fill,minimum size=1mm, inner sep=0pt]  at (-3,0) {};
\draw[blue, very thick] (-3,-2) -- (-3,0);
\node at (-2,-1) {$  =$};
\node at (-1,-1) {$\color{blue}  \alpha$};
\end{tikzpicture}
\end{aligned}
\end{equation}

\begin{equation}{\label{picture5}}
\begin{aligned}
  \begin{tikzpicture}[scale=0.5]
\draw[blue,very thick] (-3,-2) -- (-3,0);
\node at (-3.5,-1) {$ f $};
\node at (-1,-1) {$ = $};
\draw[blue,very thick] (1,-2) -- (1,0);
\node at (1.8,-1) {${\color{blue} s}  f $};
\node at (3,-1) {$ +$};
\node[circle,blue,draw,fill,minimum size=1mm, inner sep=0pt]  at (4.5,-1.5) {};
\node[circle,blue,draw,fill,minimum size=1mm, inner sep=0pt]  at (4.5,-0.5) {};
\draw[blue,very thick] (4.5,-2) -- (4.5,-1.5);
\draw[blue,very thick] (4.5,-0.5) -- (4.5,0);
\node at (4.5,-1) {$ {\color{blue} \partial} f $};
\end{tikzpicture}
\end{aligned}
\end{equation}

\begin{equation}{\label{picture6}}
\begin{aligned}
  \begin{tikzpicture}[scale=0.4]
\draw[blue,very thick] (0,0) -- (1,1);
\draw[blue,very thick] (2,0) -- (1,1);
\draw[blue,very thick] (1,1) -- (1,2);
\draw[blue,very thick] (1,2) -- (0,3);
\draw[blue,very thick] (1,2) -- (2,3);
\node at (3, 1.5) {$ = $};
\draw[blue,very thick] (4,0) -- (4,1) -- (4.5,1.5) -- (5,1.5) --(5.5,1) --(5.5,0);
\draw[blue,very thick] (4,3) -- (4,2) -- (4.5,1.5) -- (5,1.5) --(5.5,2) --(5.5,3);
\end{tikzpicture}
\end{aligned}
\end{equation}

\begin{equation}{\label{picture7}}
\begin{aligned}
  \begin{tikzpicture}[scale=0.4]
\draw[blue,very thick] (0,0) -- (0,3);
\draw[blue,very thick] (0,1.5) -- (1,1.5);
\node[circle,blue,draw,fill,minimum size=1mm, inner sep=0pt]  at (1,1.5) {};
\node at (2,1.5) {$ =$};
\draw[blue,very thick] (3,0) -- (3,3);
\end{tikzpicture}
\end{aligned}
\end{equation}

\begin{equation}{\label{picture8}}
\begin{aligned}
  \begin{tikzpicture}[scale=0.4]
\draw [blue,very thick,  very thick] (0,0) circle [radius=0.5];;
\node at (1,0) {$ =$};
\node at (2,0) {$ 0$};
\draw[blue, very thick](0,0.5) -- (0,1);
\draw[blue, very thick](0,-0.5) -- (0,-1);
\end{tikzpicture}
\end{aligned}
\end{equation}
where ${\color{blue} \partial}: R \rightarrow R $ is the Demazure operator given by
$${\color{blue} \partial}(f) := \frac{f- {\color{blue} s}(f)}{{\color{blue} \alpha}}.$$
\noindent
We do not explain the remaining relations, but instead refer the reader to
{\color{black}Definition 5.2 of} \cite{EW1}.

% \medskip
% Although it is of minor importance for the present paper we point out 
% that $ \cal D $ is a graded category 
% as can be easily checked (the relations are all homogeneous).
\medskip
Let us now suppose that $ \Bbbk  $ is a field and that $\mathfrak h $ is
{\color{black}{a}}
\emph{Soergel realization} of $ (W,S) $, in the sense of 
\cite{EW1}, meaning that $ \mathfrak h $ is a faithful representation of $ W $ such that 
Theorem {\ref{SoergelCategorificationBimodules}} holds with respect {\color{black}{to}} $\mathfrak h $. 
Then by \cite{EW1} there is an equivalence of categories 
$ {\cal F }: {\cal D} \rightarrow  {\mathbb B}{\mathbb S}\rm Bim $, defined on objects by 
${\cal F }( \underline{w} ) := B_{\underline{w}}  $, 
where $ B_{ \emptyset} := R$. 
On one-colour morphisms $ \cal F $ is defined as follows
\begin{equation}{\label{picture12}}
\begin{aligned}
  \begin{tikzpicture}[scale=0.5]
  \node[circle,blue,draw,fill,minimum size=1mm, inner sep=0pt]  at (-1,2) {};
\draw[blue, very thick] (-1,0) -- (-1,2);
\node at (0,1) {$  \mapsto $};
\node at (3,1) {$ f \otimes g \mapsto fg $};
\end{tikzpicture}
\end{aligned}
\end{equation}

\begin{equation}{\label{picture13}}
\begin{aligned}
  \begin{tikzpicture}[scale=0.5]
\node[circle,blue,draw,fill,minimum size=1mm, inner sep=0pt]  at (-1,0) {};
\draw[blue, very thick] (-1,0) -- (-1,2);
\node at (0,1) {$  \mapsto $};
\node at (5,1) {$ 1 \mapsto \frac{1}{2}({\color{blue} \alpha} \otimes 1 + 1 \otimes {\color{blue} \alpha} )$};
\end{tikzpicture}
\end{aligned}
\end{equation}

\begin{equation}{\label{picture14}}
\begin{aligned}
  \begin{tikzpicture}[scale=0.5]
\draw[blue,very thick] (0,0) -- (1,1);
\draw[blue,very thick] (2,0) -- (1,1);
\draw[blue,very thick] (1,1) -- (1,2);
\node at (2.5,1) {$  \mapsto $};
\node at (6,1) {$  1 \otimes g \otimes 1 \mapsto {\color{blue}\partial}  g \otimes 1  $};
  \end{tikzpicture}
  \end{aligned}
\end{equation}

\begin{equation}{\label{picture15}}
\begin{aligned}
  \begin{tikzpicture}[scale=0.5]
\draw[blue,very thick] (0,0) -- (0,1);
\draw[blue,very thick] (0,1) -- (1,2);
\draw[blue,very thick] (0,1) -- (-1,2);
\node at (2.5,1) {$  \mapsto $};
\node at (6,1) {$  1 \otimes 1   \mapsto 1  \otimes 1 \otimes 1 .  $};
  \end{tikzpicture}
  \end{aligned}
\end{equation}
(Here, if $ \Bbbk $ is of characteristic 2, 
one should replace 
$ \frac{1}{2}({\color{blue} \alpha} \otimes 1 + 1 \otimes {\color{blue} \alpha} )$ in (\ref{picture13}) by 
$ \Delta_{ \color{blue} s} := {\color{blue} \delta }  \otimes 1 -  1 \otimes {\color{blue} s }( {\color{blue} \delta}) $ where ${\color{blue}  \delta}  $ is defined by 
$ {\color{blue} \partial }( {\color{blue} \delta}) = 1$; it exists because of Demazure surjectivity).
On polynomials, 
$ \cal F $ is defined by mapping $ f \in R $ to multiplication by $ f $ in the slot of $ B_{\underline{w}}$ associated with the region associated with  $f$.
We do not here explain the image of $ \cal F $ on diagrams involving more colours
but refer
%{\color{green}{\sout{{once again}}}}
the reader to {\color{black}{Definition 5.12}} of \cite{EW1}. 

% \medskip
% In order to explain the image under $ \cal F $ of the $2 m$-valent vertex we 
% suppose that $ \underline{w}_1 = {\color{blue} s_{i_1}} {\color{red}  s_{i_2}}  \ldots \in {\bf{exp}}_s$ 
% and $ \underline{w}_2 = {\color{red} s_{i_2}}{ \color{blue} s_{i_1}}  \ldots \in  {\bf{exp}}_s$ satisfy 
% $ w_1 = w_2$. That is $ \underline{w}_1 $ and $\underline{w}_2 $
% express the same element $ w = w_1 =w_2 $ of length $ m= m_{blue, red}  $.
% Then, as already indicated above, $ \cal F $ is defined on the $2 m$-valent vertex as the projection, 
% in $ {\Hom}_{{\mathbb BS}\rm Bim}( B_{\underline{w}}, B_{\underline{w_1}})$,
% onto the indecomposable Soergel bimodule associated with $ w$.

%Thus, by the above statements on the product of $ m$-valent vertices,  the Jones-Wenzl idempotents are mapped to the projector
%n $ {\rm End}_{{\mathbb B} {\mathbb S}\rm Bim}(  B_{\underline{w}})$ on the indecomposable Soergel bimodule $ B_{w}$.
%(We wonder if this is also true for more general rings like for example $\Bbbk = \mathbb Z $).       
This equivalence of categories allows us to answer questions related to ${\mathbb B}{\mathbb S}\rm Bim $ in $ \cal D $.

%By general theory it corresponds to decomposing $ 1 \in  {\rm End}_{\cal D}(  {\underline{w}})$ 
%into its primitive idempotents $ 1 = e_1 + e _2 +\ldots + e_k $, a  problem which on the other hand is known to be intractable in general.

\medskip
We now explain the diagrammatic category $ {\cal D}^{std} $. It is an auxiliary category 
with a particularly simple structure that we shall rely on for our main results.
%% The relationship between $ \cal D$ and $ {\cal D}^{std} $ is comparable to the relationship between the categories of finite 
%% dimensional $ {\mathbb F}_p {\mathfrak S}_n $-modules and finite 
%% dimensional $ {\mathbb Q} S_n $-modules, for $ {\mathfrak S}_n $ the symmetric group on $n$ letters.

\begin{defi}
Let $  \Bbbk $ be a commutative ring. Then 
$ {\cal D}^{std}$ is the additive $ \Bbbk $-linear monoidal category whose objects are direct sums of elements of
$ {\bf{exp}}_s$.
{\color{black} For $ \underline{w}, \underline{z} \in  {\bf{exp}}_s$
the morphism space}
$ \Hom_{{\cal D}^{std}} (\underline{w}, \underline{z} )$ is the $\Bbbk $-module spanned by 
standard diagrams with bottom sequence $ \underline{w} $ and top sequence $ \underline{z}$ modulo a number of 
relations, {\color{black}{that we explain shortly, and for direct sums of such objects we use the
linear extension}}. The composition of morphisms
and the monoidal structure on $ {\cal D}^{std}$ are defined as in $ \cal D$. 
\end{defi}

Let us briefly explain the relations. One should think of them as the relations in $ \cal D $ with the 'lower terms' deleted.
% The tacit relations from $ \cal D $ carry over to $ {\cal D}^{std} $. For example, the diagrams
% $
% \begin{tikzpicture}[densely dotted, very thick, scale=0.3]
% \draw[blue](0,0) to[out=90,in=180] (1,1) to[out=0,in=90] (2,0);
% \end{tikzpicture}
% $
% and
% $
% \begin{tikzpicture}[densely dotted, very thick, scale=0.3]
% \draw[blue](0,0) to[out=-90,in=180] (1,-1) to[out=0,in=-90] (2,0);
% \end{tikzpicture}
% $
% form an adjoint pair such that the cyclicity relation holds with respect to them.
More precisely, the one-colour relations are the following 
\begin{equation}{\label{picture16}}
\begin{aligned}
  \begin{tikzpicture}[densely dotted,very thick, scale=0.6]
\draw[densely dotted, blue,very thick] (0,-2) -- (0,0);
\node at (0.3,-1) {$ f $};
\node at (1,-1) {$ = $};
\draw[blue,very thick] (3,-2) -- (3,0);
\node at (2.5,-1) {$ {\color{blue}s}  f $};
  \end{tikzpicture}
  \end{aligned}
\end{equation}

\begin{equation}{\label{picture17}}
\begin{aligned}
  \begin{tikzpicture}[very thick, densely dotted, scale=0.6]
\draw [densely dotted, blue,very thick] (0,0) circle [radius=0.5];;
\node at (1,0) {$ =$};
\node at (2,0) {$ 1$};
%\draw[blue, very thick](0,0.5) -- (0,1);;
\end{tikzpicture}
\end{aligned}
\end{equation}
\begin{equation}{\label{picture18}}
\begin{aligned}
  \begin{tikzpicture}[densely dotted,very thick, scale=0.35]
\draw[blue](0,3) to[out=-90,in=180] (1,2) to[out=0,in=-90] (2,3);
\draw[blue](0,0) to[out=90,in=180] (1,1) to[out=0,in=90] (2,0);
\node at (3,1.5) {$ =$};
\draw[densely dotted, blue,very thick] (4,0) -- (4,3);
\draw[densely dotted, blue,very thick] (5,0) -- (5,3);
\end{tikzpicture}  \, \, . \end{aligned}
\end{equation}
The two-colour relations are as follows (in the cases $m=2, m= 3$) 
\begin{equation}{\label{picture19}}
\begin{aligned}
  \begin{tikzpicture}[very thick, densely dotted, scale=0.4]
\draw[blue, very thick] (0,0) -- (1,1) -- (0,2);
\draw[red, very thick] (1,0) -- (0,1) -- (1,2);
\node at (2,1) {$=   $};
\draw[blue, very thick] (3,0)  -- (3,2);
\draw[red, very thick] (3.5,0) -- (3.5,2);
\node at (4,0) {$,   $};
\draw[red, very thick] (5,0) -- (6,1) -- (6,3) -- (5,4); 
\draw[red, very thick] (7,0) -- (6,1);
\draw[red, very thick] (6,3) -- (7,4);
\draw[blue, very thick] (6,0) -- (6,1) -- (5,2) -- (6,3) --(6,4);
\draw[blue, very thick] (6,1) --(7,2) -- (6,3)  ;
\node at (8,2) {$ =    $};
\draw[red, very thick] (9,0) -- (9,4); 
\draw[blue, very thick] (10,0) -- (10,4); 
\draw[red, very thick] (11,0) -- (11,4); 
  \end{tikzpicture}
  \end{aligned}
\end{equation}
which implies that the $2m$-valent vertices are idempotents in ${ \cal D }^{std}$.
For an explanation of the remaining relations of $ { \cal D }^{std}$, we
%{\color{green}{\sout{{once again}}}}
refer the reader to {\color{black}{section 4 of}} \cite{EW1}.

\medskip
An important remark is that $ {\cal  D}^{std} $ is a simpler category than $ \cal D$. In fact we have that  
\begin{equation}{\label{standard}}\Hom_{ {\cal D}^{std}}( \underline{x},\underline{y}) = \left\{ \begin{array}{ll}
\Bbbk   & { \rm if } \, \, \, \, x = y \\
0 & { \rm otherwise,} \end{array} \right.
\end{equation}
see the remarks following Theorem 4.8 of \cite{EW1} 
and \cite{EW2}. Indeed, for 
$ \underline{y}, \underline{y^{\prime}} \in {\bf{exp}}_s$ with $ y= y^{\prime} $ there 
is a unique diagram $ Std(\underline{y}, \underline{y^{\prime}}) $ in $ {\cal D}^{std}  $ from $\underline{y} $ to $ \underline{y^{\prime}} $
and the basis for $\Hom_{ {\cal D}^{std}}( \underline{y},\underline{y^{\prime}})$
is given by that diagram.

% \medskip
% In the case of $ \mathfrak h $ being a reflection vector faithful representation of $ (W,S) $ there is a 
% bimodule category equivalent to $ {\cal  D}^{std} $. Indeed, for 
% $ \underline{w} \in {\bf{exp}}_s $ we define the {\it standard} $R$-bimodule $ R_w $ as $ R $ itself
% as a ring, 
% where $f \in  R $ acts on the left via multiplication by $ f$, and acts on the right via  
% multiplication by $ w(f)$. Note that $ R_w $ is not a Soergel bimodule. 
% We define $ {\it Std}{\rm Bim} $ to be the $ \Bbbk$-linear, monoidal, additive category whose objects are sums of standard $R$-bimodules. One can then show that
% $$ \Hom_{ {\it Std}{\rm Bim}   }( R_{\underline{x}}, R_{\underline{y}}    ) = \left\{ \begin{array}{ll}
% \Bbbk   & { \rm if } \, \, \, \, x = y \\
% 0 & { \rm otherwise.} \end{array} \right.$$
% Combining this with ({\ref{standard}}) it follows that  $ {\cal D}^{std}$ and $ {\it Std}{\rm Bim} $ are equivalent categories 
% via the functor that sends $ \underline{x} $ to $ R_{ \underline{x}}$.
\medskip
Let $ Q := Q(R)$ be the quotient field of $ R  $.
We are interested in the corresponding 
category ${\cal D}_{Q} $, in which
elements of $ Q $ are allowed in the regions of the diagrams. 
There is a canonical {\emph{localization}} functor $ {\cal D}  \rightarrow {\cal D}_{Q} $ which
is injective on the $ \Hom$-spaces, since $ \Hom_{\cal D} ({\underline w} ,{\underline y} ) $ is free
over $ R$ (as we shall see shortly). 
The following Theorem, {\color{black}see Definition 5.15 and Theorem 5.16 of} \cite{EW1},
gives a diagrammatical description of the Karoubian envelope of
${\cal D}_{Q} $.
% Indeed, in the bimodule setting we have, 
% working over $ Q $,  that $ B_s = Q_s \oplus Q_{are as1} $. 
% Let $\iota: Q_s \rightarrow B_s $ and $\pi:  B_s \rightarrow Q_s $ denote the splitting maps. 
\begin{teo}
$Kar ({\cal D}_{Q}) $ is  
the additive $ Q $-linear monoidal category whose objects are direct sums
of expressions $ \underline{w} $ over $ S$ as before, but where
each letter $ s \in S $ may now be a normal index or a {\it reflection index}. 
The morphism of the normal indices are drawn normally, whereas the morphisms of the reflection indices are drawn 
using dashed arcs. The morphisms in $ \Hom_{ Kar ({\cal D}_{Q})   }( {\underline{w}}, {\underline{y}}    ) $ are 
spanned over $ Q $ by Soergel diagrams from 
$ {\underline{w}} $ to $  {\underline{y}} $ 
allowing normal arcs and dashed arcs where bivalent vertices, involving the same colour, 
like the following ones, are allowed:
\begin{equation}{\label{picture20}}
\begin{aligned}
  \begin{tikzpicture}[scale=0.4]
    \node at (-2,1) {$ \iota= $};
\draw[densely dotted, blue, very thick](-1,0) -- (-1,1);
\draw[blue, very thick](-1,1) -- (-1,2);
\node at (-2+3,1) {$ \pi= $};
\draw[ blue, very thick](-1+3,0) -- (-1+3,1);
\draw[densely dotted, blue, very thick](-1+3,1) -- (-1+3,2);
\end{tikzpicture} \, \, .
\end{aligned}
\end{equation}
%Thus in particular the cyclicity condition holds for dashed or undashed diagrams. On the other hand, 
The relations in $Kar ({\cal D}_{Q})  $ involving undashed morphisms are the same as in $ {\cal D}_{Q} $ and 
the relations involving dashed morphisms are those of $ {\cal D}_{Q}^{std} $. The relations
involving bivalent vertices are the following ones (in the case
of the colour blue)
% We define $Kar ({\cal D}_{Q}) $ as the 
% We represent them by the following diagrams, containing {\it bivalent} vertices:
% These diagrams are the starting point for 
% the monoidal category $Kar ({\cal D}_{Q}) $, also introduced in [EW1]. The name is justified by the Theorem to follow.
% The objects of $Kar ({\cal D}_{Q}) $ are expressions $ \underline{w} $ over $ S$ just as before, but where
% each letter $ s \in S $ may now be a normal index or a {\it reflection index}. 
% The morphism of the normal indices are drawn normally, whereas the morphisms of the reflection indices are drawn 
% using dashed arcs.
%this time each letter $ s $ either corresponds to the bimodule $ B_s $ or to $ Q_s$. 
%If it corresponds to $ B_s $ the corresponding arc is drawn as usual, 
%if it corresponds to $ Q_s$ the arc is drawn dotted.
% \medskip
% The morphisms in $ \Hom_{ Kar ({\cal D}_{Q})   }( {\underline{w}}, {\underline{y}}    ) $ are 
% spanned over $ Q $ by Soergel diagrams from 
% $ {\underline{w}} $ to $  {\underline{y}} $ 
% allowing normal arcs and dashed arcs.  Bivalent vertices, 
% as illustrated above, are allowed.

% The relations in $Kar ({\cal D}_{Q})  $ involving undashed morphisms are the same as in $ {\cal D}_{Q} $ and 
% the relations involving dashed morphisms are as in $ {\cal D}_{Q}^{std} $.
%Thus in particular the cyclicity condition holds for dashed or undashed diagrams. On the other hand, 
%for diagrams involving bivalent vertices we do not require cyclicity.
\begin{equation}{\label{onecolour 1}}
\begin{aligned}
  \begin{tikzpicture}[scale=0.2]
\draw[densely dotted, blue, very thick](-1,0) -- (-1,2);
\draw[blue, very thick](-1,2) -- (-1,4);
\node[circle,blue,draw,fill,minimum size=1mm, inner sep=0pt]  at (-1,4) {};
\node at (0.9, 2) {$ = 0 $};
\draw[ blue, very thick](-1+6,0) -- (-1+6,2);
\draw[ densely dotted, blue, very thick](-1+6,2) -- (-1+6,4);
\node[circle,,draw,blue,fill,minimum size=1mm, inner sep=0pt]  at (-1+6,0) {};
\node at (6.9, 2) {$ = 0 $};
  \end{tikzpicture}
\end{aligned}
\end{equation}
\begin{equation}{\label{onecolour 2}}
\begin{aligned}
  \begin{tikzpicture}[scale=0.45]
\draw[ blue, very thick](-1,0) -- (-1,2);
\node at (-0.5,1) {$ = $};
\node at (0.5,1) {$  \frac{1}{{\color{blue}\alpha}} $};
\draw[blue, very thick](1,0) -- (1,0.5);
\node[circle,draw,blue,fill,minimum size=1mm, inner sep=0pt]  at (1,0.5) {};
\draw[blue, very thick](1,2) -- (1,1.5);
\node[circle,draw,blue,fill,minimum size=1mm, inner sep=0pt]  at (1,1.5) {};
\node at (2.3,1) {$ \frac{1}{\color{blue}{\alpha}} $};
\node at (1.4,1) {$ + $};
\draw[blue, very thick](3,0) -- (3,0.5);
\draw[densely dotted, blue, very thick](3,1.5) -- (3,0.5);
\draw[blue, very thick](3,2) -- (3,1.5);
\end{tikzpicture}
\end{aligned}
\end{equation}
\begin{equation}{\label{onecolour 3}}
\begin{aligned}
  \begin{tikzpicture}[scale=0.3]
\draw[densely dotted, blue, very thick](-1,0) -- (-1,1);
\draw[blue, very thick](-1,1) -- (-1,2);
\draw[densely dotted, blue, very thick](-1,2) -- (-1,3);
\node at (0,1.5) {$ = $};
\node at (1.5,1.5) {$ {\color{blue} \alpha} $};
\draw[densely dotted, blue, very thick](2.3,0) -- (2.3,3);
  \end{tikzpicture}
 \end{aligned} 
\end{equation}
%% \begin{equation}{\label{onecolour 4}}
%% \begin{tikzpicture}[scale=0.3]
%% \draw[ blue, very thick](-1,0) -- (-1,1);
%% \draw[densely dotted, blue, very thick](-1,1) -- (-1,2);
%% \draw[ blue, very thick](-1,3) -- (-1,2);
%% \node at (0,1.5) {$ = $};
%% \node at (1.5,1.5) {$ \color{blue} \alpha $};
%% \draw[ blue, very thick](2.3,0) -- (2.3,3);
%% \node at (3,1.5) {$ - $};
%% \draw[ blue, very thick](4,0) -- (4,1);
%% \draw[ blue, very thick](4,2) -- (4,3);
%% \node[circle,draw,blue,fill,minimum size=1mm, inner sep=0pt]  at (4,1) {};
%% \node[circle,draw,blue,fill,minimum size=1mm, inner sep=0pt]  at (4,2) {};
%% \end{tikzpicture}
%% \end{equation}

\begin{equation}{\label{onecolour 5}}
\begin{aligned}
  \begin{tikzpicture}[scale=0.3]
\draw[densely dotted, blue, very thick](-1,0) -- (-1,1);
\draw[densely dotted, blue,very thick](-1,1) to[out=90,in=180] (0,2) to[out=0,in=90] (1,1);
\draw[densely dotted, blue, very thick](1,1) -- (1,0);
\node at (2,1) {$ = $};
\node at (3,1) {$ \frac{1}{\color{blue}\alpha} $};
\draw[densely dotted, blue, very thick](-1+5,0) -- (-1+5,1);
\draw[blue,very thick](-1+5,1) to[out=90,in=180] (0+5,2) to[out=0,in=90] (1+5,1);
\draw[densely dotted, blue, very thick](1+5,1) -- (1+5,0);
  \end{tikzpicture}
\end{aligned}  
\end{equation}

\begin{equation}{\label{onecolour 6}}
\begin{aligned}
  \begin{tikzpicture}[scale=0.3]
\draw[densely dotted, blue, very thick](-1,2) -- (-1,1);
\draw[densely dotted, blue,very thick](-1,1) to[out=-90,in=180] (0,0) to[out=0,in=-90] (1,1);
\draw[densely dotted, blue, very thick](1,1) -- (1,2);
\node at (2,1) {$ = $};
\node at (3,1) {$ \frac{1}{\color{blue}\alpha} $};
\draw[densely dotted, blue, very thick](-1+5,2) -- (-1+5,1);
\draw[blue,very thick](-1+5,1) to[out=-90,in=180] (0+5,0) to[out=0,in=-90] (1+5,1);
\draw[densely dotted, blue, very thick](1+5,1) -- (1+5,2);
\end{tikzpicture} \, \, .
\end{aligned}
\end{equation}
\end{teo}

% \medskip
% For the $m$-valent vertex we impose the following relation (in the case $ m= 3$)
% \begin{equation}{\label{mvalent}}
% \begin{tikzpicture}[scale=0.3]
% \draw[densely dotted, red, very thick] (6,0) -- (6,1) -- (7,2) -- (7,4);
% \draw[densely dotted, red, very thick] (8,0) -- (8,1) -- (7,2);
% \draw[densely dotted, blue, very thick] (7,0) -- (7,2) -- (8,3) -- (8,4); 
% \draw[densely dotted, blue, very thick] (7,2) -- (6,3) -- (6,4);
% \node at (9,2) {$ = $};
% \node at (10,2) {$ \frac{1}{\rho} $};
% \draw[densely dotted, red, very thick] (6+5,0) -- (6+5,1);
% \draw[red, very thick] (6+5,1)-- (7+5,2) -- (7+5,3);
% \draw[ red, very thick] (6+5,1)-- (7+5,2) -- (7+5,3);
% \draw[densely dotted, red, very thick](7+5,3) -- (7+5,4);
% \draw[densely dotted,red, very thick] (8+5,0) -- (8+5,1);
% \draw[red, very thick] (8+5,1) -- (7+5,2);
% \draw[blue, very thick] (7+5,1) -- (7+5,2) -- (8+5,3);
% \draw[blue, very thick] (7+5,2) -- (6+5,3);
% \draw[densely dotted, blue, very thick] (7+5,0) -- (7+5,1);
% \draw[densely dotted, blue, very thick] (6+5,3) -- (6+5,4);
% \draw[densely dotted, blue, very thick] (8+5,3) -- (8+5,4); 
% \end{tikzpicture}
% \end{equation}
% where $ \rho $ is the product of the positive roots of the corresponding dihedral group.
% Actually, this can be viewed as the definition of the dashed $ m$-valent vertex. By the balancedness condition mentioned above it 
% is stable under color-preserving rotations.

\medskip
There is a simple description of $Kar ({\cal D}_{Q})  $. 
Indeed, let $ D $ be a diagram in $Kar ({\cal D}_{Q})  $. 
Then applying ({\ref{onecolour 2}}) to all the top and bottom arcs of $ D$ and using the relations
of $Kar ({\cal D}_{Q})  $, together with the relation (\ref{standard}), one gets that
$ D $ becomes a linear combination of diagrams $ S(\underline{y}, \underline{y^{\prime}}) $ in 
which each $S(\underline{y}, \underline{y^{\prime}}) $ is of the form $ Std(\underline{y}, \underline{y^{\prime}}) $
but with the top (dashed) arcs extended to the upper border via
\begin{tikzpicture}[scale=0.2]
\draw[densely dotted, blue, very thick](-1,0) -- (-1,1);
\draw[blue, very thick](-1,1) -- (-1,2);
\end{tikzpicture}
and the bottom (dashed) arcs extended to the lower border via
\begin{tikzpicture}[scale=0.2]
\draw[ blue, very thick](-1,0) -- (-1,1);
\draw[densely dotted, blue, very thick](-1,1) -- (-1,2);
\end{tikzpicture}.
Morever, top boundary vertices 
\begin{tikzpicture}[scale=0.2]
\node[circle,blue,draw,fill,minimum size=1mm, inner sep=0pt]  at (-1,0) {};
\draw[blue, very thick] (-1,0) -- (-1,2);
\end{tikzpicture}
and bottom boundary vertices
\begin{tikzpicture}[scale=0.2]
\node[circle,blue,draw,fill,minimum size=1mm, inner sep=0pt]  at (-1,2) {};
\draw[blue, very thick] (-1,0) -- (-1,2);
\end{tikzpicture}
are allowed.
The details of this description are given in section 5 of \cite{EW1}.

\section{Cellularity of $ {\cal D} $.}
Let us now return to the category $ \cal D$. Elias and Williamson showed in Proposition 6.23 of \cite{EW1} that it is a cellular category, in the the sense
of \cite{Wes}. In particular, for any $ \underline{w} \in {\bf{exp}}_s $ we have that
$ A_{\underline{w}} = { \rm End}_{ \cal D} (B_{\underline{w}}) $ is a cellular algebra.
In this section we explain the various ingredients of the cellular structure of $ \cal D $ according to \cite{EW1}.
The cellular basis itself is a diagrammatical version  
of Libedinsky's light leaves basis, see \cite{Li}.

Let us first recall the original definition of a cellular algebra, as formulated by Graham and Lehrer in \cite{GL}. 
\begin{defi}
Let $ A $ be a finite dimensional algebra 
over a commutative ring $ {\color{black}{ \Bbbk    }}  $. 
Then a cellular basis for $A$ is a triple $ (\Lambda, {\tab},  C ) $ where $ \Lambda $ is a poset, 
$\rm  Tab $ is a function on $ \Lambda $ with values in finite sets
and $ C:\coprod_{\lambda \in \Lambda} {\rm Tab}(\lambda)  \times {\rm Tab}(\lambda) \rightarrow A$
is an injection 
such that 
$$ \{ C_{st}^{\lambda} \, |\,  s,t \in {\rm Tab}(\lambda), \,  \lambda \in \Lambda \} $$
is a {\color{black}{$ \Bbbk$}}-basis for $ A $: the cellular basis for $ A$. The
rule $ (C_{st}^{\lambda})^{\ast} := C_{ts}^{\lambda}$ defines 
a {\color{black}{$ \Bbbk$}}-linear antihomomorphism of $ A$ and, 
finally, the structure constants for $A$ with respect to $ \{ C_{st}^{\lambda} \}$ 
satisfy the following condition with respect to the partial order: for all $ a \in A $ we have
$$  a C_{st}^{\lambda} = \sum_{u \in \tab(\lambda) }
r_{usa} C_{ut}^{\lambda} + \mbox{ lower  terms} $$
where lower terms means a linear combination of $ C^{\mu}_{ab} $ where $ \mu < \lambda $ 
and where $ r_{usa} \in {\color{black}{ \Bbbk    }} $.
\end{defi}

Let us now explain the various ingredients of the cellular basis for  $A_{\underline{w}}$, given in 
\cite{EW1}. For $ {\color{black}{ \Bbbk    }} $ we choose $ R $ itself. 
Suppose that $ \underline{w} = s_{ i_1} s_{ i_2} \cdots s_{ i_k} \in {\bf{exp}}_s $. We then define a {\it subexpression} of $ \underline{w}$
to be a sequence $ \e = (e_1, e_2, \ldots, e_k  ) \in  \{0,1 \}^k $. For each such subexpression 
$ \e $ we define 
$ \underline{w}^{\bf e} := s_{ i_1}^{e_1} s_{ i_2}^{e_2} \cdots s_{ i_k}^{e_k} {\color{black}\in {\bf{exp}}_s}$
and we define 
$ w^{\bf e} \in W $ as {\color{black}{the}} corresponding group element. If $  w^{\bf e} = y $ 
then $\e$ is said to \emph{express} $ y $ and in this case we write $ y \le \underline{w}$.
Note that if $ \underline{w} \in {\bf{rexp}}_s $, then we have that $ y \le \underline{w}$ iff
$ y \le {w}$. 

\medskip
Let $ \underline{w}^{\bf e} = s_{ i_1}^{e_1} s_{ i_2}^{e_2} \cdots s_{ i_k}^{e_k} $. 
For $ \Lambda $ we choose the subset of $ W $ given by
$ \{ {w}^{\bf e} | \, \e \,  \mbox{ is a subexpression of } \underline{w} \} $ with 
the poset structure being induced from the Bruhat order on $ W $. 
For the Tab function we 
choose
\begin{equation} {\rm  Tab}(y) := \{ \e \in \{0,1 \}^k |\,  \e \,\,  \mbox{is a subexpression of } \underline{w} 
  \mbox{ expressing } y \}.
\end{equation}
We are now only left with the definition of the cellular basis itself. We need some combinatorial concepts related to 
subexpressions. 
For $ \underline{w} = s_{ i_1}  \cdots   s_{ i_j} \cdots  s_{ i_k}$ we define 
$ \underline{w}^{\le j} = s_{ i_1}  \cdots s_{ i_j} \in {\bf{exp}}_s $.
Any subexpression $ \bf e $ of $ \underline{w} $ defines a sequence $ (\underline{w}_0, \underline{w}_1, \ldots, \underline{w}_k ) $ 
in $ {\bf{exp}}_s$ via $ \underline{w}_0 = 1 $ and
recursively
\begin{equation}
\underline{ w}_j = 
\left\{ \begin{array}{ll} 
\underline{w}_{j-1} s_j  & \mbox{   if } e_j = 1 \\  \underline{w}_{j-1} &   \mbox{   otherwise. } \end{array} \right.
\end{equation}
This gives rise to a series of {\it symbols} $ T = (t_1, t_2, \ldots, t_k ) \in \{ U, D \}^k $ ($ U$= up, $D $ = down) defined as follows
\begin{equation}
 t_j = 
\left\{ \begin{array}{ll} 
U   & \mbox{   if } w_{j-1} s_j >  w_{j-1} \\  D &   \mbox{   otherwise. } \end{array} \right.
\end{equation}
In particular, we always have that $ t_1 = U$. We merge the sequences $ T $ and $ \e $ into one sequence $ M = (m_1, m_2, \ldots, m_k ) $ 
by concatenating the symbols, that is $ M_j := t_j e_j $. 

With this notation at hand, 
we construct a series of 
morphism $  {\LL}_{\underline{w}, \e, \le j}  \in {\rm Hom}_{\cal D} (B_{\underline{w}^{\le j}}, B_{\underline{w}_j} ) $ as follows. 
We first let $ {\LL}_{\underline{w}, \e, \le 0} $ be the empty diagram. 
Suppose recursively that $ \alpha := {\LL}_{\underline{w}, \e, \le j-1} $ has already been
constructed. Then ${\LL}_{\underline{w}, \e, \le j} $ is obtained from 
$ \alpha $ by first adding on the right a vertical arc of colour $ i_j$. This arc is then further manipulated in a way that depends on 
the value of $ M_j $. The rules are as follows.

\medskip
$ \bullet\, \,$If $ M_j = U0$, then the new arc is terminated with a dot.

\medskip
$ \bullet\,\, $If $ M_j = U1$, then the new arc is continued to the top. 

\medskip
$ \bullet\, \,$If $ M_j = D0$, then $ s_{i_j} $ is in the right descent set of $ w_{k-1} $ and one applies a series of $m$-valent vertices to 
$ \alpha$ such that the result has an arc of colour $ i_j $ to the right of $ \alpha$. Finally a trivalent vertex of
colour $ i_j $ is applied to the final two 
strands.

\medskip
$ \bullet\,\, $If $ M_j = D1$, then we proceed as in case $D0$ but finish with a cap of colour $i_j $.

\medskip
\noindent
The four cases have the following diagrammatical representations: 
\begin{equation}{\label{fourpossibilities1}}
\begin{aligned}
  \begin{tikzpicture}[scale=0.4]
\node at (-2+1,1) {$ U0:$};
\draw (0,0) rectangle (2,2);
\node at (1,1) {$ \alpha $};
\draw[black, very thick] (0.5,-1)-- (0.5, 0);
\draw[black, very thick] (1.5,-1)-- (1.5, 0);
\draw[black, very thick] (1,-1)-- (1, 0);
\draw[black, very thick] (0.5,-1+3)-- (0.5, 0+3);
\draw[black, very thick] (1.5,-1+3)-- (1.5, 0+3);
\draw[black, very thick] (1,-1+3)-- (1, 0+3);
\draw[blue, very thick] (2.5,-1)-- (2.5, 1);
\node[circle,draw,blue,fill,minimum size=1mm, inner sep=0pt]  at (2.5,1) {};
\node at (-1+7,1) {$ U1:$};
\draw (0+7,0) rectangle (2+7,2);
\node at (1+7,1) {$ \alpha $};
\draw[black, very thick] (0.5+7,-1)-- (0.5+7, 0);
\draw[black, very thick] (1.5+7,-1)-- (1.5+7, 0);
\draw[black, very thick] (1+7,-1)-- (1+7, 0);
\draw[black, very thick] (0.5+7,-1+3)-- (0.5+7, 0+3);
\draw[black, very thick] (1.5+7,-1+3)-- (1.5+7, 0+3);
\draw[black, very thick] (1+7,-1+3)-- (1+7, 0+3);
\draw[blue, very thick] (2.5+7 ,-1)-- (2.5+7, 3);
\end{tikzpicture}
\end{aligned}
\end{equation}

\begin{equation}{\label{fourpossibilities2}}
\begin{aligned}
  \begin{tikzpicture}[scale=0.4]
\node at (-2+1,1) {$ D0:$};
\draw (0,0) rectangle (2,2);
\node at (1,1) {$ \alpha $};
\draw[blue, very thick] (0.5,-1)-- (0.5, 0);
\draw[black, very thick] (1.5,-1)-- (1.5, 0);
\draw[black, very thick] (1,-1)-- (1, 0);
\draw[black, very thick] (0.5,-1+3)-- (0.5, 0+3.5);
\draw[blue, very thick] (1.5,+3-1)-- (1.5, 0+3-1);
\draw[black, very thick] (1,-1+3)-- (1, 0+3.5);
\draw[blue, very thick] (2.5,-1)-- (2.5, 3-1);
\draw[blue,very thick](1.5, 0+3-1) to[out=90,in=180] (2,3.5-1) to[out=0,in=90] (2.5,3-1);
\draw[blue, very thick] (2,3.5-1)-- (2, 4.5-1);
\node at (-2+1+7,1) {$ D1:$};
\draw (0+7,0) rectangle (2+7,2);
\node at (1+7,1) {$ \alpha $};
\draw[blue, very thick] (0.5+7,-1)-- (0.5+7, 0);
\draw[black, very thick] (1.5+7,-1)-- (1.5+7, 0);
\draw[black, very thick] (1+7,-1)-- (1+7, 0);
\draw[black, very thick] (0.5+7,-1+3)-- (0.5+7, 0+3.5);
\draw[blue, very thick] (1.5+7,+3-1)-- (1.5+7, 0+3-1);
\draw[black, very thick] (1+7,-1+3)-- (1+7, 0+3.5);
\draw[blue, very thick] (2.5+7,-1)-- (2.5+7, 3-1);
\draw[blue,very thick](1.5+7, 0+3-1) to[out=90,in=180] (2+7,3.5-1) to[out=0,in=90] (2.5+7,3-1);
\end{tikzpicture} \, \, . \end{aligned} 
\end{equation}
Finally we set $ {\LL}_{\underline{w}, \e} := {\LL}_{\underline{w}, \e, \le k} $: this is the {\it light leaves morphism} associated with $ \e $. 
It is a diagrammatical version of the bimodule homomorphism introduced by Libedinsky, see \cite{Li}.

\medskip
For all $ \e \in {\rm  Tab}(y)  $
our chosen $  {\LL}_{\underline{w}, \e}  $ belongs to
$ {\rm Hom}_{\cal D} (B_{\underline{w}}, B_{\underline{w}^{\e}} ) $ where $ {\underline{w}^{\e}} = y $.
We choose a fixed expression $ \underline{y}$ for $ y $ and fix for all appearing $ {\underline{w}^{\e}} $ 
a set of braid moves transforming $ {\underline{w}^{\e}} $ to $ \underline{y}$. We then modify the 
$ {\LL}_{\underline{w}, \e} $ by 
multiplying with the corresponding 
set of $ m$-valent vertices. In this way, all chosen $ {\LL}_{\underline{w}, \e} $ now belong to the same
$ {\rm Hom}_{\cal D} (B_{\underline{w}}, B_{\underline{y}} ) $.

\medskip
It should be noted that although the $ {\LL}_{\underline{w}, \e} $'s depend heavily on the choices of $m$-valent vertices along the way, all choices will do for our results.  

\medskip
For $ \e,  \e_1 \in {\rm Tab}(y) $ we define $\beta:= {\LL}_{\underline{w}, \e} $ and $\beta_1:= {\LL}_{\underline{w}, \e_1} $.
Let $ \beta_1^{\ast} \in  {\rm Hom}_{\cal D} (B_{\underline{y}}, B_{\underline{w}} ) $ be the morphism obtained from $ \beta_1 $ by 
reflection along a horizontal axis. We then define the last ingredient
$ C $ of the cellular basis as $ C^{y}_{ \e \e_1} = {\LL}_{\e, \e_1, y} :=\beta_1^{\ast} \beta \in A_{\underline{w}}$.
\medskip
Let us illustrate the light leaves basis {\color{black}{in}} a couple of simple examples. Let first 
$ W $ be of type $ A_1$, that is $ S = \{ {\color{blue}s}\} $. Then $ W $ has just two elements $ \{ 1,  {\color{blue}s}\} $ 
and all diagrams are
one-coloured, of colour blue. %% Let us take $ \underline{w} := \color{blue}s s$ representing $ 1 $. Then we have 
Let us take $ \underline{w} := \color{blue}s s s$ representing $ \color{blue}s $. Then we have 
$$ \begin{array}{c}{ \rm Tab}(1) = \{ (0,0,0), (1,1,0),  (0,1,1), (1,0,1) \} \\
{ \rm Tab}({\color{blue}s}) = \{ (0,0,1), (1, 1,1), (0,1,0),  (1,0,0)  \}.
\end{array}
$$
The corresponding symbols $ M $ are
\begin{equation}{\label{symbols}}
 \begin{array}{c} \{ (U0,U0,U0), (U1,D1,U0),  (U0,U1,D1), (U1,D0,D1)\} \\
\{(U0,U0,U1), (U1, D1,U1), (U0,U1,D0), (U1,D0,D0)  \}
\end{array}
\end{equation}
and so the corresponding light leaves diagrams in $ { \rm Tab}(1) $ are 
\begin{equation}{\label{lightleavesA1}}
\begin{aligned}
  \begin{tikzpicture}[scale=0.3]
\draw[blue, very thick] (0,0)-- (0, 2);
\node[circle,draw,blue,fill,minimum size=1mm, inner sep=0pt]  at (0,2) {};
\draw[blue, very thick] (1,0)-- (1, 2);
\node[circle,draw,blue,fill,minimum size=1mm, inner sep=0pt]  at (1,2) {};
\draw[blue, very thick] (2,0)-- (2, 2);
\node[circle,draw,blue,fill,minimum size=1mm, inner sep=0pt]  at (2,2) {};

\draw[blue, very thick] (0+4,0)-- (0+4, 2);
\draw[blue, very thick] (1+4,0)-- (1+4, 2);
\draw[blue,very thick](0+4,2) to[out=90,in=180] (0+4+0.5, 2+0.5) to[out=0,in=90] (1+4, 2);
\draw[blue, very thick] (2+4,0)-- (2+4, 2);
\node[circle,draw,blue,fill,minimum size=1mm, inner sep=0pt]  at (2+4,2) {};

\draw[blue, very thick] (0+8+4,0)-- (0+8+4, 2-0.5);
\draw[blue, very thick] (1+8+4,0)-- (1+8+4, 3-0.5);
\draw[blue, very thick] (2+8+4,0)-- (2+8+4, 2-0.5);
\draw[blue,very thick](0+8+4,2-0.5) to[out=90,in=180] (1+8+4, 3-0.5) to[out=0,in=90] (2+8+4, 2-0.5);

\draw[blue, very thick] (0+12-4,0)-- (0+12-4, 2);
\node[circle,draw,blue,fill,minimum size=1mm, inner sep=0pt]  at (0+12-4,2) {};
\draw[blue, very thick] (1+12-4,0)-- (1+12-4, 2);
\draw[blue, very thick] (2+12-4,0)-- (2+12-4, 2);
\draw[blue,very thick](1+12-4,2) to[out=90,in=180] (1.5+12-4, 2.5) to[out=0,in=90] (2+12-4, 2);
  \end{tikzpicture}
  \end{aligned}
\end{equation}
and in $ { \rm Tab}({\color{blue}s}) $ they are 

\begin{equation}{\label{lightleavesA1s}}
\begin{aligned}
  \begin{tikzpicture}[scale=0.3]
\draw[blue, very thick] (0+4,0)-- (0+4, 2);
\draw[blue, very thick] (1+4,0)-- (1+4, 2);
\draw[blue,very thick](0+4,2) to[out=90,in=180] (0+0.5+4, 2+0.5) to[out=0,in=90] (1+4, 2);
\draw[blue, very thick] (2+4,0)-- (2+4, 3);

\draw[blue, very thick] (0+8,0)-- (0+8, 2);
\node[circle,draw,blue,fill,minimum size=1mm, inner sep=0pt]  at (0+8,2) {};
\draw[blue, very thick] (1+8,0)-- (1+8, 2);
\draw[blue, very thick] (2+8,0)-- (2+8, 2);
\draw[blue,very thick](1+8,2) to[out=90,in=180] (1.5+8, 2.5) to[out=0,in=90] (2+8, 2);
\draw[blue, very thick] (1.5+8, 2.5)-- (1.5+8, 2.5+0.5);

\draw[blue, very thick] (0+12-12,0)-- (0+12-12, 2);
\node[circle,draw,blue,fill,minimum size=1mm, inner sep=0pt]  at (0+12-12,2) {};
\draw[blue, very thick] (1+12-12,0)-- (1+12-12, 2);
\node[circle,draw,blue,fill,minimum size=1mm, inner sep=0pt]  at (1+12-12,2) {};
\draw[blue, very thick] (2+12-12,0)-- (2+12-12, 2+1);

\draw[blue, very thick] (0+4+8,0)-- (0+4+8, 1);
\draw[blue, very thick] (1+4+8,0)-- (1+4+8, 1);
\draw[blue,very thick](0+4+8,1) to[out=90,in=180] (0+0.5+4+8, 1+0.5) to[out=0,in=90] (1+4+8, 1);
\draw[blue, very thick] (1+4+8-0.5,1.5)-- (1+4+8-0.5, 3);
\draw[blue,very thick](1+4+8-0.5,2.5) to[out=0,in=90](1+4+8+1,0.5);
\draw[blue, very thick] (1+4+8+1,0.5)-- (1+4+8+1,0);

\end{tikzpicture} \, \, . \end{aligned}
\end{equation}
Thus the cellular basis for $ {\rm End}_{\cal D} ( B_{\underline{w}} ) $ has 32 elements.
Here are four of them corresponding to $ { \rm Tab}(1) \times { \rm Tab}(1)  $
\begin{equation}{\label{lightleavesA11}}
\begin{aligned}
  \begin{tikzpicture}[xscale=0.3, yscale=0.25]
\draw[blue, very thick] (0,0)-- (0, 2);
\node[circle,draw,blue,fill,minimum size=1mm, inner sep=0pt]  at (0,2) {};
\draw[blue, very thick] (1,0)-- (1, 2);
\node[circle,draw,blue,fill,minimum size=1mm, inner sep=0pt]  at (1,2) {};
\draw[blue, very thick] (2,0)-- (2, 2);
\node[circle,draw,blue,fill,minimum size=1mm, inner sep=0pt]  at (2,2) {};

\draw[blue, very thick] (0+4,0)-- (0+4, 2);
\draw[blue, very thick] (1+4,0)-- (1+4, 2);
\draw[blue,very thick](0+4,2) to[out=90,in=180] (0+4+0.5, 2+0.5) to[out=0,in=90] (1+4, 2);
\draw[blue, very thick] (2+4,0)-- (2+4, 2);
\node[circle,draw,blue,fill,minimum size=1mm, inner sep=0pt]  at (2+4,2) {};

\draw[blue, very thick] (0+8,0)-- (0+8, 2-0.5);
\draw[blue, very thick] (1+8,0)-- (1+8, 3-0.5);
\draw[blue, very thick] (2+8,0)-- (2+8, 2-0.5);
\draw[blue,very thick](0+8,2-0.5) to[out=90,in=180] (1+8, 3-0.5) to[out=0,in=90] (2+8, 2-0.5);

\draw[blue, very thick] (0+12,0)-- (0+12, 2);
\node[circle,draw,blue,fill,minimum size=1mm, inner sep=0pt]  at (0+12,2) {};
\draw[blue, very thick] (1+12,0)-- (1+12, 2);
\draw[blue, very thick] (2+12,0)-- (2+12, 2);
\draw[blue,very thick](1+12,2) to[out=90,in=180] (1.5+12, 2.5) to[out=0,in=90] (2+12, 2);

\draw[blue, very thick] (0,0+5)-- (0, -2+5);
\node[circle,draw,blue,fill,minimum size=1mm, inner sep=0pt]  at (0,-2+5) {};
\draw[blue, very thick] (1,0+5)-- (1, -2+5);
\node[circle,draw,blue,fill,minimum size=1mm, inner sep=0pt]  at (1,-2+5) {};
\draw[blue, very thick] (2,0+5)-- (2, -2+5);
\node[circle,draw,blue,fill,minimum size=1mm, inner sep=0pt]  at (2,-2+5) {};

\draw[blue, very thick] (0+4,0+5)-- (0+4, -2+5);
\node[circle,draw,blue,fill,minimum size=1mm, inner sep=0pt]  at (0+4,-2+5) {};
\draw[blue, very thick] (1+4,0+5)-- (1+4, -2+5);
\node[circle,draw,blue,fill,minimum size=1mm, inner sep=0pt]  at (1+4,-2+5) {};
\draw[blue, very thick] (2+4,0+5)-- (2+4, -2+5);
\node[circle,draw,blue,fill,minimum size=1mm, inner sep=0pt]  at (2+4,-2+5) {};

\draw[blue, very thick] (0+8,0+5)-- (0+8, -2+5);
\node[circle,draw,blue,fill,minimum size=1mm, inner sep=0pt]  at (0+8,-2+5) {};
\draw[blue, very thick] (1+8,0+5)-- (1+8, -2+5);
\node[circle,draw,blue,fill,minimum size=1mm, inner sep=0pt]  at (1+8,-2+5) {};
\draw[blue, very thick] (2+8,0+5)-- (2+8, -2+5);
\node[circle,draw,blue,fill,minimum size=1mm, inner sep=0pt]  at (2+8,-2+5) {};

\draw[blue, very thick] (0+12,0+5)-- (0+12, -2+5);
\node[circle,draw,blue,fill,minimum size=1mm, inner sep=0pt]  at (0+12,-2+5) {};
\draw[blue, very thick] (1+12,0+5)-- (1+12, -2+5);
\node[circle,draw,blue,fill,minimum size=1mm, inner sep=0pt]  at (1+12,-2+5) {};
\draw[blue, very thick] (2+12,0+5)-- (2+12, -2+5);
\node[circle,draw,blue,fill,minimum size=1mm, inner sep=0pt]  at (2+12,-2+5) {};
\end{tikzpicture} \, \, . \end{aligned}
\end{equation}
It should be mentioned that in \cite{Li1} Libedinsky gives a non-recursive description of
$ {\rm End}_{\cal D} ( B_{\underline{w}} ) $, for all $ \underline{w} \in  {\bf{rexp}}_s$, but only in 
the $ A_1$ and {\color{black}{$ \tilde{A}_1$-cases}}. 

Let us next consider type $A_2 $, that is $ S := \{ {\color{blue} \alpha} , {\color{ red} \beta } \} $ 
with $ m_{ {\color{blue} \alpha}  {\color{ red} \beta }} = 3 $. Let us take $  {\color{blue} s_1} :={\color{blue} \alpha} $ and 
$ {\color{red} s_2} :={\color{red}\beta} $ and let $ \underline{w} := {\color{blue} s_1} {\color{red} s_2} {\color{blue} s_1} $. 
Then $  {\rm Tab}(1) = \{ (0,0,0), (1,0,1) \} $ with 
symbols $\{ (U0,U0,U0), (U1,U0,D1) \} $ and the corresponding light leaves morphisms are
\begin{equation}{\label{lightleavesA2}}
\begin{aligned}
  \begin{tikzpicture}[xscale=0.4, yscale=0.3]
\draw[blue, very thick] (0,0)-- (0,2);
\node[circle,draw,blue,fill,minimum size=1mm, inner sep=0pt]  at (0,2) {};
\draw[red, very thick] (1,0)-- (1,2);
\node[circle,draw,red,fill,minimum size=1mm, inner sep=0pt]  at (1,2) {};
\draw[blue, very thick] (2,0)-- (2,2);
\node[circle,draw,blue,fill,minimum size=1mm, inner sep=0pt]  at (2,2) {};

\draw[blue, very thick] (0+4,0)-- (0+4,2);
\draw[red, very thick] (1+4,0)-- (1+4,2);
\node[circle,draw,red,fill,minimum size=1mm, inner sep=0pt]  at (1+4,2) {};
\draw[blue,very thick](0+4,2) to[out=90,in=180] (1+4, 3) to[out=0,in=90] (2+4,2);
\draw[blue, very thick] (2+4,0)-- (2+4,2);
\end{tikzpicture} \, \, . \end{aligned}
\end{equation}
We have $  {\rm Tab}({\color{blue} s_1}) = \{ (1,0,0), (0,0,1) \} $ with symbols
$\{ (U1,U0,D0), (U0,U0,U1) \} $ and light leaves morphisms
\begin{equation}{\label{lightleavesA21}}
\begin{aligned}
  \begin{tikzpicture}[xscale=0.4, yscale=0.3]
\draw[blue, very thick] (0+0,0)-- (0+0,2);
\draw[red, very thick] (1+0,0)-- (1+0,2);
\node[circle,draw,red,fill,minimum size=1mm, inner sep=0pt]  at (1+0,2) {};
\draw[blue,very thick](0+0,2) to[out=90,in=180] (1+0, 3) to[out=0,in=90] (2+0,2);
\draw[blue, very thick] (2+0,0)-- (2+0,2);
\draw[blue, very thick] (1,3)-- (1,4);
\draw[blue, very thick] (0+4,0)-- (0+4,2);
\node[circle,draw,blue,fill,minimum size=1mm, inner sep=0pt]  at (0+4,2) {};
\draw[red, very thick] (1+4,0)-- (1+4,2);
\node[circle,draw,red,fill,minimum size=1mm, inner sep=0pt]  at (1+4,2) {};
\draw[blue, very thick] (2+4,0)-- (2+4,4);
\end{tikzpicture} \, \, .  \end{aligned}
\end{equation}
We have $  {\rm Tab}({\color{blue} s_1} {\color{red} s_2}  ) = \{ (1,1,0) \} $ with symbols $ \{ (U1,U1,U0) \} $ and
light leaf 
\begin{equation}
\begin{aligned}  
\begin{tikzpicture}[xscale=0.4, yscale=0.25]
\draw[blue, very thick] (0+0,0)-- (0+0,4);
\draw[red, very thick] (1+0,0)-- (1, 4);
\node[circle,draw,blue,fill,minimum size=1mm, inner sep=0pt]  at (2,2) {};
\draw[blue, very thick] (2,0)-- (2,2);
\end{tikzpicture} \, \, .  \end{aligned}
\end{equation}
Let us finally mention $  {\rm Tab}( {\color{blue} s_1} {\color{red} s_2} {\color{blue} s_1})   = \{ (1,1,1) \} $ that corresponds to the identity map
\begin{equation}
 \begin{aligned} 
\begin{tikzpicture}[xscale=0.4, yscale=0.25]
\draw[blue, very thick] (0+0,0)-- (0+0,4);
\draw[red, very thick] (1+0,0)-- (1, 4);
\draw[blue, very thick] (2,0)-- (2,4);
\end{tikzpicture} \, \, . \end{aligned}
\end{equation}

\section{Jucys-Murphy elements}
Let us fix an expression $ \underline{w} = s_{i_1} s_{i_2} \cdots s_{i_k}  \in {\bf{exp}}_s$.
In this section we consider for $ i \in \{1, 2, \ldots, k \} $ 
the diagram $ L_j \in A_{\underline{w}}$ given by 
\begin{equation}{\label{JM}}
\begin{aligned}  
\begin{tikzpicture}[xscale=0.3, yscale=0.6]
\draw[black, very thick](0,0) -- (0,2);
\draw[black, very thick](1,0) -- (1,2);
 \node[circle,draw,black,fill,minimum size=0.5mm, inner sep=0pt]  at (2,1) {};
 \node[circle,draw,black,fill,minimum size=0.5mm, inner sep=0pt]  at (3,1) {};
 \node[circle,draw,black,fill,minimum size=0.5mm, inner sep=0pt]  at (4,1) {};
\draw[blue, very thick](5,0) -- (5,0.5);
\draw[blue, very thick](5,1.5) -- (5,2);
 \node[circle,draw,blue,fill,minimum size=1mm, inner sep=0pt]  at (5,0.5) {};
 \node[circle,draw,blue,fill,minimum size=1mm, inner sep=0pt]  at (5,1.5) {};
\node at (5,3) {\small $ j $'th  position};
 \node[circle,draw,black,fill,minimum size=0.5mm, inner sep=0pt]  at (6,1) {};
 \node[circle,draw,black,fill,minimum size=0.5mm, inner sep=0pt]  at (7,1) {};
 \node[circle,draw,black,fill,minimum size=0.5mm, inner sep=0pt]  at (8,1) {};
\draw[black, very thick](9,0) -- (9,2);
\draw[black, very thick](10,0) -- (10,2);
\node at (-2,1) {$ L_j =  $};
\end{tikzpicture} \end{aligned}
\end{equation}
where the first colour black refers to $ i_1 $, the second colour black refers to $ i_2 $ and so on, except that 
the blue colour refers to $ i_j$.
In this section we study the left multiplication by $ L_j $ on the light leaves basis $ \{ {\LL}_{ \e, \e_1,  y } \} $. 
It is interesting to note that the elements $ L_j $ already appear, somewhat hidden, in Elias and Williamson's 
proof of Soergel's conjecture, see \cite{EW}. Indeed, 
suppose that $ \underline{w} = s_{i_1} s_{i_2} \cdots s_{i_k} $ and 
set $ \overline{B_{\underline{w} }} := B_{\underline{w} } \otimes_R \mathbb R $ 
where $ \mathbb R $ is made into an $ R $ algebra by mapping {\color{black}{elements}}
of positive degree to zero. 
Then the Lefschetz element $ \rho $ introduced in \cite{EW}, acting {\color{black}{on}}
$ \overline{B_{\underline{w} }}$, satisfies by
Lemma 3.4 of \cite{EW} the formula
\begin{equation}\label{Lefsc}
\rho( - ) = \sum_{j=1}^k  (s_{i_{j-1}} \ldots   s_{i_1} \rho) ( \alpha_{ i_j}^{\vee}) \, {\color{black}{ \chi_j}} \circ \phi_j .
\end{equation}
But \begin{tikzpicture}[scale=0.2]
\node[circle,blue,draw,fill,minimum size=1mm, inner sep=0pt]  at (-1,0) {};
\draw[blue, very thick] (-1,-2) -- (-1,0);
\end{tikzpicture}
is the diagrammatical version of $  \phi_j $
and \begin{tikzpicture}[scale=0.2]
\node[circle,blue,draw,fill,minimum size=1mm, inner sep=0pt]  at (-1,-2) {};
\draw[blue, very thick] (-1,-2) -- (-1,0);
\end{tikzpicture}
is the diagrammatical version of 
${ \color{black}{\chi_j}}$ and thus, in the case where $ \underline{w} \in {\bf{rexp}}_s $ 
we have that $ \rho $ acts in the Bott-Samelson bimodule $ \overline{B_{\underline{w} }}$, 
via an $ {\mathbb R}^+ $-linear combination of our $ L_j$'s.
%{\color{green}{\sout{{We find this observation astonishing.}}}}
The main result of this section is that the set $ \{L_j \mid j=1, \ldots, k \} $
verifies the axiomatic condition, introduced by Mathas in \cite{Ma}, 
for being a set of JM-elements
for $A_{\underline{w}} $ with respect to the light leaves basis.
{\color{black}{On the other hand, the results of \cite{EW} rely on positivity properties for $ \rho $ over the real
    field $ \mathbb R $.
    From this point of view, the formula (\ref{Lefsc}) may indicate 
    a connection between Lefschetz elements and JM-elements that could hold not only for
    Soergel bimodules, but also for other graded cellular algebras. 
    It would be interesting to investigate 
 this possible connection.}} 

\medskip
Let us recall {\color{black}{Mathas' axiomatic condition for JM-elements.}}
\begin{defi}{\label{JMdefinition}}
  Let $ A $ be a cellular algebra over {\color{black}{the commutative ring $ \Bbbk    $}}
  with datum $ (\Lambda, {\rm Tab}, C)$. For each $ \lambda \in \Lambda $, 
suppose that $ {\rm Tab}(\lambda) $ is endowed with a poset structure with order relation $  < $ (depending on $ \lambda$).
Let $ {\bf L } =\{ L_1, L_2, \ldots, L_k \} $ be a commutative family of elements of $A$, such that $ L_i^{\ast} = L_i $ for all $i$.
Suppose that there is a function $ c_a:  \{ 1, 2, \ldots, k \} \rightarrow  {\color{black}{ \Bbbk   }}  $
for each $ a  \in {\rm Tab}(\lambda) $.
Then ${ \bf L}$ is said to be a family 
of JM-elements for $ A $ with content functions $  \{ c_a \} $ if for all $ i $ we have that 
\begin{equation}{\label{JMupppertriangular}}
\ L_i C_{ab}^{\lambda} = c_{a}(i) C_{ab}^{\lambda} + \mbox{lower terms} 
\end{equation}
where lower terms means a linear combination of 
elements from $ \{ C_{a_1 b}^{\lambda} | a_1 < a \} \cup  \{ C_{a_2 b_2}^{\mu} | \mu  < \lambda\} $.
\end{defi}
    {\color{black}{Let us now explain how our cellular algebra $ A_{\underline{w}}$ fits into the
      framework of Definition \ref{JMdefinition}.}}
    For the partial order on $ {\rm Tab}(y ) $
%{\color{green}{\sout{{of our cellular basis for $ A_{\underline{w}}$ }}}}
we use 
the path dominance order $ \preceq $ introduced by Elias and Williamson. Let us explain it. Recall that 
a subexpression $ \bf e $ of $ \underline{w} $ defines a sequence $ (\underline{w}_1, \underline{w}_2, \ldots, \underline{w}_k ) $  
in ${\bf{exp}}_s $ 
and hence also a sequence $ ({w}_1, {w}_2, \ldots, {w}_k ) $ of elements in $ W$.
If $   {\bf f } $ is another subexpression of $ \underline{w} $ with corresponding sequence 
$ ({v}_1, {v}_2, \ldots, {v}_k ) $ in $ W $ then 
we say that $ {\bf e} \preceq {\bf f } $ if $ w_i \leq v_i $ for all $i$.
If $ {\bf e} \preceq {\bf f } $ we also say $ {\LL}_{\underline{w}, \bf e} \preceq {\LL}_{\underline{w}, \bf f} $.

%% \medskip
%% In the $ A1 $ examples ({\ref{lightleavesA1simple_ss}}) and ({\ref{lightleavesA1simple_ss1}}) above we have that 
%% the first diagram is less than the second, with respect to the path dominance order $ \preceq $.  

\medskip
Let us consider the diagrams of ({\ref{lightleavesA1}}). Let us denote them $ D1, D2, D3, D4 $ from left to right and let similarly the diagrams of
({\ref{lightleavesA1s}}) be denoted $ E1, E2, E3, E4 $. Then these diagrams are related via path dominance $ \preceq $ as follows 
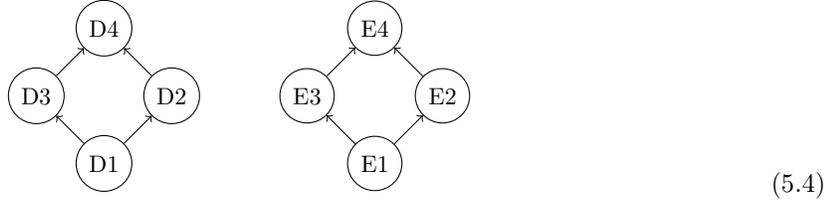
\begin{equation}{\label{pathdominance}}
\begin{aligned}
  \begin{tikzpicture}[scale=.9, transform shape]
\tikzstyle{every node} = [circle, radius=0.6,draw]
\node (a) at (0, 0) {D1};
\node (b) at (1,1) {D2};
\node (c) at (-1,1) {D3};
\node (d) at (0,2) {D4};
\foreach \from/\to in {a/b, a/c, c/d, b/d}
\draw [->] (\from) -- (\to);

\tikzstyle{every node} = [circle, radius=0.6,draw]
\node (a1) at (0+4, 0) {E1};
\node (b1) at (1+4,1) {E2};
\node (c1) at (-1+4,1) {E3};
\node (d1) at (0+4,2) {E4};
\foreach \from/\to in {a1/b1, a1/c1, c1/d1, b1/d1}
\draw [->] (\from) -- (\to);
  \end{tikzpicture}
\end{aligned}  
\end{equation}
with $ D1 $ and $ E1 $ being the smallest and so on.

\medskip
Recall that $ \underline{w} = s_{i_1} s_{i_2} \cdots s_{i_k}  $ is our fixed expression for $ w $. Suppose that 
$ \e = (e_1, e_2, \ldots, e_k ) \in  { \rm Tab}(y) $ and that 
$ T $ is the $a$'th coordinate of the string of symbols given by $ \e $, that is 
$ T \in \{  U0,U1, D0, D1\}$.  We then define for $ a =1, 2, \ldots ,k $ our candidate for the content
function $c_{\e}: \{1,2,\ldots, k \} \rightarrow {\color{black}{{ \Bbbk} =R   }}  $ as
follows
\begin{equation}{\label{candidatcontent}}
c_{\e}(a) := \left\{
\begin{array}{ll}
0 &  \mbox{ if } T= U1, D1  \\
s_{i_1}^{e_1}   \ldots  s_{i_{a-1}}^{e_{a-1}} \alpha_{i_a} & \mbox{ if } T= U0, D0. \\
\end{array} \right.
\end{equation}

We are now in position to state and prove the promised Theorem. The proof uses Proposition 6.6 of
\cite{EW1} which is a main ingredient for showing
the linear independence of the light leaves basis.
\begin{teo}{\label{firstlightleavesJM}}
The set $ {\bf L} := \{ L_1, \ldots, L_k \} $ defined in (\ref{JM}) is a set of JM-elements for $ A_{\underline{w}} $
with respect to the light leaves basis and path dominance order $ \preceq$.
\end{teo}

\begin{dem}
Clearly the $L_i$'s commute and satisfy $ L_i^{\ast} = L_i $ and hence we only need verify the lower triangularity condition 
(\ref{JMupppertriangular}).  
Recall that we have fixed a target object 
$ \underline{y} = s_{j_1} \cdots s_{j_l} \in {\bf{rexp}}_s $ for all light leaves morphisms
in $ \{ \LL_{ \underline{w}, \e } | \e \in {\rm Tab}(y) \} $.
Recall also the localization functor 
$ {\cal D} \rightarrow {\cal D}_Q $; it induces a faithful functor $ {\cal D} \rightarrow Kar({\cal D}_Q) $.
%% For any $ L_{ \underline{w}, \e } $ with $ \e \in {\rm Tab}(y) $ we have that $ G( L_{ \underline{w}, \e }) $ is a $ Q$-linear 
%% combination of standard morphisms $ Std(\underline{z}, \underline{z^{\prime}}) $
%% such that $ {z} = z^{\prime} $ and such that,
%% by the recursive construction of the light leaves, 
%% if $ \underline{z} = \underline{w}^{\bf f} $ then 
%% $ {\bf f} \preceq {\bf e} $ and if $ \underline{z}^{\prime} = \underline{y}^{\bf f_1} $ then ${\bf f_1}$ is a subexpression of $\underline{y} $.
%% In the proof of Proposition 6.6 of \cite{EW} it is shown that 
%% if $ \underline{z} =   \underline{w}^{\bf e}   $ and if $ \underline{z}^{\prime} = \underline{y}$ then 
%% the coeffcient of $ Std(  \underline{z}, \underline{z}^{\prime}  ) $ is a unit in $ Q $. Hence 
%% the $ Q$-spans of the following two sets coincide
%% \begin{equation}{\label{twosets}}
%%   \{  G({\LL}_{ \underline{w}, \bf f })| {\bf f } \preceq \e \} \mbox{       and         }
%%  \{    Std(\underline{w}^{\bf f},\underline{y}^{\bf f_1})    | {\bf f  } \preceq \e, \,  {\bf f_1  } \mbox{ subexpression of } \underline{y} \}. 
%% \end{equation}
Using the description of $ Kar({\cal D}_Q) $ given in the last paragraph of section 3 we have an expansion
\begin{equation}{\label{expanding}}
 {\LL}_{ \underline{w}, \bf e} = \sum_{\substack{ {\f, \f_1 }   \\ {\f }  \preceq \e}} a_{\bf f, f_1 }  S(\underline{w}^{\bf f},\underline{y}^{\bf f_1})      
\end{equation}
where $ {\bf f  } \in \{0,1\}^k $,  
$ {\bf f_1  } \in \{0,1\}^l $ and 
$ a_{\bf f,  f_1 } \in Q $ whereas 
$ S(\underline{w}^{\bf f},\underline{y}^{\bf f_1})  $ is as in section 3.
The condition $ {\bf f }  \preceq \e $ comes from the fact that the target object of
{\color{black}{$ \underline{y}^{\bf f_1}$}} of 
    $ S(\underline{w}^{\bf f},\underline{y}^{\bf f_1}) $ is a 
subexpression of $\underline{y}$ and {\color{black}{$ y = w^\e$}}.
The coefficient $ a_{\e, (1,\ldots,1) } $ is calculated in \cite{EW1} and is shown to be
    nonzero.

We now get that 
\begin{equation}\label{fromthisweget}
 {\LL}_{ \underline{w}, \bf e}  {\color{black}{ L_i}}  = \sum_{\substack{\f, \f_1 \\ \f   \preceq \e }} a_{\bf f, f_1 }  
  S(\underline{w}^{\bf f},\underline{y}^{\bf f_1}) {\color{black}{ L_i}}.
\end{equation}  
On the other hand, for all $ {\bf f } $ and $ {\bf f_1} $ we have that 
\begin{equation}{\label{expanding1}}
 S(\underline{w}^{\bf f},\underline{y}^{\bf f_1}) {\color{black}{ L_i}} = \left\{ \begin{array}{ll}  s_{j_1}^{f_1} s_{j_2}^{f_2}
\ldots s_{j_{i-1}}^{f_{i-1}} \alpha_{i} S(\underline{w}^{\bf f},\underline{y}^{\bf f_1}) &
\mbox{if } f_i = 0  \\
0 & \mbox{if } f_i = 1 \end{array} \right.  
\end{equation}
Indeed if $ f_i = 0 $ then $ S(\underline{w}^{\bf f},\underline{y}^{\bf f_1}) $ 
has a bottom boundary dot \begin{tikzpicture}[scale=0.2]
\node[circle,blue,draw,fill,minimum size=1mm, inner sep=0pt]  at (-1,0) {};
\draw[blue, very thick] (-1,-2) -- (-1,0);
\end{tikzpicture}
at the $i$'th position and so the multiplication by $ L_i $ gives rise to the scalar $ \alpha_i $
which is pulled to the left whereas if $f_i = 1 $ then 
$ S(\underline{w}^{\bf f},\underline{y}^{\bf f_1}) $ has a bottom 
\begin{tikzpicture}[scale=0.2]
  \draw[blue, very thick] (-1,-2) -- (-1,-1);
  \draw[blue,densely dotted, very thick] (-1,-1) -- (-1,0);
\end{tikzpicture}
at the $ i$'th position which is mapped to zero under the multiplication with $ L_i $.
We conclude from this that 
\begin{equation}\label{using}
{\LL}_{ \underline{w}, \bf e} {\color{black}{L_i}} \in {\rm{span}}_Q \{ S(\underline{w}^{\bf f},\underline{y}^{\bf f_1}) \mid
\f    \preceq \e \}.
\end{equation}
Let us now fix an $ \e_0 \in \tab(y) $ and consider
$ {\color{black}{L_i}} {\LL}_{ \e, \e_0, y }     \in A_{\underline{w}}$. Using ({\ref{expanding}}) and (\ref{using}) we have 
that
\begin{equation}\label{beforeappearing}
  L_i {\LL}_{ \e, \e_0, y } = {\color{black} (  {\LL}_{ \underline{w}, \e } L_i)^{\ast}  {\LL}_{\underline{w}, \e_0 }  }
  \in {\rm{span}}_Q \{ S(\underline{w}^{\bf f},\underline{w}^{\bf f_1}) \mid
 \f, \f_1   \preceq \e\}.
  \end{equation}
Let us now consider the $R$-expansion of $ L_i {\LL}_{ \e, \e_0, y } $ in terms of the cellular basis: 
\begin{equation}\label{appearing}
 L_i {\LL}_{ \e, \e_0, y } = \sum_{\e_1 \in \tab(y) } a_{\e_1}{\LL}_{ \e_1, \e_0, y } +  \mbox{lower terms.} 
\end{equation}  
We must prove that for all $ \e_1 $ appearing in (\ref{appearing}) with $ a_{\e_1} \neq 0 $ we have that 
$ \e_1 \preceq \e$. Suppose that this is not the case and choose a maximal counterexample 
$ \e_{max} $ with respect to $ \prec$. 
When expanding $  {\LL}_{ \e_{max}, \e_0, y } $ in terms of 
$ \{ S(\underline{w}^{\bf f},\underline{w}^{\bf f_1}) \} $ we then have that 
the coefficient of $  S(\underline{w}^{\e_{max}},\underline{w}^{\e_0}) $ is nonzero, by the remarks
after ({\ref{expanding}}). On the other hand, when expanding any $ {\LL}_{ \e_1, \e_0, y } $ 
in terms of $ \{ S(\underline{w}^{\bf f},\underline{w}^{\bf f_1}) \} $ only terms with $ \f \preceq \e_1 $
and $ \f_1 \preceq \e_0 $ appear with nonzero coefficient, as follows once again by the remarks 
after ({\ref{expanding}}). Hence, these $ {\LL}_{ \e_1, \e_0, y } $ do not contribute to the coefficient of 
$  S(\underline{w}^{\e_{max}},\underline{w}^{\e_0}) $. 
Finally, for the basis elements $ {\LL}_{ \e_2, \e_3, z } $ appearing in the lower terms of (\ref{appearing}) 
we have that $ z < y $ and so these elements also 
do not contribute to the coefficient of 
$  S(\underline{w}^{\e_{max}},\underline{w}^{\e_0}) $. All in all we get that the
coefficient of $  S(\underline{w}^{\e_{max}},\underline{w}^{\e_0}) $ {\color{black}{on}} the 
right hand side of (\ref{appearing}) is nonzero{\color{black}{,}} which is the desired contradiction, by (\ref{beforeappearing}).
Thus indeed we have proved that 
\begin{equation}\label{appearing1}
 L_i {\LL}_{ \e, \e_0, y } = \sum_{\substack{\e_1 \in \tab(y) \\ \e_1 \preceq \e}} a_{\e_1}{\LL}_{ \e_1, \e_0, y } +  \mbox{lower terms.} 
  \end{equation}  
In order to determine the coefficient $ a_{\bf e }  $ of this expansion, 
we consider the coefficient of $ S(\underline{w}^{\e},\underline{w}^{\e_0}) $
{\color{black}{on}} both sides of (\ref{appearing1}) 
and use ({\ref{expanding1}}).
We find that $ a_{\bf e } = c_{\e}(i) $ and so 
the Theorem is proved.
\end{dem}

\medskip
Let us illustrate the Theorem {\color{black}{with}} the example $ \underline{w} := \color{blue}s s s $, with $ y = \color{blue}s  $ and 
$ {\bf e } = ( 1,0, 0) $. The corresponding light leaves diagram 
\begin{equation}{\label{E4check}}
  \begin{aligned}
\begin{tikzpicture}[scale=0.3]
\draw[blue, very thick] (0+4+8,0)-- (0+4+8, 1);
\draw[blue, very thick] (1+4+8,0)-- (1+4+8, 1);
\draw[blue,very thick](0+4+8,1) to[out=90,in=180] (0+0.5+4+8, 1+0.5) to[out=0,in=90] (1+4+8, 1);
\draw[blue, very thick] (1+4+8-0.5,1.5)-- (1+4+8-0.5, 3);
\draw[blue,very thick](1+4+8-0.5,2.5) to[out=0,in=90](1+4+8+1,0.5);
\draw[blue, very thick] (1+4+8+1,0.5)-- (1+4+8+1,0);
\end{tikzpicture}   \end{aligned}
\end{equation}
is the last one of 
({\ref{lightleavesA1s}}), denoted $ E4 $ in ({\ref{pathdominance}}) and
therefore the maximal {\color{black}{diagram}} appearing in 
({\ref{lightleavesA1s}}) with respect to path dominance. Now, applying $ L_1 $ to it gives the diagram denoted $ E3 $, that is 
\begin{equation}
  \begin{aligned}
  \begin{tikzpicture}[scale=0.3]
\draw[blue, very thick] (0+8,0)-- (0+8, 2);
\node[circle,draw,blue,fill,minimum size=1mm, inner sep=0pt]  at (0+8,2) {};
\draw[blue, very thick] (1+8,0)-- (1+8, 2);
\draw[blue, very thick] (2+8,0)-- (2+8, 2);
\draw[blue,very thick](1+8,2) to[out=90,in=180] (1.5+8, 2.5) to[out=0,in=90] (2+8, 2);
\draw[blue, very thick] (1.5+8, 2.5)-- (1.5+8, 2.5+0.5);
\end{tikzpicture}   \end{aligned}
\end{equation}
which is in accordance with the Theorem since we have $c_{\e}(1) =0 $ for the content function.

\medskip
Let us also calculate the action of $ L_2 $ on ({\ref{E4check}}). Let 
$ {\color{blue} \delta} \in {\mathfrak h}^{\ast}  $
be the element coming from Demazure surjectivity.
Then, using the polynomial and one-colour relations
({\ref{picture4}}), ({\ref{picture5}}), ({\ref{picture7}})
we get that 
the action of $ L_2 $ on ({\ref{E4check}}) is
\begin{equation}
  \begin{aligned}  
\begin{tikzpicture}[xscale=0.32, yscale=0.28]
\draw[blue, very thick] (0+8+4,0)-- (0+8+4, 2-0.5);
\draw[blue, very thick] (1+8+4,0)-- (1+8+4, 1);
\draw[blue, very thick] (2+8+4,0)-- (2+8+4, 2-0.5);
\draw[blue, very thick] (1+8+4,2.5)-- (1+8+4, 4);
\node[circle,draw,blue,fill,minimum size=1mm, inner sep=0pt]  at (1+8+4, 1) {};
\draw[blue,very thick](0+8+4,2-0.5) to[out=90,in=180] (1+8+4, 3-0.5) to[out=0,in=90] (2+8+4, 2-0.5);

\node at (16,2) {$ = $};

\draw[blue, very thick] (0+8+4+6,0)-- (0+8+4+6, 2-0.5);
\draw[blue, very thick] (1+8+4+6,0)-- (1+8+4+6, 0.7);
\draw[blue, very thick] (2+8+4+6,0)-- (2+8+4+6, 2-0.5);
\draw[blue, very thick] (1+8+4+6,2.5)-- (1+8+4+6, 4);
\node[circle,draw,blue,fill,minimum size=1mm, inner sep=0pt]  at (1+8+4+6, 0.7) {};
\draw[blue,very thick](0+8+4+6,2-0.5) to[out=90,in=180] (1+8+4+6, 3-0.5) to[out=0,in=90] (2+8+4+6, 2-0.5);
\draw[blue, very thick](8+4+6, 2-0.5) -- (8+4+6+1, 2-0.5);
\node[circle,draw,blue,fill,minimum size=1mm, inner sep=0pt]  at (8+4+6+1, 2-0.5) {};

\node at (22,2) {$ = $};

\draw[blue, very thick] (0+8+4+6+6,0)-- (0+8+4+6+6, 2-0.5);
\draw[blue, very thick] (1+8+4+6+6,0)-- (1+8+4+6+6, 1.5);
\draw[blue, very thick] (2+8+4+6+6,0)-- (2+8+4+6+6, 2-0.5);
\draw[blue, very thick] (1+8+4+6+6,2.5)-- (1+8+4+6+6, 4);
\draw[blue,very thick](0+8+4+6+6,2-0.5) to[out=90,in=180] (1+8+4+6+6, 3-0.5) to[out=0,in=90] (2+8+4+6+6, 2-0.5);
\draw[blue, very thick](8+4+6+6, 2-0.5) -- (8+4+6+1+6, 2-0.5);
\node at (0+8+4+6+6+0.5, 0.7){${\color{blue} \delta} $};

\node at (28,2) {$ - $};

\draw[blue, very thick] (0+8+4+6+6+6,0)-- (0+8+4+6+6+6, 2-0.5);
\draw[blue, very thick] (1+8+4+6+6+6,0)-- (1+8+4+6+6+6, 1.5);
\draw[blue, very thick] (2+8+4+6+6+6,0)-- (2+8+4+6+6+6, 2-0.5);
\draw[blue, very thick] (1+8+4+6+6+6,2.5)-- (1+8+4+6+6+6, 4);
\draw[blue,very thick](0+8+4+6+6+6,2-0.5) to[out=90,in=180] (1+8+4+6+6+6, 3-0.5) to[out=0,in=90] (2+8+4+6+6+6, 2-0.5);
\draw[blue, very thick](8+4+6+6+6, 2-0.5) -- (8+4+6+1+6+6, 2-0.5);
\node at (0+8+4+6+6+0.5+6+1, 1.2){$ {\color{blue} \scriptsize s \delta} $};

\node at (28+6,2) {$ = $};

%x

\node at (16-5-0.5,2-6+1) {$ {\color{blue} \scriptsize s \delta} -  {\color{blue} \scriptsize  \delta}  $};
\draw[blue, very thick] (0+8+4,0-6)-- (0+8+4, 2-0.5-6);
\draw[blue, very thick] (1+8+4,0-6)-- (1+8+4, 1-6+0.5);
\draw[blue, very thick] (1+8+4-1, 1-6+0.5)-- (1+8+4, 1-6+0.5);
\draw[blue, very thick] (2+8+4,0-6)-- (2+8+4, 2-0.5-6);
\draw[blue, very thick] (1+8+4,2.5-6)-- (1+8+4, 4-6);
\draw[blue,very thick](0+8+4,2-0.5-6) to[out=90,in=180] (1+8+4, 3-0.5-6) to[out=0,in=90] (2+8+4, 2-0.5-6);

\node at (16,2-6) {$ - $};

\draw[blue, very thick] (0+8+4+6,0-6)-- (0+8+4+6, 2-6);
\node[circle,draw,blue,fill,minimum size=1mm, inner sep=0pt]  at (0+8+4+6,2-6) {};
\draw[blue, very thick] (1+8+4+6,0-6)-- (1+8+4+6, 2-6);
\draw[blue, very thick] (2+8+4+6,0-6)-- (2+8+4+6, 2-6);
\draw[blue,very thick](1+8+4+6,2-6) to[out=90,in=180] (1.5+8+4+6, 2.5-6) to[out=0,in=90] (2+8+4+6, 2-6);
\draw[blue, very thick] (1.5+8+4+6, 2.5-6)-- (1.5+8+4+6, 2.5+0.5-6);

\node at (22,2-6) {$ - $};

\draw[blue, very thick] (0+4+20,0-6)-- (0+4+20, 2-6);
\draw[blue, very thick] (1+4+20,0-6)-- (1+4+20, 2-6);
\draw[blue,very thick](0+4+20,2-6) to[out=90,in=180] (0+0.5+4+20, 2+0.5-6) to[out=0,in=90] (1+4+20, 2-6);
\draw[blue, very thick] (2+4+20,0-6)-- (2+4+20, 3-6);

\end{tikzpicture}   \end{aligned}
\end{equation}
which is in accordance with the Theorem since 
$ {\color{blue} \scriptsize s \delta} -  {\color{blue} \scriptsize  \delta} =- {\color{blue} \scriptsize  \alpha} = {\color{blue} \scriptsize s \alpha}  $.

\medskip
Returning to the general situation of Definition {\ref{JMdefinition}}, we define
\begin{equation}\label{tabb} {\color{black}{{\rm Tab}(\Lambda) :=
      \bigcup_{\lambda \in \Lambda} {\rm Tab}(\lambda)}}.
\end{equation}  
We then extend the partial orders $ \le $ on the $ {\rm Tab}(\lambda)$'s to a partial order on ${\rm Tab}(\Lambda) $
by the rule: $ s < t $ if either $ s, t \in {\rm Tab}(\lambda) $ and $ s  < t $ of if $ s \in {\rm Tab}(\lambda), t \in {\rm Tab}(\mu) $ and $ \lambda  < \mu $.
With this notation, Mathas formulated in \cite{Ma} the following {\it separation} condition:
\begin{defi}{\label{JMseparation}}
Let $ {\bf L } =\{ L_1, L_2, \ldots, L_k \} $ be a family of JM-elements as in Definition {\ref{JMdefinition}}. Then 
$  {\bf L } $ is said to \emph{separate} ${\rm Tab}(\Lambda) $ if for all $ s, t \in
{\rm Tab}({\color{black}{\Lambda}}) $ satisfying  $s  < t$, there exist 
$ i $ such that $ c_s(i) \neq c_t(i)$.
\end{defi}
 
An important consequence of the separation condition in the following Proposition, see \cite{Ma}.
\begin{pro}
  Suppose that ${\bf L }$ is a family of JM-elements for a cellular algebra $A$. If ${\bf L }$ satisfies
  the separation condition, then
$ A $ is semisimple.
\end{pro}

Our next result is an application of this Proposition. This result is also implicit in \cite{EW1}.
\begin{pro}\label{isanapp}
  The family of JM-elements ${\bf L }$ from Theorem \ref{firstlightleavesJM} satisfies the separation condition
  over $ Q$. In particular 
$ A_{\underline{w}}$ is semisimple over $Q$.
\end{pro}
\begin{dem}
Let $ \e  =(e_1, \ldots, e_k) \in {\rm Tab}(y) $ and $ {\bf f}= (f_1, \ldots, f_k) \in {\rm Tab}(y_1)$. Let $ i $ be minimal such 
that $ e_i \neq f_i$. Then the symbols for $ \e $ and $ {\bf f} $ at position $ i $ is either $ U $ in both cases or $ D $ in both cases.
Thus, either one of them is $ U0$ and the other $ U1$ or one of them is $ D0 $ and the other $D1$. 
It now follows from (\ref{candidatcontent}) that the contents are different.
\end{dem}

{\color{black}{
\medskip
\noindent 
{\bf Remark}. Note that the proof of the Proposition actually shows a separation condition which is
stronger than the one given in Definition {\ref{JMdefinition}}. Namely, for \emph{any distinct} $ \e, \f 
\in {\rm Tab}(\Lambda)$ there is an $ i $ such that $ c_\e(i ) \neq c_\f(i ) $. }}

%% \medskip 
%% \noindent
%% {\bf Remark}. Even though the diagram categories are not Abelian and 
%% semisimplicity therefore is not defined in them, 
%% the Proposition can also be deduced directly from Theorem {\ref {standarddiagram}} 
%% and ({\ref{standard}}). 

\section{Determinant formula}
Let $ A $ be a cellular algebra over {\color{black}{$ \Bbbk  $}} with cell datum $ (\Lambda, {\rm Tab},  C ) $. 
By the general theory developed in \cite{GL}
there is a canonical family of $A$-modules 
$ \{ \Delta(\lambda) | \lambda \in \Lambda \} $, called the \emph{cell modules} for $A$.
Each cell module $\Delta(\lambda) $
is equipped with a natural symmetric, bilinear, 
$A$-invariant form $ \langle \cdot, \cdot \rangle_{\lambda} $, given in terms of
the cell datum.  
If {\color{black}{$ \Bbbk $}} is a field, the form  
can be used to classify the irreducible modules for $A$. To be precise, 
the irreducible $A$-modules are in correspondence with the subset of $ \Lambda $ given by 
\begin{equation} \Lambda_0 := \{ \lambda \in \Lambda 
  | \langle \cdot, \cdot \rangle_{\lambda} \neq 0 \}
\end{equation}  
where for $ \lambda \in \Lambda_0 $ the corresponding irreducible 
module is given by $ L(\lambda) := \Delta(\lambda)/ {\rm rad }  \langle \cdot, \cdot \rangle_{\lambda}  $.
Here $ {\rm rad }  \langle \cdot, \cdot \rangle_{\lambda} $ is the radical in the usual sense of a bilinear form; it is an $A$-submodule of $ \Delta(\lambda) $ because
of the $ A$-invariance of $   \langle \cdot, \cdot \rangle_{\lambda} $.

\medskip
In the case we are interested in, that is the $R$-algebra $ A_{\underline{w}} $ endowed with the cell
datum explained in the
previous sections, the cell module 
$ \Delta_{\underline{w}}(y) $ has basis $ \{ {\LL}_{ \underline{w}, \bf e} | \, {\bf e} \in {\rm Tab}(y) \} $. 
Let $ \underline{y} \in  {\bf{rexp}}_s  $  be a fixed reduced expression for $ \underline{w}^{ \bf e}  $, the
target object of the elements of $ \{ {\LL}_{ \underline{w}, \bf e} | \, {\bf e} \in {\rm Tab}(y) \} $. 
Then 
the value of the form 
$\langle {\LL}_{ \underline{w}, \bf e}  , {\LL}_{ \underline{w}, \bf e_1} \rangle_{y }  $ is by definition the coefficient of 
$  {\LL}_{ \underline{y}, \bf (1,\ldots, 1), (1,\ldots, 1) } $ when 
$ {\LL}_{ \underline{w}, \bf e} {\LL}_{ \underline{w}, \bf e_1}^{\ast} $ is expanded in the light leaves basis.
Recall that $  {\LL}_{ \underline{y}, \bf (1,\ldots, 1), (1,\ldots, 1) } $ is not uniquely defined,
but since $ \underline{y} \in  {\bf{rexp}}_s $ a possible simple choice is the identity morphism.
Note that the radical of the form is independent of the particular choice
of $  {\LL}_{ \underline{y}, \bf (1,\ldots, 1), (1,\ldots, 1) } $. 

\medskip
Let us illustrate the form on a couple of examples. 
We consider the one-colour case $ A_{\underline{ {\color{blue} s} {\color{blue} s} }}$ where the basis of the cell module 
$ \Delta_{\underline{ {\color{blue} s} {\color{blue} s} }}(1) $ is 
\begin{equation}
   \begin{aligned} 
\begin{tikzpicture}[scale=0.2]
\draw[blue, very thick] (0,0)-- (0, 2);
\draw[blue, very thick] (1,0)-- (1, 2);
\node[circle,draw,blue,fill,minimum size=1mm, inner sep=0pt]  at (0,2) {};
\node[circle,draw,blue,fill,minimum size=1mm, inner sep=0pt]  at (1,2) {};
\draw[blue, very thick] (4,0)-- (4, 2);
\draw[blue, very thick] (5,0)-- (5, 2);
\draw[blue,very thick](4,2) to[out=90,in=180] (4+0.5, 2+0.5) to[out=0,in=90] (5, 2);
\end{tikzpicture}\, \, .   \end{aligned}
\end{equation}
We then get the following values of the bilinear form by taking the coefficient of the empty diagram
\begin{equation}
  \begin{aligned}  
\begin{tikzpicture}[scale=0.4]
\node at (-1+0.5,0.5) {$ \langle  $};
\draw[blue, very thick] (0,0)-- (0, 1);
\draw[blue, very thick] (0.5,0)-- (0.5, 1);
\node[circle,draw,blue,fill,minimum size=1mm, inner sep=0pt]  at (0,1) {};
\node[circle,draw,blue,fill,minimum size=1mm, inner sep=0pt]  at (0.5,1) {};
\node at (1.2,0) {$ ,  $};
\draw[blue, very thick] (0+2,0)-- (0+2, 1);
\draw[blue, very thick] (0.5+2,0)-- (0.5+2, 1);
\node[circle,draw,blue,fill,minimum size=1mm, inner sep=0pt]  at (0+2,1) {};
\node[circle,draw,blue,fill,minimum size=1mm, inner sep=0pt]  at (0.5+2,1) {};
\node at (-1+4,0.5) {$\, \rangle_{1}   $};
\node at (-1+5+0.3,0.7) {$ = { { \color{blue} \alpha^2} } $};

\node at (-1+0.5+8,0.5) {$ \langle  $};
\draw[blue, very thick] (0+8,0)-- (0+8, 1);
\draw[blue, very thick] (0.5+8,0)-- (0.5+8, 1);
\draw[blue, very thick](0+8, 1) to[out=90,in=180] (0+8+0.25, 1+0.25) to[out=0,in=90] (0.5+8, 1);
\node at (1.2+8,0) {$ ,  $};
\draw[blue, very thick] (0+8+2,0)-- (0+8+2, 1);
\draw[blue, very thick] (0.5+8+2,0)-- (0.5+8+2, 1);
\draw[blue, very thick](0+8+2, 1) to[out=90,in=180] (0+8+0.25+2, 1+0.25) to[out=0,in=90] (0.5+8+2, 1);
\node at (-1+4+8,0.5) {$\,  \rangle_1   $};
\node at (-1+8+0.3+5,0.7) {$ = { { \color{blue} 0} } $};
\node at (-1+0.5+8+8,0.5) {$ \langle  $};
\draw[blue, very thick] (0+8+8,0)-- (0+8+8, 1);
\draw[blue, very thick] (0.5+8+8,0)-- (0.5+8+8, 1);
\draw[blue, very thick](0+8+8, 1) to[out=90,in=180] (0+8+0.25+8, 1+0.25) to[out=0,in=90] (0.5+8+8, 1);
\node at (1.2+8+8,0) {$ ,  $};
\draw[blue, very thick] (0+8+2+8,0)-- (0+8+2+8, 1);
\draw[blue, very thick] (0.5+8+2+8,0)-- (0.5+8+2+8, 1);
\node at (-1+4+8+8,0.5) {$ \, \rangle_1   $};
\node at (-1+8+0.3+5+8,0.7) {$ = { { \color{blue} \alpha.} } $};
\node[circle,draw,blue,fill,minimum size=1mm, inner sep=0pt]  at (0+8+2+8, 1) {};
\node[circle,draw,blue,fill,minimum size=1mm, inner sep=0pt]  at (0.5+8+2+8, 1) {};
\end{tikzpicture}    \end{aligned}
\end{equation}
Thus the determinant $ \det  \langle \cdot, \cdot \rangle_{1} $ of $   \langle \cdot, \cdot \rangle_{1} $ is 
$\color{blue} -\alpha^2 $. 
In particular,  
since $ \color{blue} -\alpha^2 $ is a unit in $ Q $ we get 
that $ L_{\underline{ {\color{blue} s} {\color{blue} s} }}(1)   = \Delta_{\underline{ {\color{blue} s} {\color{blue} s} }}(1) $
which we of course already knew
%{\color{green}{\sout{{this}}}}
from the semisimplicity results of the previous section. 
On the other hand, if we were working over a ground field {\color{black}{$ \Bbbk $}} in which 
$ {\color{blue} \alpha}  $ is specialized to $ 0 $ then $   \langle \cdot, \cdot \rangle_{1} =0  $ 
and we would have $ 1 \not\in \Lambda_0$.

We next consider the cell module 
$ \Delta_{\underline{ {\color{blue} s} {\color{blue} s} }}({ \color{blue} s} ) $ with diagram basis 
\begin{equation}{\label{lightleavesA1simple_ss1}}
\begin{aligned}
  \begin{tikzpicture}[scale=0.2]
\draw[blue, very thick] (0+4,0)-- (0+4, 2);
\draw[blue, very thick] (1+4,0)-- (1+4, 2);
\draw[blue,very thick](0+4,2) to[out=90,in=180] (0+4+0.5,2.5) to[out=0,in=90] (0+4+1,2);
\draw[blue, very thick] (0+4+0.5,2.5)-- (0+4+0.5,3);
\draw[blue, very thick] (4-4,0)-- (4-4, 2);
\draw[blue, very thick] (5-4,0)-- (5-4, 3);
\node[circle,draw,blue,fill,minimum size=1mm, inner sep=0pt]  at (4-4,2) {};
\end{tikzpicture} \, \, .\end{aligned}
\end{equation}

We then get the following values of the bilinear form 
\begin{equation} \begin{aligned}
\begin{tikzpicture}[scale=0.4]
\node at (-1+0.5,0.5) {$ \langle  $};
\draw[blue, very thick] (0,0)-- (0, 1);
\draw[blue, very thick] (0.5,0)-- (0.5, 1.5);
\node[circle,draw,blue,fill,minimum size=1mm, inner sep=0pt]  at (0,1) {};
%\node[circle,draw,blue,fill,minimum size=1mm, inner sep=0pt]  at (0.5,1) {};
\node at (1.2,0) {$ ,  $};
\draw[blue, very thick] (0+2,0)-- (0+2, 1);
\draw[blue, very thick] (0.5+2,0)-- (0.5+2, 1.5);
\node[circle,draw,blue,fill,minimum size=1mm, inner sep=0pt]  at (0+2,1) {};
%\node[circle,draw,blue,fill,minimum size=1mm, inner sep=0pt]  at (0.5+2,1) {};
\node at (-1+4,0.5) {$\, \, \rangle_{ {\color{blue} s_1} }   $};
\node at (-1+5+0.3,0.7) {$ = { { \color{blue} \alpha} } $};

\node at (-1+0.5+8,0.5) {$ \langle  $};
\draw[blue, very thick] (0+8,0)-- (0+8, 1);
\draw[blue, very thick] (0.5+8,0)-- (0.5+8, 1);
\draw[blue, very thick](0+8, 1) to[out=90,in=180] (0+8+0.25, 1+0.25) to[out=0,in=90] (0.5+8, 1);
\draw[blue, very thick](0+8+0.25, 1+0.25) --(0+8+0.25, 1+0.25+0.25); 
\node at (1.2+8,0) {$ ,  $};
\draw[blue, very thick] (0+8+2,0)-- (0+8+2, 1);
\draw[blue, very thick] (0.5+8+2,0)-- (0.5+8+2, 1);
\draw[blue, very thick](0+8+2, 1) to[out=90,in=180] (0+8+0.25+2, 1+0.25) to[out=0,in=90] (0.5+8+2, 1);
\draw[blue, very thick] (0+8+0.25+2, 1+0.25) -- (0+8+0.25+2, 1+0.25+0.25);
\node at (-1+4+8,0.5) {$\, \,  \rangle_{\color{blue} s_1}   $};
\node at (-1+8+0.3+5,0.7) {$ = { { \color{blue} 0} } $};
\node at (-1+0.5+8+8,0.5) {$ \langle  $};
\draw[blue, very thick] (0+8+8,0)-- (0+8+8, 1);
\draw[blue, very thick] (0.5+8+8,0)-- (0.5+8+8, 1);
\draw[blue, very thick](0+8+8, 1) to[out=90,in=180] (0+8+0.25+8, 1+0.25) to[out=0,in=90] (0.5+8+8, 1);
\draw[blue, very thick] (0+8+0.25+8, 1+0.25) -- (0+8+0.25+8, 1+0.25+0.25); 
\node at (1.2+8+8,0) {$ ,  $};
\draw[blue, very thick] (0+8+2+8,0)-- (0+8+2+8, 1);
\draw[blue, very thick] (0.5+8+2+8,0)-- (0.5+8+2+8, 1 +0.5);
\node at (-1+4+8+8,0.5) {$ \, \,  \rangle_{ {\color{blue} s_1}}  $};
\node at (-1+8+0.3+5+8,0.7) {$ = { {  1.} } $};
\node[circle,draw,blue,fill,minimum size=1mm, inner sep=0pt]  at (0+8+2+8, 1) {};
%\node[circle,draw,blue,fill,minimum size=1mm, inner sep=0pt]  at (0.5+8+2+8, 1) {};
\end{tikzpicture} \end{aligned}
\end{equation}
Thus the determinant 
$ \det  \langle \cdot, \cdot \rangle_{ \color{blue} s} $ of $   \langle \cdot, \cdot \rangle_{ \color{blue} s} $ is 1
and $ \Delta_{\underline{ {\color{blue} s} {\color{blue} s} }}({ \color{blue} s} ) $ is irreducible, even in specializations.

\medskip
The billinear forms we are here considering coincide with the intersection forms
from \cite{JW} and \cite{Wi}, see also \cite{EW} in the bimodule 
setting. In these references, the forms are further decomposed using the grading. Our diagonalization methods, 
to be explained next, 
only work in the ungraded case.

\medskip
In complete generality the form $  \langle \cdot, \cdot \rangle_{\lambda} $ is
very difficult to control, but 
for a cellular algebra with a separating family $ {\bf L} $ of JM-elements, there is at least a
diagonalization strategy which can be used to calculate the determinant 
$ \det  \langle \cdot, \cdot \rangle_{\lambda} $ of $  \langle \cdot, \cdot \rangle_{\lambda} $.
Let us briefly explain this strategy, {\color{black}{using section 3 of}}
\cite{Ma} as a reference. Let
{\color{black}{$ A $}}
be the cellular algebra in {\color{black}{question}}, defined over {\color{black}{a}} field
{\color{black}{$ \Bbbk $}}, and admitting 
a separating family ${ \bf L} \, {\color{black}{=\{  L_1, L_2, \ldots, L_k \} }} $ of JM-elements. Let 
the cell modules be {\color{black}{$ \Delta( \lambda) $}} for $ \lambda \in \Lambda$.
{\color{black}{For $ i = 1, \ldots, k $ set $ \mathcal{C}(i) := \{ c_t(i) \mid t \in {\rm Tab}(\Lambda) \}  $ where
$ {\rm Tab}(\Lambda) $ is as in (\ref{tabb})}}. 
Suppose that $ s \in {\rm Tab}(\lambda)$ and set 
\begin{equation}
F_s := \prod_{i{\color{black}{=1, \ldots, k}}} \prod_{ \substack{{\color{black}{ c \in \mathcal{C}(i) }}\\ c \neq c_s(i)  }} \frac{ L_i -c }{ c_s(i) -c } \in {\color{black}{A}}, 
\end{equation}
{\color{black}{see Definition 3.1 of \cite{Ma}}}.
Then we have that {\color{black}{$ F_s \neq 0$}} and 
\begin{equation}\label{YORTF}
  \begin{array}{lcr}
 L_i F_s = c_s(i)  F_s,   & 
    \sum_{ s \in {\rm Tab}(\Lambda) } F_s = 1  &\mbox{and} \,\,\,\, F_s F_t = \delta_{st} F_s, 
\end{array}
\end{equation}
where $ \delta_{st} $ is the Kronecker delta. Define the {\it seminormal basis} elements as 
$ f_{st} := F_s C_{st}^{\lambda} F_t \in {\color{black}{A}}$ for $ s, t \in {\rm Tab}(\lambda)$. 
For each $ \lambda \in \Lambda $ fix some $ t_0 \in {\rm Tab}(\lambda)$ and define 
$ C_{s}^{\lambda} = C_{st_0}^{\lambda} $ and $ f_{s}^{\lambda} =  f_{st_0} $. Then 
$ \{ f_{s}^{\lambda} | s \in {\rm Tab}(\lambda) \} $ is the seminormal basis for $ \color{black}{\Delta(\lambda)}$.
Defining $ \gamma_s \in \color{black}{\Bbbk} $ by 
\begin{equation}{\label{gamma}} \gamma_s :=
%\mbox{{\color{green}{\sout{$\langle f_{s}^{\lambda},  f_{s}^{\lambda}  \rangle_{\lambda} = $}}}}
\langle f_{s}^{\lambda},  C_{s}^{\lambda}  \rangle_{\lambda},
\end{equation}
we arrive at 
the following {\color{black}{orthogonality statement}
\begin{equation}
  \langle f_{s}^{\lambda},  f_{t}^{\lambda}  \rangle_{\lambda} = \left\{ \begin{array}{ll} \gamma_s & \mbox{if  } s = t \\
    0 & \mbox{otherwise}  \end{array} \right.
\end{equation}
and} determinant formula 
\begin{equation}
\det  \langle \cdot, \cdot \rangle_{\lambda} = \prod_{ s \in {\rm Tab}(\lambda)} \gamma_s.
\end{equation}

Let us now apply this theory to our cellular algebra $ A_{ \underline{w}}$, defined in terms of 
$ \underline{w} = s_{i_1}  s_{i_2} \cdots s_{i_k}  \in {\bf{exp}}_s$. 
As we have already seen, {\color{black}{$ A_{\underline{w}}$}} is a cellular algebra over $R$. 
Let $ A_{\underline{w}, Q} :=  A_{\underline{w}} \otimes_R Q $. 
In this particular case there is a
%{\color{green}{\sout{{very}}}}
direct diagrammatic way of
constructing the eigenspace idempotents $F_{\bf e} \in A_{ \underline{w}, Q} $ which can be easily explained 
using the previous sections and already appears implicitly in \cite{EW1}. 
Indeed, each $ \e \in {\rm Tab}(\Lambda)= \{0,1\}^k $ gives rise to the following diagram 
$ F_{\bf e}^{diag} $ in $ Kar({\cal D}_Q) $
\begin{equation}
  \begin{aligned}  
  \begin{tikzpicture}[xscale=0.4, yscale=0.6]
\node at (0-0.5,1) {\small $ \frac{1}{\alpha_{i_1}} $};
\draw[black, very thick](0,0) -- (0,0.5);
\draw[black, very thick](0,1.5) -- (0,2);
\draw[black, very thick, densely dotted](0,0.5) -- (0,1.5);
\node at (0+1,1) {\small $ \frac{1}{\alpha_{i_2}} $};
\draw[black, very thick](0+1.5,0) -- (0+1.5,0.5);
\draw[black, very thick](0+1.5,1.5) -- (0+1.5,2);
\node[circle,draw,black,fill,minimum size=1mm, inner sep=0pt]  at (0+1.5,0.5) {};
\node[circle,draw,black,fill,minimum size=1mm, inner sep=0pt]  at (0+1.5,1.5) {};
\node at (0,-0.6) {\small $ e_{1 }$}; 
\node at (1.5,-0.6) {\small $ e_{2} $};
\node at (4.5+0.5,1) {\small $ \color{blue} \frac{1}{\alpha_{i_j}} $};
\draw[blue, very thick](5+0.5,0) -- (5+0.5,0.5);
\draw[blue, very thick](5+0.5,1.5) -- (5+0.5,2);
\node at (5+0.5,-0.6) {\small $ e_{j} $}; 
 \node[circle,draw,blue,fill,minimum size=1mm, inner sep=0pt]  at (5+0.5,0.5) {};
 \node[circle,draw,blue,fill,minimum size=1mm, inner sep=0pt]  at (5+0.5,1.5) {};
 \node[circle,draw,black,fill,minimum size=0.5mm, inner sep=0pt]  at (2.5,1) {};
 \node[circle,draw,black,fill,minimum size=0.5mm, inner sep=0pt]  at (3,1) {};
 \node[circle,draw,black,fill,minimum size=0.5mm, inner sep=0pt]  at (3.5,1) {};
 \node[circle,draw,black,fill,minimum size=0.5mm, inner sep=0pt]  at (6.5,1) {};
 \node[circle,draw,black,fill,minimum size=0.5mm, inner sep=0pt]  at (7,1) {};
 \node[circle,draw,black,fill,minimum size=0.5mm, inner sep=0pt]  at (7.5,1) {};
\draw[black, very thick](1+0+0+9,0) -- (1+0+0+9,0.5);
\draw[black, very thick](1+0+0+9,1.5) -- (1+0+0+9,2);
\draw[black, very thick, densely dotted](1+0+0+9,0.5) -- (1+0+0+9,1.5);
\node at (0-1+0.5+9+0.5,1) {\tiny $ \frac{1}{\alpha_{i_{k-1}}} $};
\node at (1+0+1+9,1) {\small $ \frac{1}{\alpha_{i_k}} $};
\draw[black, very thick](1+0+1.5+9,0) -- (1+0+1.5+9,0.5);
\draw[black, very thick](1+0+1.5+9,1.5) -- (1+0+1.5+9,2);
\node[circle,draw,black,fill,minimum size=1mm, inner sep=0pt]  at (1+0+0+1.5+9,0.5) {};
\node[circle,draw,black,fill,minimum size=1mm, inner sep=0pt]  at (1+0+0+1.5+9,1.5) {};
\node at (-3,1) {$ F_{\bf e}^{diag} :=  $};
\node at (1+9,-0.6) {\small $ e_{k-1} $}; 
\node at (1+0+0+10+0.5,-0.6) {\small $ e_{k} $};
  \end{tikzpicture}
  \end{aligned}  
\end{equation}
\noindent
where $ e_{j} = 1 $ corresponds to
{\color{black}{\begin{tikzpicture}[scale=0.4]
\node at (0-0.7,0) {\small $ \frac{1}{\alpha_{i_j}} $};
\draw[black, very thick](0,-1) -- (0,0.5-1);
\draw[black, very thick](0,1.5-1) -- (0,2-1);
\draw[black, very thick, densely dotted](0,0.5-1) -- (0,1.5-1);
\end{tikzpicture}}
and $ e_{j} = 0 $ corresponds to 
\begin{tikzpicture}[scale=0.4]
\node at (4.5+0.3,1) {\small $ \color{black} \frac{1}{\alpha_{i_j}} $};
\draw[black, very thick](5+0.5,0) -- (5+0.5,0.5);
\draw[black, very thick](5+0.5,1.5) -- (5+0.5,2);
 \node[circle,draw,black,fill,minimum size=1mm, inner sep=0pt]  at (5+0.5,0.5) {};
 \node[circle,draw,black,fill,minimum size=1mm, inner sep=0pt]  at (5+0.5,1.5) {};
\end{tikzpicture}
at the $ j$'th position}, and the colours of the lines are adjusted to $ \underline{w}$.
For example, in the above picture $e_{1}= 1, e_{2} = 0 $ and $ s_{j} $ is blue.
%{\color{green}{\sout{{Using ({\ref{onecolour 2}}) we get that $ \sum_{\bf e}  F_{\bf e}^{diag} = 1 $ and therefore we have }}}}
{\color{black}{Let us show the equality $  F_{\bf e}^{diag}= F_{\bf e} $:
 \begin{lem}
For all $  \e \in {\rm Tab}(\Lambda) $ we have that 
\begin{equation}\label{needtoshow} F_{\bf e}^{diag}= F_{\bf e}. \end{equation}
 \end{lem}     
\begin{dem}
We first check that $ F_\e^{diag} \neq 0 $.
Multiplying $  F_\e^{diag} $ at the $ j$'th position on top and bottom
with 
\begin{tikzpicture}[scale=0.4]
\draw[black, very thick, densely dotted](0,0) -- (0,0.7);
\end{tikzpicture}
if $ e_j = 1 $, and on top with 
\begin{tikzpicture}[scale=0.4]
\draw[black, very thick](5+0.5,0) -- (5+0.5,0.5);
%\draw[black, very thick](5+0.5,1.5) -- (5+0.5,2);
 \node[circle,draw,black,fill,minimum size=1mm, inner sep=0pt]  at (5+0.5,0.5) {};
 %\node[circle,draw,black,fill,minimum size=1mm, inner sep=0pt]  at (5+0.5,1.5) {};
\end{tikzpicture}
and on bottom with
\begin{tikzpicture}[scale=0.4]
%\draw[black, very thick](5+0.5,0) -- (5+0.5,0.5);
\draw[black, very thick](5+0.5,1.5) -- (5+0.5,2);
% \node[circle,draw,black,fill,minimum size=1mm, inner sep=0pt]  at (5+0.5,0.5) {};
\node[circle,draw,black,fill,minimum size=1mm, inner sep=0pt]  at (5+0.5,1.5) {};
\end{tikzpicture} if $ e_j = 0 $, and reducing the result
via relations ({\ref{picture4}}), ({\ref{picture16}}) and ({\ref{onecolour 3}}), 
we arrive at $ \lambda \cdot 1_{\underline{w}^{\e}}$ where $ \lambda \in Q^{\times} $
and where $ 1_{\underline{w}^{\e}} $ is the identity element of $ {\rm End}_{{\cal D}^{std}} (\underline{w}^{\e})$. 
But by ({\ref{standard}}) and the comments following ({\ref{standard}}), we have that 
$ {\rm End}_{{\cal D}^{std}} (\underline{w}^{\e}) \neq 0 $, 
and so also $ 1_{\underline{w}^{\e}} \neq 0 $. Hence 
$ F_{\bf e}^{diag} \neq 0 $, as claimed.

\medskip

Let $ {\mathcal L }_Q $ be the subalgebra of $ A_{ \underline{w}, Q} $ generated by  
$ \{ L_1,L_2 \ldots, L_k \} $. Then it follows from the relations ({\ref{onecolour 2}})
and ({\ref{picture16}}) 
that $ F_{\bf e}^{diag} \in {\mathcal L }_Q $. 
On the other hand, from (\ref{YORTF}) we have that $ L_i = \sum_{\e \in {\rm Tab}(\Lambda)} c_{\e}(i) F_{\e}$
and so $ {\mathcal L }_Q $ can also be described as the subalgebra of $ A_{ \underline{w}, Q} $
generated by the $  F_{\e}$'s. From this we obtain an expansion
\begin{equation}\label{expansionYSFX}
F_{\bf e}^{diag} = \sum_{\f \in {\rm Tab}(\Lambda) } \lambda_{\f}  F_{\f}, \, \lambda_{\f} \in Q. 
\end{equation}  
We next observe that the formulas in
(\ref{YORTF}) also hold when $ F_{\bf e}$ is replaced by $ F_{\bf e}^{diag}$. 
Indeed, the first formula in (\ref{YORTF}), involving the content function, is a consequence of the relations
({\ref{picture4}}), ({\ref{picture16}}) and ({\ref{onecolour 1}}), whereas the two other formulas 
follow from the relations ({\ref{onecolour 1}}) and ({\ref{onecolour 2}}).

Now by Proposition \ref{isanapp} and the Remark following it, there is for each
$ \f \neq \e $ an $ i_\f \in \{1, 2, \ldots, k\} $ such that $ c_\e( i_\f ) \neq  c_\f( i_\f ) $. Hence multiplying each side of
(\ref{expansionYSFX}) by $ \prod_{\f \neq \e } (L_i - c_\f( i_\f ) )$
we get 
\begin{equation}\label{expansionYSF}
\prod_{\f \neq \e } (c_\e( i_\f ) - c_\f( i_\f ) ) F_{\bf e}^{diag} = \prod_{\f \neq \e } (c_\e( i_\f ) - c_\f( i_\f ) ) \lambda_{\e} F_{\bf e}
\, \, \, \Longleftrightarrow \, \, \, F_{\bf e}^{diag} = \lambda_{\e} F_{\bf e}.
\end{equation}  
But both $ F_{\bf e}^{diag} $  and $ F_{\bf e} $ are nonzero idempotents and so
$ \lambda_{\e} = 1$, which shows 
(\ref{needtoshow}).
\end{dem}  

}}

\medskip

We may now use $F_{\e}^{diag} $ to calculate $ \gamma_\e \in Q $. 
Recall that we have fixed $ \underline{w} = s_{j_1} s_{j_2} \cdots s_{j_{k}} $
and $ \e \in {\rm Tab}(y)$. Let $ i \in \{ 1, \ldots, k \} $ and define
$ w^{<i} := s_{j_1}^{e_1} s_{j_2}^{e_2} \ldots s_{j_{i-1}}^{e_{i-1}}$. Let $M$ be the symbol of $ \e $ at the $i$'th position 
and define for $ i \in \{1, \ldots, k \} $ the element $\epsilon_{\e}^i \in Q $ by the formula
\begin{equation}{\label{gammai}}
\epsilon_{\e}^i = \left\{ \begin{array}{ll}  
w^{<i} \alpha_{j_i}  &\mbox{if } T = U0  \\
(w^{<i} \alpha_{j_i})^{-1}  &\mbox{if } T = D0  \\
1  &\mbox{if } T = U1 \\
-1  &\mbox{if } T =  D1.  \\
\end{array} \right.  
\end{equation}

We now claim that $ \epsilon_{\e}^i = \gamma_{\e}^i $. 
%{\color{green}{\sout{{Believing }}}}
{\color{black}{Assuming}}
this claim {\color{black}{is true}} we would obtain the promised determinant formula. 
It expresses the determinant $ \det  \langle \cdot, \cdot \rangle_{y} $ 
in terms of a purely Coxeter group combinatorial calculation, and in particular avoids  
any diagram calculations. Here is the precise formula and the proof of the claim: 
\begin{pro}{\label{determinantformula1}}
In the above notation we have that $ \epsilon_{\e}= \gamma_{\e} $ and $ \gamma_{\e} = \prod_i \epsilon_{\e}^i $ and so 
\begin{equation}\det  \langle \cdot, \cdot \rangle_{y} = \prod_{ \e \in {\rm Tab}(y), i  \in \{ 1, \ldots, k \} } \epsilon_\e^i.
\end{equation}  
\end{pro}
\begin{dem}
Let us warm up by considering the 
$ A_2$-example in which $ \e = (1,0,0) $, 
$ S := \{ {\color{blue} \alpha} , {\color{ red} \beta } \} $, 
$ m_{ {\color{blue} \alpha}  {\color{ red} \beta }} = 3 $ and 
$ \underline{w} := {\color{blue} s_{\alpha}} {\color{red} s_{\beta}} {\color{blue} s_{\alpha}} $. 
The corresponding light leaves diagram is given by $ (U1,U0,D0) $, it is 
\begin{equation}
 \begin{aligned} 
\begin{tikzpicture}[xscale=0.35, yscale=0.2]
\draw[blue, very thick] (0+0,0)-- (0+0,2);
\draw[red, very thick] (1+0,0)-- (1+0,2);
\node[circle,draw,red,fill,minimum size=1mm, inner sep=0pt]  at (1+0,2) {};
\draw[blue,very thick](0+0,2) to[out=90,in=180] (1+0, 3) to[out=0,in=90] (2+0,2);
\draw[blue, very thick] (2+0,0)-- (2+0,2);
\draw[blue, very thick] (1,3)-- (1,4);
\end{tikzpicture} \, \, .  \end{aligned} 
\end{equation}
By the definitions, in order to calculate $ \gamma_{\e} $
we should
%{\color{green}{\sout{{then }}}}
expand the following diagram in terms of the light leaves basis
and take the coefficient of the identity diagram
\begin{equation}  \begin{aligned} 
\begin{tikzpicture}[xscale=0.3, yscale=0.2]
\draw[blue, very thick] (0+0,0)-- (0+0,2);
\draw[red, very thick] (1+0,0)-- (1+0,2);
\node[circle,draw,red,fill,minimum size=1mm, inner sep=0pt]  at (1+0,2) {};
\draw[blue,very thick](0+0,2) to[out=90,in=180] (1+0, 3) to[out=0,in=90] (2+0,2);
\draw[blue, very thick] (2+0,0)-- (2+0,2);
\draw[blue, very thick] (1,3)-- (1,4);
\node[circle,draw,red,fill,minimum size=1mm, inner sep=0pt]  at (1+0,0) {};
\node[circle,draw,blue,fill,minimum size=1mm, inner sep=0pt]  at (2+0,0) {};
\node at (-0.5,-0.5) { $   \color{blue} \tiny \frac{1}{\alpha} $};
\node at (0.5,-0.5) { $   \color{red} \tiny \frac{1}{\beta} $};
\node at (1.5,-0.5) { $   \color{blue} \tiny \frac{1}{\alpha} $};
\draw[blue, very thick, densely dotted](0,0) -- (0,-1-1);
\draw[blue, very thick] (0+0,-1-1)-- (0+0,-3-1);
\draw[red, very thick] (1+0,-1-1)-- (1+0,-3-1);
\node[circle,draw,red,fill,minimum size=1mm, inner sep=0pt]  at (1+0,-3-1) {};
\draw[blue,very thick](0+0,-3-1) to[out=-90,in=180] (1+0, -4-1) to[out=0,in=-90] (2+0,-3-1);
\draw[blue, very thick] (2+0,-1-1)-- (2+0,-3-1);
\draw[blue, very thick] (1,-4-1)-- (1,-5-1);
\node[circle,draw,red,fill,minimum size=1mm, inner sep=0pt]  at (1+0,-3+1) {};
\node[circle,draw,blue,fill,minimum size=1mm, inner sep=0pt]  at (2+0,-3+1) {};
\end{tikzpicture}.   \end{aligned} 
\end{equation}
But applying the relations for $ Kar({\cal D}_Q) $ we get that 
\begin{equation}{\label{a2-example}}
 \begin{aligned}   
\begin{tikzpicture}[xscale=0.3, yscale=0.2]
\draw[blue, very thick] (0+0,0)-- (0+0,2);
\draw[red, very thick] (1+0,0)-- (1+0,2);
\node[circle,draw,red,fill,minimum size=1mm, inner sep=0pt]  at (1+0,2) {};
\draw[blue,very thick](0+0,2) to[out=90,in=180] (1+0, 3) to[out=0,in=90] (2+0,2);
\draw[blue, very thick] (2+0,0)-- (2+0,2);
\draw[blue, very thick] (1,3)-- (1,4);
\node[circle,draw,red,fill,minimum size=1mm, inner sep=0pt]  at (1+0,0) {};
\node[circle,draw,blue,fill,minimum size=1mm, inner sep=0pt]  at (2+0,0) {};
\node at (-0.5,-0.5) { $   \color{blue} \tiny \frac{1}{\alpha} $};
\node at (0.5,-0.5) { $   \color{red} \tiny \frac{1}{\beta} $};
\node at (1.5,-0.5) { $   \color{blue} \tiny \frac{1}{\alpha} $};
\draw[blue, very thick, densely dotted](0,0) -- (0,-1-1);
\draw[blue, very thick] (0+0,-1-1)-- (0+0,-3-1);
\draw[red, very thick] (1+0,-1-1)-- (1+0,-3-1);
\node[circle,draw,red,fill,minimum size=1mm, inner sep=0pt]  at (1+0,-3-1) {};
\draw[blue,very thick](0+0,-3-1) to[out=-90,in=180] (1+0, -4-1) to[out=0,in=-90] (2+0,-3-1);
\draw[blue, very thick] (2+0,-1-1)-- (2+0,-3-1);
\draw[blue, very thick] (1,-4-1)-- (1,-5-1);
\node[circle,draw,red,fill,minimum size=1mm, inner sep=0pt]  at (1+0,-3+1) {};
\node[circle,draw,blue,fill,minimum size=1mm, inner sep=0pt]  at (2+0,-3+1) {};

\node at (4,-1) { $  = $};

\draw[densely dotted, blue, very thick] (7,-5) -- (7,3);
\draw[blue, very thick] (7,-5) -- (7,-6);
\draw[ blue, very thick] (7,4) -- (7,3);
\node at (6.5,-0.5) { $   \color{blue} \tiny \frac{1}{\alpha} $};
\node at (8,-0.5) { $   \color{red} \tiny {\beta} $};
\node at (9,-0.5) { $   \color{blue} \tiny \frac{1}{\alpha} $};

\node at (-1+12,-1) { $  = $};

\draw[ blue, very thick] (-1+7+7,-5) -- (-1+7+7,3);
\draw[blue, very thick] (-1+7+7,-5) -- (-1+7+7,-6);
\draw[ blue, very thick] (-1+7+7,4) -- (-1+7+7,3);
\node at (-1+7+8,-0.5) { $   \color{red} \tiny {\beta} $};
\node at (-1+7+9,-0.5) { $   \color{blue} \tiny \frac{1}{\alpha} $};

\node at (8+-3+7+8-3,-0.5) { $   - $};

\node at (8+-3+7+8,-0.5) { $   \color{blue} \tiny \frac{1}{\alpha} $};
\node[circle,draw,blue,fill,minimum size=1mm, inner sep=0pt]  at (8+-1+7+7,-5) {};
\node[circle,draw,blue,fill,minimum size=1mm, inner sep=0pt]  at (8+-1+7+7,3) {};
%\draw[densely dotted, blue, very thick] (8+-1+7+7,-5) -- (8+-1+7+7,3);
\draw[blue, very thick] (8+-1+7+7,-5) -- (8+-1+7+7,-6);
\draw[ blue, very thick] (8+-1+7+7,4) -- (8+-1+7+7,3);
\node at (8+-1+7+8,-0.5) { $   \color{red} \tiny {\beta} $};
\node at (8+-1+7+9,-0.5) { $   \color{blue} \tiny \frac{1}{\alpha} $};
\end{tikzpicture} \, \,  \end{aligned} 
\end{equation}
where we used ({\ref{onecolour 2}}) for the second equality. 
Applying the polynomial relation ({\ref{picture5}}) once to this, we find that $ \gamma_\e$ is equal to 
$ {\color{blue} \tiny s} ({\color{red} \tiny \beta}  
{\color{blue} \tiny \frac{1}{\alpha}}) $, as claimed.

\medskip
Let us now turn to the general case. 
Let us consider an $ \e \in \tab(y) $ and let $ e_i $ be the $ i$'th coordinate of $ \e$. If the corresponding symbol is $ M=U0$, then the situation is exactly as for the $e_2 $-term of 
({\ref{a2-example}}), and the contribution to $ \gamma_\e$ is $ w^{<i} \alpha_{i} $ as claimed. 
If the corresponding symbol is $ M=D0$, then the situation is exactly as for the $e_3 $-term of 
({\ref{a2-example}}), adding a trivalent vertex, and the contribution to $ \gamma_\e$ is $ (w^{<i} \alpha_{i})^{-1} $ as claimed.

If $ T=U1$ then there are two possibilities. Either the situation is as for 
$e_1 $ in ({\ref{a2-example}}) where the dashed line extends to the top and bottom of the diagram and 
so the contribution to $\gamma_{\e} $ will be $1 $, as claimed. Or alternatively, the dashed line turns around and eventually runs into a later $ D1$. This situation is not 
represented in ({\ref{a2-example}}) but appears 
in the example $ \underline{w} := {\color{blue} s}  {\color{blue} s} $,  $ \e = (1,1) $ with 
symbols $ (U1, D1) $. The light leaves morphism is in this case
\begin{tikzpicture}[scale=0.3]
\draw[blue, very thick] (4,0)-- (4, 1);
\draw[blue, very thick] (5,0)-- (5, 1);
\draw[blue,very thick](4,1) to[out=90,in=180] (4+0.5, 1+0.5) to[out=0,in=90] (5, 1);
\end{tikzpicture}
and so we should take the coefficient of the empty diagram in
\begin{equation}  \begin{aligned} 
\begin{tikzpicture}[xscale=0.4, yscale=0.3]
\draw[blue, very thick] (4,0)-- (4, 1);
\draw[blue, very thick] (5,0)-- (5, 1);
\draw[blue,very thick](4,1) to[out=90,in=180] (4+0.5, 1+0.5) to[out=0,in=90] (5, 1);

\draw[densely dotted, blue, very thick] (4,0)-- (4, -2);
\draw[densely dotted, blue, very thick] (5,0)-- (5, -2);

\draw[blue, very thick] (4,-1)-- (4, -2);
\draw[blue, very thick] (5,-1)-- (5, -2);
\draw[blue,very thick](4,-2) to[out=-90,in=180] (4+0.5, -2-0.5) to[out=0,in=-90] (5,-2);

\node at (4-0.5,0) { $   \color{blue} \tiny \frac{1}{\alpha} $};
\node at (4-0.5+1,0) { $   \color{blue} \tiny \frac{1}{\alpha} $};
\end{tikzpicture} \, \, .  \end{aligned} 
\end{equation}
By 
({\ref{picture16}}), ({\ref{picture17}}), ({\ref{onecolour 5}}) and ({\ref{onecolour 6}}) the value of this is $-1$, 
which is therefore the combined contribution of the $ U1 $ and $ D1 $ symbols in this case, as claimed.
The general case is done the same way.
\end{dem}

\medskip
\noindent
{\bf Remark}. Note that $ \gamma_\e^i \in Q $, but even so $\det  \langle \cdot, \cdot \rangle_{y} \in R $, of course. 

\medskip \noindent
{\bf Example}.
Let us take $ \underline{w} := \color{blue} sss $. In (\ref{lightleavesA1}) we have given the basis 
for $ \Delta(1) $ and in (\ref{lightleavesA1s}) the basis for $ \Delta(s) $. But using the symbols given in ({\ref{symbols}}) we get easily via the Proposition, 
without drawing diagrams, 
that $ \det  \langle \cdot, \cdot \rangle_{1} = {\color{blue}\alpha}^4 $ and 
$ \det  \langle \cdot, \cdot \rangle_{\color{blue}s} = {\color{black}1} $.

\medskip
\noindent
{\bf Remark}. It follows from the Proposition that $ \det  \langle \cdot, \cdot \rangle_{y} $ is always a product 
of roots. 

\medskip
\noindent
    {\bf Remark}.
    Note that by general principles the zero'th homogeneous part $ A_{\underline{w}}^0 $
    of $ A_{\underline{w}} $ is
      also a cellular algebra. 
    Unfortunately, as already mentioned in the introduction, 
    our formula is not compatible with the grading on $ A_{\underline{w}} $, and in particular
    is does not give rise to a determinant formula for the cell modules for
$ A_{\underline{w}}^0 $. 
      The reason is that the JM-elements $ \{ L_i \} $
      are homogeneous of degree 2 and hence cannot be used to diagonalize
      the cell modules for $ A_{\underline{w}}^0 $.

\section{A Shapovalov-like expression for $\det  \langle \cdot, \cdot \rangle_{y}$}
The usefulness of determinant expressions in representation theory depends on the property that they
can be rewritten in terms 
of characters or dimensions of other cell modules. This rewriting was first achieved by
Shapovalov and  Jantzen for Verma modules for complex 
semisimple Lie algebras, see \cite{Sh} and \cite{Ja}.  Later, similar 
rewritings were found for Weyl modules for algebraic groups by Jantzen and Andersen,  for Specht modules 
for the symmetric groups by Schaper, James, Murphy and so on, see for example \cite{A}, \cite{A1}, \cite{JM}.
This section is devoted to a rewriting of 
$ \det  \langle \cdot, \cdot \rangle_{y } $ in this spirit. A major difference between our case and the above mentioned cases is 
that we do not have a closed formula for the dimension of the cell module.

\medskip
Our argument will be an induction on the length of $ \underline{w}$ using Proposition {\ref{determinantformula1}}. 
In this sense, it is closest to the argument given by James and Murphy for the symmetric group, see
\cite{JM}.
Just like in \cite{JM} we also need
a branching rule for $ \Delta_{\underline{w}}( y) $. For $ A_{\underline{w}}$ this branching rule was
discovered by D. Plaza in \cite{P} 
in the bimodule setting, but, as we shall shortly see, it carries over to the diagrammatical
setting.
%On the other hand, 
%D. Plaza worked in the algebraic setting of bimodules
%so let us here translate his results into the diagrammatical language of $ \cal D$. 

\medskip
Recall that we have fixed the expression $ \underline{w} = s_{j_1}  \cdots s_{j_k} $. Set 
$ {\color{blue}\alpha}:= \alpha_{j_k} $ and
$ {\color{blue}s}:= s_{j_k} $, where 
$ j_k $ corresponds to the color blue, say. For any $ x \in W $ we define
$ {x}^{\prime} := x\color{blue}s$. 
There is a straightforward embedding $ A_{\underline{w}^{\prime}} \subseteq A_{\underline{w}} $
defined diagrammatically by adding a through blue line to the right of a diagram in $A_{\underline{w}^{\prime}}$, 
as illustrated below
\begin{equation}   \begin{aligned} 
\begin{tikzpicture}[xscale=0.3, yscale=0.3]
\draw [very thick, dotted] (0,0) rectangle (4,4);
\draw [very thick] (1,1) .. controls (1.5,2) and (2.5,2) .. (3,3);
\draw [very thick, ] (1,3) .. controls (1.5,2) and (2.5,2) .. (3,1);
\node at (6,2) {$  \mapsto $};
\draw [very thick,dotted] (0+8,0) rectangle (4+8,4);
\draw [very thick] (1+8,1) .. controls (1.5+8,2) and (2.5+8,2) .. (3+8,3);
\draw [very thick] (1+8,3) .. controls (1.5+8,2) and (2.5+8,2) .. (3+8,1);
\draw [blue, very thick] (5+8,0) -- (5+8,4);
\node at (13.5,0) {$  . $};
\end{tikzpicture}  \end{aligned} 
\end{equation} 
Letting $A_{\underline{w}^{\prime}} \mbox{-mod} $ and $A_{\underline{w}} \mbox{-mod} $
denote the module categories for $A_{\underline{w}^{\prime}}$ and $A_{\underline{w}} $, 
we get an associated restriction functor 
\begin{equation}
A_{\underline{w}} \mbox{-mod}   \rightarrow A_{\underline{w}^{\prime}} \mbox{-mod}, \, \,  M \mapsto \Res M.
\end{equation}
In general, we use from now on 
the convention that $ \Delta_{ \underline{z}}(y ) := 0 $ if $ y \nleq z $.
The following branching rule was found by D. Plaza in the bimodule setting, see \cite{P1}.
We now state and prove it in the diagrammatical setting. 
\begin{teo}{\label{branching}}
Assume that $ y \le w $.

\medskip
\noindent
1)
If $ y^{\prime} > y$ then there is a short exact sequence of $ A_{\underline{w}^{\prime}}$-modules
\begin{equation}{\label{firstseq}}
0 \rightarrow \Delta_{\underline{w}^{\prime}}(y) \rightarrow \Res \Delta_{\underline{w}^{}}(y) \rightarrow  \Delta_{\underline{w}^{\prime}}(y^{\prime}) \rightarrow 0. 
\end{equation}

\medskip
\noindent
2) 
If $ y^{\prime} < y$ then there is a short exact sequence of $ A_{\underline{w}^{\prime}}$-modules
\begin{equation}{\label{secondseq}}
0 \rightarrow \Delta_{\underline{w}^{\prime}}(y^{\prime}) \rightarrow \Res \Delta_{\underline{w}}(y) 
\rightarrow  \Delta_{\underline{w}^{\prime}}(y) \rightarrow 0.
\end{equation}
\end{teo}
\noindent
{\bf Remark}. In \cite{P} the sequences are considered as sequences of graded modules over the graded algebra  
$ A_{\underline{w}^{\prime}}$. In view of the last remark of the previous section we here ignore the grading.

\medskip
\noindent
\begin{dem}
  Let first $M$ be the symbol at the $k$'th position of any $ \e \in \tab(y)$. By definition, if
%{\color{green}{\sout{{we have that}}}} 
$ y^{\prime} >  y $ then $M$ is either $ U0$ or $ D1 $ whereas if
$ y^{\prime} <  y $ then $M$ is either $ U1$ or $ D0 $.

\medskip
Let us first consider the case $ y^{\prime} > y $. The light leaves diagrams $ {\LL}_{ \underline{w}, \bf e} $ satisfying 
$ M= U0$ are exactly those with 
\begin{tikzpicture}[xscale=0.4, yscale=0.6]
 \draw[blue, very thick](5,-0.20) -- (5,0.5);
 \node[circle,draw,blue,fill,minimum size=1mm, inner sep=0pt]  at (5,0.5) {};
\end{tikzpicture} on the right. They define an $ A_{\underline{w}^{\prime}}$-module isomorphic to 
$ \Delta_{\underline{w}^{\prime}}(y) $ where the isomorphism is given by adding 
\begin{tikzpicture}[xscale=0.4, yscale=0.6]
\draw[blue, very thick](5,-0.20) -- (5,0.5);
\node[circle,draw,blue,fill,minimum size=1mm, inner sep=0pt]  at (5,0.5) {};
\end{tikzpicture} to the right.
On the other hand, the light leaves diagrams $ {\LL}_{ \underline{w}, \bf e} $ satisfying 
$ M= D1$ are in bijection with the light leaves basis for $ \Delta_{\underline{w}^{\prime}}(y^{\prime}) $ where 
the bijection can be described by bending up the last line, 
thus transforming a
\begin{tikzpicture}[very thick, scale=0.2]
\draw[blue](0,0) to[out=90,in=180] (1,1) to[out=0,in=90] (2,0);
\end{tikzpicture}
{\color{black}{into a vertical blue line}} \begin{tikzpicture}[very thick, scale=0.2]
\draw[blue](0,0) -- (0,2);
\end{tikzpicture}.
For example, in the $ A1$-case given in ({\ref{lightleavesA1}}) the first two diagrams end in $M=U0$ and the last two 
diagrams end in $ M= D1 $ and these become, after bending up, the diagrams of ({\ref{lightleavesA1simple_ss1}}).

\medskip
The bending up also works well at the module level. Indeed, let us extend it to a linear map $ \pi:
\Delta_{\underline{w}^{}}(y) \rightarrow  \Delta_{\underline{w}^{\prime}}(y^{\prime}) $ 
given diagrammatically via 
\begin{equation}   \begin{aligned} 
\begin{tikzpicture}[xscale=0.2, yscale=0.2]
\node at (-3,2) {$  \pi:  $};

\draw [very thick, dotted] (0,0) rectangle (4,4);
\draw [very thick] (1,1) .. controls (1.5,2) and (2.5,2) .. (3,3);
\draw [very thick, ] (1,3) .. controls (1.5,2) and (2.5,2) .. (3,1);
\draw [blue, very thick] (5,-1) -- (5,2);
\node[circle,draw,blue,fill,minimum size=1mm, inner sep=0pt]  at (5,2) {};
\node at (7,2) {$  \mapsto $};
\node at (9,2.2) {$  0 $};
\draw [very thick] (1,-1) -- (1,0);
\draw [very thick] (2,-1) -- (2,0);
\draw [very thick] (3,-1) -- (3,0);
\draw [very thick] (1,4) -- (1,5);
\draw [very thick] (2,4) -- (2,5);
%\draw [very thick] (3,4) -- (3,5);

\node at (10,0) {$  , $};

\draw [very thick, dotted] (0+12,0) rectangle (4+12,4);
\draw [very thick] (1+12,1) .. controls (1.5+12,2) and (2.5+12,2) .. (3+12,3);
\draw [very thick, ] (1+12,3) .. controls (1.5+12,2) and (2.5+12,2) .. (3+12,1);
\draw [blue, very thick] (5+12,-1) -- (5+12,4);
\node at (7+12,2) {$  \mapsto $};
\draw [very thick] (1+12,-1) -- (1+12,0);
\draw [very thick] (2+12,-1) -- (2+12,0);
\draw [very thick] (3+12,-1) -- (3+12,0);
\draw [very thick] (1+12,4) -- (1+12,5);
\draw [very thick] (2+12,4) -- (2+12,5);
%\draw [blue, very thick] (3+12,4) -- (3+12,5);
\draw[blue, very thick](3+12,4) to[out=90,in=180] (4+12,4+1) to[out=0,in=90] (5+12,4);

\draw [very thick, dotted] (0+12+9,0) rectangle (4+12+9,4);
\draw [very thick] (1+12+9,1) .. controls (1.5+12+9,2) and (2.5+12+9,2) .. (3+12+9,3);
\draw [very thick, ] (1+12+9,3) .. controls (1.5+12+9,2) and (2.5+12+9,2) .. (3+12+9,1);
%\draw [blue, very thick] (5+12+9,0) -- (5+12+9,4);
\draw [very thick] (1+12+9,-1) -- (1+12+9,0);
\draw [very thick] (2+12+9,-1) -- (2+12+9,0);
\draw [very thick] (3+12+9,-1) -- (3+12+9,0);
\draw [very thick] (1+12+9,4) -- (1+12+9,5);
\draw [very thick] (2+12+9,4) -- (2+12+9,5);
%\draw [blue, very thick] (3+12,4) -- (3+12,5);
%\draw[blue, very thick](3+12+9,4) to[out=90,in=180] (4+12+9,4+1) to[out=0,in=90] (5+12+9,4);
\draw[blue, very thick](3+12+9,4) -- (3+12+9,5); 
\end{tikzpicture} \, \, .  \end{aligned} 
\end{equation} 
Then $ \ker \pi = \Delta_{\underline{w}^{\prime}}(y) $, $ \im \pi = \Delta_{\underline{w}^{\prime}}(y) $. Moreover, 
$ \pi $ is an $ A_{\underline{w}^{\prime}}$-linear map as follows from the diagrammatical identity 

\begin{equation}   \begin{aligned} 
\begin{tikzpicture}[xscale=0.2, yscale=0.2]
\draw [very thick, dotted] (0+12,0) rectangle (4+12,4);
\draw [very thick] (1+12,1) .. controls (1.5+12,2) and (2.5+12,2) .. (3+12,3);
\draw [very thick, ] (1+12,3) .. controls (1.5+12,2) and (2.5+12,2) .. (3+12,1);
\draw [blue, very thick] (5+12,-1) -- (5+12,4);
\node at (7+12+3,0) {$  =$};
\draw [very thick] (1+12,-1) -- (1+12,0);
\draw [very thick] (2+12,-1) -- (2+12,0);
\draw [very thick] (3+12,-1) -- (3+12,0);
\draw [very thick] (1+12,4) -- (1+12,5);
\draw [very thick] (2+12,4) -- (2+12,5);
\draw[blue, very thick](3+12,4) to[out=90,in=180] (4+12,4+1) to[out=0,in=90] (5+12,4);

\draw [very thick, dotted] (0+12,-5) rectangle (4+12,-5+4);
\draw [very thick] (1+12,1-5) .. controls (1.5+12,2-5) and (2.5+12,2-5) .. (3+12,3-5);
\draw [very thick, ] (1+12,3-5) .. controls (1.5+12,2-5) and (2.5+12,2-5) .. (3+12,1-5);
\draw [blue, very thick] (5+12,-1-5+1) -- (5+12,4-5+1);
\draw [very thick] (1+12,-1-5) -- (1+12,0-5);
\draw [very thick] (2+12,-1-5) -- (2+12,0-5);
\draw [very thick] (3+12,-1-5) -- (3+12,0-5);
\draw [very thick] (1+12,4-5) -- (1+12,5-5);
\draw [very thick] (2+12,4-5) -- (2+12,5-5);

\draw[blue, very thick](5+12,-1-5+1+1) to[out=-90,in=180](5+12+1,-1-5-1+1)to[out=0,in=-90](5+12+1+1,-1-5+1); 
\draw[blue, very thick](5+12+1+1,-1-5+1) -- (5+12+1+1,-1-5+10+1);

\draw [very thick, dotted] (0+12+13,0) rectangle (4+12+13,4);
\draw [very thick] (1+12+13,1) .. controls (1.5+12+13,2) and (2.5+12+13,2) .. (3+12+13,3);
\draw [very thick, ] (1+12+13,3) .. controls (1.5+12+13,2) and (2.5+12+13,2) .. (3+12+13,1);
%\draw [blue, very thick] (5+12+13,-1) -- (5+12+13,4);
%\node at (7+12+3+13,0) {$  =$};
\draw [very thick] (1+12+13,-1) -- (1+12+13,0);
\draw [very thick] (2+12+13,-1) -- (2+12+13,0);
\draw [very thick] (3+12+13,-1) -- (3+12+13,0);
\draw [very thick] (1+12+13,4) -- (1+12+13,5);
\draw [very thick] (2+12+13,4) -- (2+12+13,5);
\draw[blue, very thick](3+12+13,4) -- (3+12+13,4+1);

 % to[out=90,in=180] (4+12+13,4+1) to[out=0,in=90] (5+12+13,4);

\draw [very thick, dotted] (0+12+13,-5) rectangle (4+12+13,-5+4);
\draw [very thick] (1+12+13,1-5) .. controls (1.5+12+13,2-5) and (2.5+12+13,2-5) .. (3+12+13,3-5);
\draw [very thick, ] (1+12+13,3-5) .. controls (1.5+12+13,2-5) and (2.5+12+13,2-5) .. (3+12+13,1-5);
%\draw [blue, very thick] (5+12+13,-1-5+1) -- (5+12+13,4-5+1);
\draw [very thick] (1+12+13,-1-5) -- (1+12+13,0-5);
\draw [very thick] (2+12+13,-1-5) -- (2+12+13,0-5);
\draw [very thick] (3+12+13,-1-5) -- (3+12+13,0-5);
\draw [very thick] (1+12+13,4-5) -- (1+12+13,5-5);
\draw [very thick] (2+12+13,4-5) -- (2+12+13,5-5);

%\draw[blue, very thick](5+12+13,-1-5+1+1) to[out=-90,in=180](5+12+1+13,-1-5-1+1)to[out=0,in=-90](5+12+1+1+13,-1-5+1); 
%\draw[blue, very thick](5+12+1+1+13,-1-5+1) -- (5+12+1+1+13,-1-5+10+1);
\end{tikzpicture} \, \,   \end{aligned} 
\end{equation} 
where the left hand side represents $ \pi ( a v ) $ for $ a \in  A_{\underline{w}^{\prime}}$ and 
$ v \in \Delta_{\underline{w}^{}}(y) $ and the right hand side $ a \pi (  v ) $.
Upon expanding the left hand side in terms of the light leaves basis for $\Delta_{\underline{w}^{}}(y) $
the terms with a blue dot in the upper right  corner are mapped to zero under $ \pi$
and are also viewed as zero in the right hand side.

\medskip
In the case $ y^{\prime} < y $, the $\e$'s with $M=U1$ correspond to diagrams with a blue line 
\begin{tikzpicture}[very thick, scale=0.2]
\draw[blue](0,0) -- (0,2);
\end{tikzpicture} 
to the right. These are in correspondence with a basis for the
$ A_{\underline{w}^{\prime}}$-module $ \Delta_{\underline{w}^{\prime}}(y^{\prime}) $, where the correspondence is given by 
deleting the blue line. The remaining diagrams, corresponding to $ M=D0$, have a trivalent vertex on the right.
These diagrams are in bijection with the light leaves basis for $ \Delta_{\underline{w}^{\prime}}(y) $, where the bijection is given by transforming 
the trivalent vertex to a line.

\medskip
Also in this case, one checks that the maps work well at the module level.
\end{dem}

\medskip
In order to give  
the Shapovalov-like reformulation of the determinant, we shall assume that $ \mathfrak h$ is
{\color{black}{of}} 
one of 
{\color{black}{the following two}} types.
% We need to 
% For convenience we shall suppose that $ \mathfrak h $ is one of the following two realizations of $ (W,S)$.
{\color{black}{For the first type}} 
$ {\mathfrak h} $ is chosen as the geometric representation of $ (W,S) $, considered in \cite{H1}. Recall that
this $ {\mathfrak h} $ 
is a realization of $ (W,S) $, in the above sense, defined over $ \mathbb R $ and faithful as a $W$-representation. 
A \emph{root} $ \alpha \in {\mathfrak h}^{\ast}  $ is by definition an element of the form $\alpha = w \alpha_s $ for $ s \in S$ and $ w \in W$.
It is said to be \emph{positive}, written $ \alpha > 0 $, if 
$ \alpha = \sum_{ s \in S} \lambda_s \alpha_s$ for $ \lambda_s \ge 0 $ for all $ s$. 
For $ \alpha  $ a root there is a well-defined associated reflection $ s_{\alpha} $ defined by 
$ s_{\alpha} = z s z^{-1} $ 
where $ z\in W $ is chosen such that $ \alpha = z \alpha_s$.

\medskip
{\color{black}{For the second type}} we
shall assume that $ (W,S) $ is the Weyl group of the root system $ \Phi$ with
Cartan matrix $ A = (a_{st})_{s,t \in S}$. We let $  {\mathfrak h}_{\Z} := \oplus_{s\in S}\Z \alpha^{\vee}_s$
be the coroot lattice. Then $ {\mathfrak h}_{\Z} $, together
with $ \{ \alpha^{\vee}_s \mid  s\in S\} \subset {\mathfrak h}_{\Z} \} $
and $ \{ \alpha_s \mid  s\in S\} \subset {\mathfrak h}_{\Z}^{\ast} \} $, defines a realization of $ (W,S) $ over $ \Z$.
Let $ p >3$ be a prime. Then our second {\color{black}{type}} of realization 
is $ {\mathfrak h}:= {\mathfrak h}_{\Z} \otimes_{\Z} \overline{\F}_p$.
This is the \emph{Cartan matrix representation} of $ W$.
In $ {\mathfrak h}_{\Z}$ we have the concept of roots and positive roots just as before.
\medskip

Although other choices of realizations $ \mathfrak  h $ are possible, at present the above ones
are sufficient for the applications we have in mind. 
We remark that both choices of realizations are Soergel realizations in the above sense.
For the first $ \mathfrak h$ this is a consequence of \cite{Li2} and for the second
$ \mathfrak h$ it follows from the Appendix of \cite{Li3}; this is where $ p>3 $ is needed. Thus in both cases we obtain 
equivalences between the diagrammatic and bimodules categories, see \cite{EW1}.

\medskip

Let us now return to the determinant $ \det  \langle \cdot, \cdot \rangle_{\lambda} $. 
We prove the following Theorem, giving
the promised Shapovalov-like reformulation of the determinant.
\begin{teo}{\label{lightleavesJM_det}}
Let $ {\mathfrak h } $ be as above. Then we have that 
\begin{equation}{\label{shapovalov}}
\det  \langle \cdot, \cdot \rangle_{y} = \pm  \prod_{  \ms \beta > 0, s_{\beta} y > y }   \beta^{ \dim \Delta_{\underline{w}}( s_{\beta} y)}.
\end{equation}
\end{teo}
\begin{dem}
The proof is a purely combinatorial transformation 
of the formula from Proposition {\ref{determinantformula1}} into 
the expression given in ({\ref{shapovalov}}), using induction on the length $ k $ of $ \underline{w} $.
The case $ k=1 $ is straightforward to check.
Let us therefore assume that ({\ref{shapovalov}}) 
holds for $ k-1$ and verify it for $ k$.
We do it simultaneously for the two choices for $ \mathfrak h$; for the second choice we 
work with $  {\mathfrak h}_{\Z} $ and deduce the formula ({\ref{shapovalov}}) for $  {\mathfrak h} $
from the $  {\mathfrak h}_{\Z} $-case
by reduction modulo {\color{black}{$ p$.}}

Suppose first that $ y^{\prime} > y $. As mentioned in the proof of the previous Theorem {\ref{branching}}, 
the $\e$'s ending in $ U0$ give rise to a basis for $ \Delta_{\underline{w}}(y) $, whereas those ending in $ D1 $ give rise to 
a basis for $ \Delta_{\underline{w}}(y^{\prime}) $, by restriction to $ \underline{w}^{\prime}$. Thus, by the inductive hypothesis and 
Proposition {\ref{determinantformula1}} we get that 
\begin{equation}{\label{thefirstproduct}}
%\begin{array}{ll}
\det  \langle \cdot, \cdot \rangle_{y}  =  \pm  ({y\color{blue} \alpha})^{ \dim  \Delta_{\underline{w}^{\prime}}(y)}
\prod_{  \ms \beta > 0,  s_{\beta} y > y }   \beta^{ \dim \Delta_{\underline{w}^{\prime}}( s_{\beta} y)} \prod_{  \ms \beta > 0, s_{\beta}  y^{\prime} > y^{\prime} }   
\beta^{ \dim \Delta_{\underline{w}^{\prime}}( s_{\beta} y^{\prime})}
%\end{array}
\end{equation}
Now using
({\ref{firstseq}}) and ({\ref{secondseq}}) we get that 
\begin{equation}{\label{andsowewould}}
\dim \Delta_{\underline{w}}( s_{\beta} y) =
\dim \Delta_{\underline{w}^{\prime}}( s_{\beta} y^{\prime}) + 
\dim \Delta_{\underline{w}^{\prime}}( s_{\beta} y)
\end{equation}
and so we would like to join the two last products of ({\ref{thefirstproduct}}) into one product,
running over $ \beta > 0 $ with
$ s_{\beta} y > y$. As long as $ \beta $ is in 
$R_1:= \{ \beta > 0 | s_{\beta} y > y  \mbox{ iff } s_{\beta} y^{\prime} > y^{\prime} \}$
this is clearly possible. 
Let us therefore analyse the sets 
$$ R_2:= \{ \beta > 0 | s_{\beta} y <  y  \mbox{ and } s_{\beta} y^{\prime} > y^{\prime} \} \mbox{ and }
R_3:= \{ \beta > 0 | s_{\beta} y > y  \mbox{ and } s_{\beta} y^{\prime} < y^{\prime} \}. $$
We first claim that $ R_2 = \emptyset$. Indeed, if $ \beta \in R_2 $, then $ s_{\beta} y <  y $ and so we get that 
$ l( s_{\beta} y) \le l(y) -1  = l( y^{\prime}) -2$, in contradiction with $ s_{\beta} y^{\prime} > y^{\prime} $.

Let us then consider $ \beta \in R_3$. Let $ y= s_{i_1} \cdots s_{i_m} $  be a reduced expression for $ y $.
Then $ y^{\prime} = s_{i_1} \cdots s_{i_m} {\color{blue}s} $ is a reduced expression for $ y^{\prime} $.
By the strong exchange condition, see Theorem 5.8 of \cite{H1},  applied to $ s_{\beta} y^{\prime} < y^{\prime}  $, 
we get that either
$ s_{\beta} y^{\prime} = s_{i_1} \ldots \widehat{s_{i_j}} \cdots s_{i_m} {\color{blue}s} $ for some $ j $, or 
$ s_{\beta} y^{\prime} = s_{i_1} \cdots s_{i_m}  $. The first case is impossible since it would imply that 
$ s_{\beta} y = s_{i_1} \cdots \widehat{s_{i_j}} \cdots s_{i_m}  < y$. The second case corresponds to 
$ s_{\beta} y  {\color{blue}s} = y  $, or equivalently $ s_{\beta}    = y {\color{blue}s} y^{-1}$ and so 
$ R_3 = \{ y {\color{blue}\alpha} \}.$

All in all, using ({\ref{andsowewould}}) we find that the contribution to $ ({\ref{thefirstproduct}}) $ from $ R_1$ is 
\begin{equation} \pm 
\prod_{  \ms \beta \in R_1,  s_{\beta} y > y }   \beta^{ \dim \Delta_{\underline{w}^{\prime}}( s_{\beta} y)} \prod_{  \ms \beta \in R_1, s_{\beta}  y^{\prime} > y^{\prime} }   
\beta^{ \dim \Delta_{\underline{w}^{\prime}}( s_{\beta} y^{\prime})} = \pm
\prod_{  \ms \beta \in R_1 }   
\beta^{ \dim \Delta_{\underline{w}}( s_{\beta}y )}.
\end{equation}
The contribution from $R_3 $ to $ ({\ref{thefirstproduct}}) $ comes from the second factor. 
But since $ s_{\beta} y = y^{\prime} $, its contribution together with the first factor of $ ({\ref{thefirstproduct}}) $ is 
\begin{equation}
({y\color{blue} \alpha})^{ \dim  \Delta_{\underline{w}^{\prime}}(y)}  ({y\color{blue} \alpha})^{ \dim  \Delta_{\underline{w}^{\prime}}(y^{\prime})}
=  ({y\color{blue} \alpha})^{ \dim  \Delta_{\underline{w}}(y)}.
\end{equation}
Finally, since $ R_2 $ does not contribute and since $ R_1 \cup R_3 = \{ \beta > 0 | s_{\beta} y > y \} $ we conclude that our 
formula ({\ref{shapovalov}}) of the Theorem holds. 

The case $ y^{\prime} < y $ is treated the same way.
\end{dem}

\section{Jantzen type filtrations and sum formulas.}
We construct in this section an analogue of the Jantzen filtration and, as an application of Theorem \ref{lightleavesJM_det}, obtain an associated sum formula. Recall that 
the original Jantzen filtration was constructed for Verma modules for 
semisimple complex Lie algebras. 
Although the combinatorial methods used here have little in common with those used in 
the representation theory of Lie algebras, ultimately the sum formulas rely in both cases on determinant expressions.

\medskip
Let us first suppose that $ {\mathfrak h}  $ is the geometric representation of $ (W,S)$, defined over the real numbers $ \mathbb R$.
Let $ {\cal R } $ be the polynomial algebra $ {\cal R }:= { \mathbb R}[x] $.
Then $ {\cal R } $ is a discrete valuation domain with maximal ideal $ (x) $ and $x$-adic valuation
{\color{black}{which}} we denote by 
$ \nu_x(\cdot)  $.
Recall that $ R = \oplus_m S^{m}({\mathfrak h}^{\ast}) $. 
It may be identified with the polynomial algebra over $ \mathbb R $ in $ |S| $ variables
and hence there is a natural algebra homomorphism 
$ \varphi: R \rightarrow \cal R  $ satisfying $ \nu_x( \varphi(\beta) ) = 1 $ for all roots $ \beta $.

\medskip
Let us consider the base change from 
$ R $ to $ \cal R $
\begin{equation}{\label{scalarextension}}
A_{\underline{w},{\cal R}}:= A_{\underline{w}} \otimes_R \cal R.
\end{equation}
By the general cellular algebra theory, we then know that $ A_{\underline{w},{\cal R}} $ 
is also a cellular algebra, with cell modules $ \Delta_{\underline{w}, \cal R}(y) := 
\Delta_{\underline{w}}(y) \otimes_R \cal R$. 
Let us denote the associated bilinear form on 
$ \Delta_{\underline{w}, \cal R}(y)$ by $ \langle \cdot, \cdot \rangle_{y} $ as well. Motivated by Jantzen's original work we now define 
\begin{equation}{\label{Jantzen}}
\Delta_{\underline{w}, \cal R}^i(y) := \{ v \in  \Delta_{\underline{w}, \cal R}(y) | \langle v, w \rangle_{y} \in x^i {\cal R} \mbox{ for all } 
w  \in  \Delta_{\underline{w}, \cal R}(y) \}.
\end{equation}
Since $ \langle \cdot, \cdot \rangle_{y}$ is $A_{\underline{w},{\cal R}} $-invariant this gives to a decreasing filtration 
\begin{equation}{\label{Jantzenfiltration}}
\Delta_{\underline{w}, \cal R}(y) = \Delta_{\underline{w},\cal R}^0(y)  \supseteq \Delta_{\underline{w},\cal R}^1(y) \supseteq \Delta_{\underline{w},\cal R}^2(y)  \supseteq \cdots
\end{equation}
of $ A_{\underline{w}, \cal R} $-modules. 

\medskip
From $ \cal R $ we can further extend scalars to $ {\mathbb R}  $ via the homomorphism $ {\cal R} \rightarrow {\cal R} /(x)  \cong {\mathbb R} $.
We denote the corresponding cellular algebra by $ A_{\underline{w},{{\mathbb R}}}$ and its cell modules by 
$ \Delta_{\underline{w},  {\mathbb R}}(y) := \Delta_{\underline{w}, \cal R}(y) \otimes_{\cal R} {\mathbb R} $. We denote the bilinear form 
on $ \Delta_{\underline{w},  {\mathbb R}}(y)  $ 
by 
$ \langle \cdot, \cdot \rangle_{ y}$, too.

\medskip
In \cite{P}, D. Plaza studied the representation theory of the
graded cellular algebra $ A_{\underline{w},{{\mathbb R}}}$. Using \cite{EW}, 
he showed that the nonzero graded decomposition numbers for $ A_{\underline{w},{\mathbb  R}}$
are Kazhdan-Lusztig polynomials for $ W$. In other words,  
$ A_{\underline{w},{\mathbb  R}}$ is a highly interesting non-semisimple algebra. 
Although it is not quasi-hereditary, it still has many features in common with category $ \cal O $ for a
complex semisimple Lie algebra.

\medskip
The inclusion $ \Delta_{\underline{w}, \cal R}^i(y) \subseteq \Delta_{\underline{w},\cal R}(y) $ induces a homomorphism
$ \Delta_{\underline{w}, {\mathbb R}}^i(y)  \rightarrow \Delta_{\underline{w},{\mathbb R}}(y) $ and we define 
$ \Delta_{\underline{w},{\mathbb R}}^i(y)  \subseteq \Delta_{\underline{w},{\mathbb R}}(y) $
as the image of this homomorphism. 
Via ({\ref{Jantzenfiltration}}) this gives rise to a filtration of $  A_{\underline{w},{{\mathbb R}}} $-modules
\begin{equation}{\label{Jantzenfiltration1}}
\Delta_{\underline{w},{\mathbb R}}(y) = \Delta_{\underline{w},{\mathbb R}}^0(y)  \supseteq \Delta_{\underline{w},{\mathbb R}}^1(y) \supseteq \Delta_{\underline{w},{\mathbb R}}^2(y)  \supseteq \ldots 
\end{equation}
This filtration of $ \Delta_{\underline{w},{\mathbb R}}(y) $ is our analogue
of the Jantzen filtration for Verma modules and 
the following Theorem is a weak analogue of Jantzen's sum formula, involving dimensions of the modules.
It follows from the Theorem that the filtration is finite.
\begin{teo}{\label{sumformula}} 
We have that
$\dim \Delta_{\underline{w}, {\mathbb R}}(y)/ \Delta_{\underline{w}, {\mathbb R}}^1(y) = \dim  \Delta_{\underline{w},{\mathbb R}}(y)/ {\rm rad}  \langle \cdot, \cdot \rangle_{ y}$. 
Moreover, the following sum formula holds
\begin{equation}{\label{sumform}} 
\sum_{i>0}  \dim \Delta_{\underline{w},{\mathbb R}}^i(y) = \sum_{\beta > 0, s_{\beta} y > y } \dim \Delta_{\underline{w}, {\mathbb R}}( s_{\beta} y ) .     
\end{equation}
\end{teo}
\begin{dem}
The proof follows closely the proof in the classical case, see 'Key Lemma' of section 5.6 of \cite{H}. 
For the reader's convenience we here sketch the argument. 
We first prove ({\ref{sumform}}).
Let $ s := {\rm rank}\,  \Delta_{\underline{w},{\cal R}}(y) $. 
Since $ \cal R $ is a principal ideal domain, the matrix $ M $ for 
$  \langle \cdot, \cdot \rangle_{y} $
is equivalent to a diagonal $ \cal R$-matrix 
$ D := {\rm diag}(d_1, \ldots, d_s) $, in other words there are invertible $ \cal R$-matrices $ A,B$ of sizes $ s \times s $
such that $ D = AMB $. 
For $ j=1, \ldots, s$ we let $ v_j $ be the basis vector corresponding to the $j$'th column of $ D$
and we set $ a_j := \nu_x(d_j ) $.
Then $ \Delta_{\underline{w}, {\cal R}}^i(y) $ is spanned by the elements $ \{ v_j | a_j \ge i \} $ together with the 
elements $ \{ x^{ i- a_j} v_j | a_j < i \} $. The last elements vanish when tensoring over $ \mathbb R $ and so the
left hand side of  (\ref{sumform}) is $ a_1 + a_2 + \ldots + a_s$.
But this is equal to $ \nu_x(\det D ) $ which is equal to the 
right hand side of (\ref{sumform}) by Theorem {\ref{lightleavesJM_det}} and the assumption on $ \varphi$. 

The equality $\dim \Delta_{\underline{w},{\mathbb R}}(y)/ \Delta_{\underline{w},{\mathbb R}}^1(y) = 
\dim \Delta_{\underline{w},{\mathbb R}}(y)/ {\rm rad}  \langle \cdot, \cdot \rangle_{y}$ 
is proved similarly.
\end{dem}

\medskip
We next explain one of several ways to obtain a sum formula in the positive characteristic situation.
Let ${ \mathfrak h}   $ be the Cartan matrix representation of a Weyl group $ W$ over
$ \overline{\F}_p$, as introduced above.
Choose this time $ {\cal R } := \Z_{(p)} $, the localization of $ \Z $ at $p$. 
Then $ {\cal R }  $ is also a discrete valuation
domain with valuation function that we denote $ \nu_p(  \cdot )$. Its maximal ideal is
$ p{\cal R } $ and we have $ {\cal R}/p {\cal R } \cong {\F_p}$. 
Let $ \varphi: {\mathfrak h}_{\Z} \rightarrow \Z $ be the homomorphism of $ \mathbb Z $-modules given 
by $ \varphi (\alpha_s) := p $ for the basis elements $ \{ \alpha_s \mid s \in S \} $ 
and denote also by $ \varphi$ the
homomorphism obtained from $ \varphi $ by composing with the inclusion $ \Z \subseteq {\cal R } $. 

\medskip
We define $ A_{\underline{w},{\cal R}}:= A_{\underline{w}} \otimes_R \cal R $ , 
$ A_{\underline{w},{{\F_p}}} := A_{\underline{w},{\cal R}} 
\otimes_{\cal R } {\F_p}$ and so on, mimicking what we did before.
We then obtain an
$ A_{\underline{w},{\cal R}} $-module filtration of $ \Delta_{\underline{w}, \cal R}(y) $ via
\begin{equation}{\label{JantzenCharp}}
\Delta_{\underline{w}, \cal R}^i(y) := \{ v \in  \Delta_{\underline{w}, \cal R}(y) | \langle v, w \rangle_{y} \in p^i {\cal R} \mbox{ for all } 
w  \in  \Delta_{\underline{w}, \cal R}(y) \}
\end{equation}
and this induces a filtration $ \{ \Delta_{\underline{w}, {\F_p}}^i(y) | i =0,1, \ldots  \} $ of $ \Delta_{\underline{w}, {\F_p}}(y)  $
and hence also a filtration $ \{ \Delta_{\underline{w}, \overline{\F}_p}^i(y) | i =0,1, \ldots  \} $ of $  \Delta_{\underline{w}, \overline{\F}_p}(y)  $
via $  \Delta_{\underline{w}, \overline{\F}_p}^i(y) :=  \Delta_{\underline{w}, \F_p}^i(y) \otimes_{\F_p} \overline{\F}_p$.
We get the following Theorem.
\begin{teo}{\label{sumformulaA}} 
We have that
$\dim \Delta_{\underline{w}, {\overline{\F}_p }}(y)/ \Delta_{\underline{w}, {\overline{\F}_p}}^1(y) = \dim  \Delta_{\underline{w},{\overline{\F}_p}}(y)/ {\rm rad}  \langle \cdot, \cdot \rangle_{ y}$. 
Moreover, the following sum formula holds
\begin{equation}{\label{sumformchar}}
\sum_{i>0}  \dim \Delta_{\underline{w},{\overline{\F}_p}}^i(y) = 
\sum_{\beta > 0, s_{\beta} y > y }  \nu_{p}( \varphi( \beta )) \dim \Delta_{\underline{w}, {\overline{\F}_p}}( s_{\beta} y ) .     
\end{equation}
\end{teo}
\begin{dem}
The proof is essentially the same as the proof of Theorem {\ref{sumformula}}. 
\end{dem}

\medskip
In the classical theory the sum formula is formulated as an equality in the Grothendieck group and thus 
we would have expected ({\ref{sumform}}) and ({\ref{sumformchar}}) to hold at this level of generality. Our next goal is to show that in fact there is a variation of 
the above constructions for which ({\ref{sumform}}) and 
({\ref{sumformchar}}) do 
hold at the Grothendieck group level. 

\medskip
As can be seen already in the one-colour calculations presented above, the algebra $ A_{\underline{w},{{\mathbb R}}}$
is not quasi-hereditary in general, in other words we have $ \Lambda \neq \Lambda_0 $ in general. Moreover, 
in general it appears to be difficult to determine $ \Lambda_0  $, and hence for the Grothendieck 
group to work well we need to change the setup.

\medskip
As mentioned above, $\cal D $ is cellular category in the sense of Westbury. Let us recall his definition from \cite{Wes}: 
\begin{defi}{\label{cellularcat}}
Let $ \Bbbk $ be a commutative ring with identity and let $ \cal C $ be a $ \Bbbk $-linear category with duality $ \ast$. 
Then $ \cal C $ is called a cellular category if there exists a poset $ \Lambda$ and for each 
$ \lambda \in \Lambda $ and each object $ n $ in $ \cal D $ a finite set $ \tab(n, \lambda) $ together with a map 
$ \tab(m, \lambda) \times  \tab(n, \lambda) \rightarrow Hom_{\cal C}(m,n), (S,T) \mapsto C_{ST}^{\lambda} $, satisfying 
$(C_{ST}^{\lambda})^{\ast} = C_{TS}^{\lambda}$. 
These data satisfy that 
\begin{equation}
\begin{split}
\{ C_{ST}^{\lambda} | S \in \tab(m, \lambda), T \in \tab(n, \lambda), \lambda \in \Lambda \} \mbox{ is a basis for } Hom_{\cal C}(m,n)  \\
\end{split}
\end{equation}
and for all $ a \in Hom_{\cal C}(n,p), S \in \tab(m, \lambda), T \in \tab(n, \lambda) $ 
\begin{equation}
a  C_{ST}^{\lambda} = \sum_{ S^{\prime} \in \tab(p, \lambda)} r_{a}(S^{\prime}, S) C_{ S^{\prime}, T}^{\lambda}  \mbox{ mod } A^{\lambda} 
\end{equation}
where $ A^{\lambda} $ is the span of $\{  C_{ST}^{\mu } | \mu < \lambda, S \in \tab(m, \mu), T \in \tab(p,\mu) \} $.
\end{defi}

As already mentioned in \cite{Wes}, each object $ m $ in a cellular category $ \cal C $ gives rises to the cellular algebra
$ {\rm End}_{\cal R}(m) $. This fact can be generalized as follows. Let $ I $ be any finite subset of the objects of $ \cal C$
and define $ {\rm End}_{\cal C}(I) $ as the $ \Bbbk$-direct sum 
\begin{equation}
 {\rm End}_{\cal C}(I)  := \oplus_{ m,n \in I } {\rm Hom}_{\cal C}(m,n). 
\end{equation}
Then $ {\rm End}_{\cal C}(I) $ has a natural $\Bbbk$-algebra structure as follows
\begin{equation}
g \cdot f := \left\{\begin{array}{ll} gf  & \mbox{ if } f \in {\rm Hom}_{\cal C}(m,n), g \in {\rm Hom}_{\cal C}(n,p) \mbox{ for some } m,n,p \\
0 & \mbox{ otherwise.} 
\end{array}
\right. 
\end{equation}
\begin{teo}
Let $ \cal C $ be a cellular category as in Definition {\ref{cellularcat}}
and define for $ \lambda \in \Lambda $ the set $ \tab(\lambda) := \cup_{ n \in I } \tab(n, \lambda) $.
Let for $ S \in \tab(\lambda), T \in \tab(\lambda) $ the element $ C_{ST}^{\lambda} \in    {\rm End}_{\cal C}(I)  $ be
defined as the corresponding inclusion of $ C_{ST}^{\lambda} \in {\rm Hom}_{\cal C}(m,n) $ in $ {\rm End}_{\cal C}(I)  $. 
Then these data define a cellular algebra structure on $ {\rm End}_{\cal C}(I) $. 
\end{teo}
\begin{dem}
This is immediate from the definitions.
\end{dem}

\medskip
For our category $ \cal D $ the objects are $ {\bf{exp}}_s$ and
the tab function is
\begin{equation} {\rm  Tab}(\underline{w}, y) := \{ \e \in \{0,1 \}^k | \e \mbox{ subexpression of } \underline{w} 
  \mbox{ expressing } y \}.
\end{equation}
Let us now 
choose an ideal $ \pi  $ in $ W $, that is $ x \in \pi, y < x \Rightarrow y \in \pi $, 
and let us also choose for all $ y \in \pi $ an arbitrary $ \underline{y} \in {\bf{rexp}}_s $ expressing $ y $.
We set $ \underline{\pi } := \{ \underline{y} | y \in \pi \} $
and define the algebra 
\begin{equation} A_{{\pi}} = A_{\underline{\pi}}   :=  {\rm End}_{\cal D} ( \oplus_{ y \in \pi } { \underline{y}}).
\end{equation} 
Then, according to the above Theorem we have that 
$ A_{{\pi}} $ is a cellular algebra
on the 
poset $ \Lambda_{\pi }:=\pi $. For $ y \in \Lambda_{\pi} $ the corresponding 
$ \tab(y) $ is 
$ \tab_{{\pi}}(y):= \bigcup_{ z \in \pi} \tab({\underline{z}}, y) $. 
For $ y \in \pi $ we denote by $ \Delta_{{\pi}}(y) $ the corresponding cell module for $ A_{{\pi}} $. Its 
basis is $ \bigcup_{ z \in \pi } \{ {\LL}_{ \underline{z}, \bf e} | \, {\bf e} \in \tab({\underline{z}}, y) \} $,
and there is an
$ R$-module decomposition 
\begin{equation}{\label{bigcellmodule}}
 \Delta_{{\pi}  }(y)  = \oplus_{ z \in \pi }  \Delta_{\underline{z} }(y) 
\end{equation}
and 
in particular 
$ \dim \Delta_{{\pi}  }(y) = \sum_{ z \in \pi } \dim \Delta_{\underline{z} }(y) $. Note that there is no similar 
formula for $  A_{{\pi}}. $ We view ({\ref{bigcellmodule}}) as a kind of weight space decomposition for 
$  \Delta_{{\pi}  }(y) $ with weight spaces $  \Delta_{\underline{z} }(y) $.

\medskip
We point out that the bilinear form $  \langle \cdot, \cdot \rangle_{ \pi, y} $
on $  \Delta_{{\pi}  }(y) $
is orthogonal
with respect to the decomposition in ({\ref{bigcellmodule}}), as follows directly from the definitions.
This is a key observation for the following. 

\medskip
Recall that the ground ring for $ A_{\pi} $ is $ R$. For a field $ \Bbbk $ that is made into an $ R$-algebra via a homomorphism $ R \rightarrow \Bbbk $, 
we obtain as usual a specialized algebra $ A_{\pi, \Bbbk } := A_{\pi } \otimes_R \Bbbk  $.
This is also a cellular algebra.
We now state the result that makes us prefer
$ A_{\pi}  $ over $ A_{\underline{w}} $. 
\begin{teo}{\label{lightleavesJM}}
$  A_{\pi}  $ and $ A_{\pi, \Bbbk}  $ are quasi-hereditary algebras. 
\end{teo}
\begin{dem}
By remark (3.10) of \cite{GL} we must show that for $ y \in  \Lambda_{\pi} $ 
we have that $  \langle \cdot, \cdot \rangle_{\pi, y} \neq 0 $. But since $ \underline{y}   $ 
is a reduced expression for $y$ we get for $ \e =(1,1,\ldots, 1) $ that 
$ {\LL}_{ \underline{y}, \bf e} $ can be chosen as the identity morphism in $ {\rm End}_{\cal D} (\underline{y})$
and so $ \langle {\LL}_{ \underline{y}, \bf e} , {\LL}_{ \underline{y}, \bf e}  \rangle_{\pi, y} = 
\langle {\LL}_{ \underline{y}, \bf e} , {\LL}_{ \underline{y}, \bf e}  \rangle_{y}  = 1$.
The Theorem follows from this. 
\end{dem}

\medskip
Let $ {\cal R } := {\mathbb R}[x] $ or $ {\cal R } := { \mathbb Z}_{(p)} $ depending 
on our choice of realization
and define $  A_{\pi, \cal R },   A_{\pi, \mathbb R } $ and $ A_{\pi, \overline{\F}_p }  $ correspondingly. 
By the Theorem, these are all quasi-hereditary algebras with cell modules that we denote 
$  \Delta_{\pi, \cal R }(y), \Delta_{\pi, \mathbb R }(y) $ and $ \Delta_{\pi, \overline{\F}_p }(y)  $.
We use the same notation $  \langle \cdot, \cdot \rangle_{ \pi, y} $ for the bilinear form 
on each of these modules and introduce 
the Jantzen type filtrations $  \Delta_{\pi, \cal R }^i(y), \Delta_{\pi, \mathbb R }^i(y) $ and $ \Delta_{\pi, \overline{\F}_p }^i(y)  $,
mimicking the previous construction.

% \medskip
% On the other hand, the results on determinants that were obtained in the previous sections are valid for $ A_{\pi} $ as well. 
% In particular we have the following generalization of Theorem {\ref{shapovalov}}.
% \begin{teo}{\label{shapovalov_version2}}
% We have the following generalization of Theorem {\ref{shapovalov}}
% \begin{equation}
% \det  \langle \cdot, \cdot \rangle_{\pi, y} = \pm  \prod_{  \ms \beta > 0, s_{\beta} y > y }   \beta^{ \dim \Delta_{\underline{\pi}} ( s_{\beta} y)}. 
% \end{equation}
% \end{teo}
% \begin{dem}
% By the definitions, the decomposition ({\ref{bigcellmodule}}) is orthogonal, that is
% $ \langle \cdot, \cdot \rangle_{\pi, y} = \bigoplus_{ y \in \pi }  \langle \cdot, \cdot \rangle_{ y} $. 
% Since $ \dim \Delta_{\pi  }(y) = \sum_{ z \in \pi } \dim \Delta_{\underline{z} }(y) $, we conclude 
% via Theorem {\ref{shapovalov}} that the fomula of the Theorem is correct.
% \end{dem}

% \medskip
% From this we get the following generalization of Theorem {\ref{sumformula}}. 
% Using the orthogonality of the decomposition ({\ref{bigcellmodule}}), the 
% proof is the same as the proof of Theorem {\ref{sumformula}}. 
% \begin{teo}{\label{sumformula1}} 
% We have that
% $\dim \Delta_{{\pi}, \Bbbk}(y)/ \Delta_{\pi, \Bbbk}^1(y) = \dim  \Delta_{{\pi},\Bbbk}(y)/ {\rm rad}  \langle \cdot, \cdot \rangle_{\pi, y}$. 
% Moreover, the following sum formula holds
% \begin{equation}{\label{sumform}}
% \sum_{i>0}  \dim \Delta_{{\pi},\Bbbk}^i(y) = \sum_{\beta > 0, s_{\beta} y > y } \dim \Delta_{{\pi}, \Bbbk}( s_{\beta} y ) .     
% \end{equation}
% \end{teo}

\medskip
Let us now turn to the Grothendieck groups. 
Let {\color{black}{$\Bbbk$}}
be either $  \mathbb R $ or $\overline{\F}_p $, depending on the choice of $ \mathfrak h $.
Let $ \langle A_{{\pi}, \Bbbk} \mbox{-mod} \rangle $ (resp. 
$ \langle A_{\underline{w}, k} \mbox{-mod} \rangle $) be the Grothendieck groups of finite dimensional 
$  A_{{\pi}, \Bbbk} $-modules (resp. finite dimensional $  A_{\underline{w}, \Bbbk} $-modules).
Because of quasi-heredity, the classes $ [ \Delta_{{\pi}, \Bbbk}(y) ]   $
in $ \langle A_{{\pi}, \Bbbk} \mbox{-mod} \rangle $
of the cell modules 
form a $ \Z$-basis for $ \langle A_{{\pi}, \Bbbk} \mbox{-mod} \rangle $.
On the other hand, the corresponding classes $ [ \Delta_{{\underline{w}}, \Bbbk }(y)]  $
in $ \langle A_{\underline{w}, \Bbbk} \mbox{-mod} \rangle $
only form a generating set for $ \langle A_{\underline{w}, \Bbbk} \mbox{-mod} \rangle $ since 
$  A_{\underline{w}, \Bbbk} $ is not quasi-hereditary in general.

Let us consider the following 
projection map $ \varphi_{\underline{w}} $:
\begin{equation}{\label{grothendieckhomo}}
\varphi_{\underline{w}}:  \langle A_{{\pi, \Bbbk}} \mbox{-mod} \rangle  \rightarrow 
\langle A_{\underline{w}, \Bbbk} \mbox{-mod} \rangle,  \, \, \,  [ \Delta_{{\pi}, \Bbbk}(y) ]
\mapsto [ \Delta_{\underline{w}, \Bbbk}(y) ].
\end{equation}
We need to describe $ \varphi_{\underline{w}} $ in a different way. 
There is a natural diagonal subalgebra $ \oplus_{ w \in \pi } A_{\underline{w}, \Bbbk} $ of $ A_{\pi, \Bbbk} $ having 
$ A_{\underline{w}, \Bbbk} $ as an algebra summand. Let 
\begin{equation}
  \tau_{\underline w}:  A_{{\pi}, \Bbbk} \mbox{-mod} \rightarrow A_{\underline{w}, \Bbbk} \mbox{-mod}
\end{equation}  
be the composition of the corresponding restriction 
and idempotent truncation functors. Then $\tau_{\underline w} $
induces a 
$ \mathbb Z$-module homomorphism between the Grothendieck groups: this is our $ \varphi_{\underline{w}} $. 
Via this description of $ \varphi_{\underline{w}}$ and the orthogonality of the decomposition ({\ref{bigcellmodule}}), 
we get the following compatibility of the Jantzen filtrations
\begin{equation}{\label{compatible}}
\varphi_{\underline{w}}( [ \Delta_{{\pi}, \Bbbk}^i(y) ]) = [ \Delta_{\underline{w}, \Bbbk}^i(y) ].
\end{equation}

Let now $ \dim_{\underline{w}}:  \langle A_{\underline{w}, \Bbbk} \mbox{-mod} \rangle \rightarrow {\mathbb Z}, [M] \mapsto \dim_{\Bbbk} M   $
be the dimension homomorphism and define a $ \mathbb Z$-module homomorphism
$ \Phi: \langle A_{{\pi}, \Bbbk} \mbox{-mod} \rangle \rightarrow \oplus_{\underline{w}\in  \underline{\pi}} {\mathbb Z} $
by setting the $\underline{w}$'th coordinate equal to the composite $ \dim_{\underline{w}} \circ \,  \varphi_{\underline{w}}$. 
\begin{lem}
In the above notation $ \Phi: \langle A_{{\pi}, \Bbbk} \mbox{-mod} \rangle \rightarrow \oplus_{\underline{w}\in  \underline{\pi}} {\mathbb Z} $ 
is an isomorphism of $ \mathbb Z$-modules.
\end{lem}
\begin{dem}
Since $ \langle A_{{\pi}, \Bbbk} \mbox{-mod} \rangle $ and 
$ \oplus_{\underline{w}\in  \underline{\pi}} {\mathbb Z} $ are free 
$ \mathbb Z$-modules
of the same rank, it is enough to show that 
$ \Phi$ is surjective. Let $ (n_{\underline{z}})_{ \underline{z} \in \underline{\pi}} $ be an element of 
$  \oplus_{\underline{z}\in  \underline{\pi}} {\mathbb Z} $. 
Choose $ \underline{z}_0 $ 
satisfying $ n_{\underline{z}_0} \neq 0 $ and $ {z}_0 $ minimal in $ W $ with respect to this.
The $ \underline{z}_0$'th component of $ \Phi(  [\Delta_{{\pi}}(z_0) ]) $ is $ \dim_{\Bbbk} \Delta_{\underline{z}_0}( z_0)= 1 $
and so the 
$ \underline{z}_0$'th component of $ \Phi( n_{ \underline{z}_0} [\Delta_{{\pi}}(z_0) ]) $ is 
$ n_{ \underline{z}_0} $. Moreover, the $ \underline{z}$'th component of 
$ \Phi(  [\Delta_{{\pi}}(z_0) ]) $ is nonzero only if $ z_0 \le \underline{z} $, or equivalently
$ z_0 \le {z} $ since $  \underline{z} \in {\bf{rexp}}_s$.
Hence, we can use induction  on 
$ (n_{\underline{z}})_{ \underline{z} \in \underline{\pi}} - \Phi( n_{ \underline{z}_0} [\Delta_{\underline{\pi}}(z_0) ]) $ 
and get that $ (n_{\underline{z}})_{ \underline{z} \in \underline{\pi}} \in im \, \Phi$ as claimed.
\end{dem}

\medskip
We are now in position to prove the main Theorems of this section. 
\begin{teo}{\label{firstmmain}}
Let $ {\mathfrak h } $ be the geometric representation over $ \mathbb R$.
For the $A_{{\pi}, \mathbb R} $-filtration $ \{ \Delta_{\pi, {\mathbb R}}^i( y) \} $ of $  \Delta_{\pi, {\mathbb R}}( y) $ we have that 
$ \Delta_{{\pi}, {\mathbb R}}(y)/ \Delta_{\underline{\pi}, {\mathbb R}}^1(y) $ is nonzero and irreducible and
the following sum formula holds 
in $ \langle A_{{\pi}, {\mathbb R}} \mbox{-mod} \rangle $ 
\begin{equation}{\label{sumformGR}}
\sum_{i>0}  [ \Delta_{{\pi},{\mathbb R}}^i(y)] = \sum_{\beta > 0, s_{\beta} y > y } [ \Delta_{{\pi}, {\mathbb R}}( s_{\beta} y) ] .     
\end{equation}
\end{teo}
\begin{dem}
By the construction of the filtration, 
the first statement follows from the quasi-heredity of $ A_{{\pi}}$.

To show the second statement we 
get by applying $ \Phi $ to the left hand side of ({\ref{sumformGR}}) the element of 
$ \oplus_{\underline{w}\in  \underline{\pi}} {\mathbb Z} $ whose $ \underline{w}$'th component is
$ \sum_{i>0}  \dim \Delta_{\underline{w},{\mathbb R}}^i(y) $. Applying  
$ \Phi $ to the right hand side of ({\ref{sumformGR}}), we get the element whose 
$ \underline{w}$'th component is 
$\sum_{\beta > 0, s_{\beta} y > y } \dim \Delta_{\underline{w}, {\mathbb R}}( s_{\beta} y ) $ and so by Theorem 
{\ref{sumformula}} the two sides are equal and the Theorem follows.

\end{dem}

Similarly, we have the following Theorem.
\begin{teo}{\label{Similarly, we have the following Theorem}}
  Let $ W $ be a Weyl group and let $ {\mathfrak h } $ be its Cartan matrix representation over $ \overline{\F}_p$.
For the $A_{{\pi}, \overline{\F}_p} $-filtration $ \{ \Delta_{\pi, {\overline{\F}_p}}^i( y) \} $ of $  \Delta_{\pi, {\overline{\F}_p}}( y)  $ we have that 
$ \Delta_{{\pi}, {\overline{\F}_p}}(y)/ \Delta_{{\pi}, {\overline{\F}_p}}^1(y) $ is nonzero and irreducible and
the following sum formula holds 
in $ \langle A_{{\pi}, {\overline{\F}_p}} \mbox{-mod} \rangle $ 
\begin{equation}{\label{sumformGRchar}}
\sum_{i>0}  [ \Delta_{{\pi},{\overline{\F}_p}}^i(y)] = \sum_{\beta > 0, s_{\beta} y > y } \nu_{p}( \varphi( \beta ))  [ \Delta_{{\pi}, {\overline{\F}_p}}( s_{\beta} y) ] .     
\end{equation}
\end{teo}

%% \medskip
%% An important aspect of the diagrammatcal category $ \cal D $ 
%% is its $ \mathbb Z$-grading. For general reasons, it induces a cellular algebra structure on its 
%% zeroth graded subcategory and then also on the zeroth graded component $ A_{ {\underline{w}}}^{0}$ of $ A_{ {\underline{w}}}$. 
%% Unfortunately, our $L_i$'s are of degree 2 and therefore do not induce JM-operators on $ A_{ {\underline{w}}}^{0}$ 
%% and hence, in particular, our methods do not give rise to a sum formula for $ A_{ {\underline{w}}}^{0}$.

\section{Applications of the sum formula.} 
We indicate in this section, via an example, how to apply formula ({\ref{sumformGR}}) to obtain decomposition numbers for 
$ A_{{\pi}}$. This is parallel to Jantzen's original calculations for Verma modules and just as in the original setting, the sum formula 
only gives complete information on the decomposition numbers in certain small cases.

\medskip
The comparison with Verma modules is strengthened 
by the following 
analogue for $ A_{{\pi}}$ of the fact that homomorphisms between Verma modules
are injective. For Verma modules this is a consequence of the PBW-Theorem
together with the definition of Verma modules as
induced modules from a Borel subalgebra. In the  
$ A_{{\pi}}$-setting, these notions and results are not available and so one must argue in
a different way to prove the injectivity. 
We shall rely on a result of D. Plaza, see \cite{P}.

\begin{teo}\label{weshallrelyPlaza}
  i). Let $ W $ be an arbitrary Coxeter group and let $ \pi $ an ideal of $W $. Let $ \mathfrak h $ be
  the geometric representation of $ W$. 
  Suppose that $ u,v \in W  $ satisfy $ u \le v $. Then there is an embedding
of $ A_{{\pi, \mathbb R}}$-modules 
$\Delta_{{\pi}, \mathbb R}(v) \subseteq \Delta_{{\pi}, \mathbb R}( u)$.
\newline
ii). Let $ W $ be a Weyl group and let $ \pi $ be an ideal of $W $. Suppose that $ u,v \in W  $ satisfy $ u \le v $.
Let $ \mathfrak h $ be the Cartan matrix representation of $ W $ with $ p > 3 $. Then there is an embedding
of $ A_{\pi, \overline{\F}_p}$-modules 
$\Delta_{\pi, \overline{\F}_p}(v) \subseteq \Delta_{\pi, \overline{\F}_p}( u)$.
\end{teo}
\begin{dem}
  As already mentioned above, in each of the settings i) and ii) we have that $ \mathfrak h $ is
  a Soergel realization and so  
  the functor $ F: {\mathcal D} \rightarrow  {\mathbb B} {\mathbb S}\rm Bim $ introduced in section 3
  defines an equivalence of categories.
  For $ \underline{w} \in {\bf{rexp}}_s $ we define
$ A_{\underline{w}}^{{\mathbb B} {\mathbb S}} := { \rm End}_{{\mathbb B} {\mathbb S}} (\underline{w}) $ 
  and then $ F $ induces 
  an algebra isomorphism $A_{\underline{w}} \rightarrow A_{\underline{w}}^{{\mathbb B} {\mathbb S}} $. 
   Let us denote by 
$ \Delta^{{\mathbb B} {\mathbb S}}_{\underline{w}}(y) $
the cell module for $ A_{\underline{w}}^{{\mathbb B} {\mathbb S}} $ with associated basis 
$ \{ {\LL}^{{\mathbb B} {\mathbb S}}_{ \underline{w}, \bf e} | \, {\bf e} \in {\rm Tab}(y) \} $.
Then we have that $ F( {\LL}_{ \underline{w}, \bf e}) =
{\LL}^{{\mathbb B} {\mathbb S}}_{ \underline{w}, \bf e}$. 
Suppose now that $ \underline{u}, \underline{v} \in {\bf{rexp}}_s $ satisfy $ u \le v \le w  $
and let $ G_v^u \in \Hom_{{\mathbb B} {\mathbb S}}(B_{\underline{v}}, B_{\underline{u}}) $ be
the light leaves morphism associated with Deodhars' 
distinguished subexpression for $ \underline{u} $ in $ \underline{v} $.
Then in \cite{P} it was proved that the map given by 
\begin{equation}
  \Delta^{{\mathbb B} {\mathbb S}}_{\underline{w}}(v) \rightarrow \Delta^{{\mathbb B} {\mathbb S}}_{\underline{w}}(u),
  l \mapsto G_v^u \circ l  
\end{equation}  
defines an injective homomorphism of $ A_{\underline{w}}^{{\mathbb B} {\mathbb S}}$-modules.
Note that although the argument in \cite{P} is presented for the ground field $ \mathbb R $, it carries
over to characteristic $ p $. 
But then by pulling back via $ F $ we get that 
\begin{equation}
  \Delta_{\underline{w}}(v) \rightarrow \Delta_{\underline{w}}(u),
  l \mapsto F^{-1}(G_v^u) \circ l  
\end{equation}  
is an injective homomorphism of $ A_{\underline{w}}$-modules. 
The Theorem follows from this by summing over all relevant summands.
\end{dem}

\medskip
\noindent
{\bf Remark}.
We believe that 
the Theorem can be proved diagrammatically, without relying on the equivalence $ F $.

\medskip
Let us finish the paper by calculating some decomposition numbers, using Theorem
{\ref{firstmmain}} and Theorem
{\ref{Similarly, we have the following Theorem}}.
We first choose $ \Bbbk := \mathbb R $ and 
$ W=\langle {\color{red} s_1}, {\color{blue} s_2}  \rangle  $ of type $ A_2 $, that 
is $$ W = \{ 1, {\color{red} s_1}, {\color{blue} s_2}, {\color{red} s_1} {\color{blue} s_2}, 
{\color{blue} s_2}{\color{red} s_1},{\color{red} s_1} {\color{blue} s_2}{\color{red} s_1}= 
{\color{blue} s_2}{\color{red} s_1} {\color{blue} s_2} \}.$$
Let us use $ \pi := W $ and 
$ \underline{\pi} := 
\{ 1, {\color{red} s_1}, {\color{blue} s_2}, {\color{red} s_1} {\color{blue} s_2}, 
{\color{blue} s_2}{\color{red} s_1},{\color{red} s_1} {\color{blue} s_2}{\color{red} s_1}  \}. $

\medskip
We associate with this the following alcove geometry
\medskip
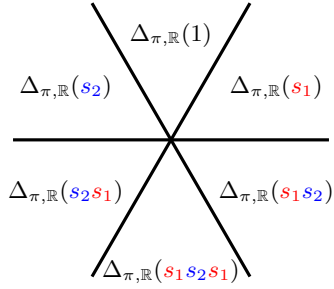
\begin{equation} \begin{aligned}
\begin{tikzpicture}[xscale=0.7, yscale=0.7]
    \foreach\ang in {60,120,...,360}{
     \draw[very thick] (0,0) -- (\ang:3cm);
    }

 \node at (0,2) {\small $\Delta_{{\pi}, \mathbb R}(1)$};
 \node at (0,-2.5) {\small $ \Delta_{{\pi}, \mathbb R} ({\color{red} s_1} {\color{blue} s_2}{\color{red} s_1})$};
 \node at (2,1) {\small $ \Delta_{{\pi}, \mathbb R}({\color{red} s_1})$};
 \node at (-2,1) {\small $ \Delta_{{\pi}, \mathbb R}({\color{blue} s_2})$};
 \node at (2,-1) {\small $ \Delta_{{\pi}, \mathbb R}({\color{red} s_1} {\color{blue} s_2})$};
 \node at (-2,-1) {\small $  \Delta_{{\pi}, \mathbb R}({\color{blue} s_2}{\color{red} s_1})$};
\end{tikzpicture} \end{aligned}
\end{equation}
The images of the simple modules $  [L_{{\pi}, \mathbb R}(y)]  $, with $ y $ running over W,  give the canonical basis of 
$ \langle A_{{\pi}, \mathbb R} \mbox{-mod} \rangle $ and by 
general cellular algebra theory we have that 
\begin{equation}{\label{cellu}}[\Delta_{{\pi}, \mathbb R}(y)] =  \sum_{y  \le u } d_{yu}^{ \, {\mathbb R}} [L_{{\pi}, \mathbb R}(u)] 
\end{equation}
where the $ d_{yu}^{\mathbb R} $'s are the decomposition numbers. Let us determine the 
expansion of all $ [\Delta_{{\pi}, \mathbb R}(y)] $'s in terms of $ [L_{{\pi}, \mathbb R}(y)] $'s, or equivalently 
the decomposition numbers $ d_{yu}^{\, \mathbb R}$.
As a starting point we get from ({\ref{cellu}}) that 
\begin{equation}
 \Delta_{{\pi}, \mathbb R}({\color{red} s_1} {\color{blue} s_2}{\color{red} s_1})= 
L_{{\pi}, \mathbb R}({{\color{red} s_1}\color{blue} s_2}{\color{red} s_1})
\end{equation}
which gives us all the decomposition numbers $ d_{yu}^{\, \mathbb R}$ with $ y = {\color{red} s_1} {\color{blue} s_2}{\color{red} s_1}$.
Let us now calculate $  [\Delta_{{\pi}, \mathbb R}( {\color{blue} s_2}{\color{red} s_1})] $. By ({\ref{cellu}}) we have 
that 
$[ \Delta_{{\pi}, \mathbb R}( {\color{blue} s_2}{\color{red} s_1})]= 
[L_{{\pi}, \mathbb R}({\color{blue} s_2}{\color{red} s_1})] + d  
[L_{{\pi}, \mathbb R}({\color{red} s_1}{\color{blue} s_2}{\color{red} s_1})]
$
for some $ d  $.
The sum formula ({\ref{sumformGR}}) reads in this case
\begin{equation}
\sum_{i>0}  [ \Delta_{{\pi},\mathbb R}^i({\color{blue} s_2}{\color{red} s_1} )] 
=  [ \Delta_{{\pi}, {\mathbb R}} ({\color{red} s_1}{\color{blue} s_2}{\color{red} s_1}) ]  = 
[ L_{{\pi}, {\mathbb R}} ({\color{red} s_1}{\color{blue} s_2}{\color{red} s_1}) ] .     
\end{equation}
Since $ \Delta_{{\pi},\mathbb R}^1({\color{blue} s_2}{\color{red} s_1} ) = {\rm rad}  \langle \cdot, \cdot \rangle_{{\pi}, {\color{blue} s_2}{\color{red} s_1}} $ 
this gives us as the only possibility $ d= 1 $, that is 
\begin{equation}{\label{onlypossi}}
[ \Delta_{{\pi}, \mathbb R}( {\color{blue} s_2}{\color{red} s_1})]= 
[L_{{\pi}, \mathbb R}({\color{blue} s_2}{\color{red} s_1})] + 
[L_{{\pi}, \mathbb R}({\color{red} s_1}{\color{blue} s_2}{\color{red} s_1})].
\end{equation}
Similarly we have 
\begin{equation}{\label{onlypossi2}}
[ \Delta_{{\pi}, \mathbb R}( {\color{red} s_1}{\color{blue} s_2})]= 
[L_{{\pi}, \mathbb R}({\color{red} s_1} {\color{blue} s_2})] +   
[L_{{\pi}, \mathbb R}({\color{red} s_1}{\color{blue} s_2}{\color{red} s_1})].
\end{equation}
Let us now turn to $[ \Delta_{{\pi}, \mathbb R}( {\color{blue} s_2})] $ which by 
({\ref{cellu}) 
can be written as 
\begin{equation}
[ \Delta_{{\pi}, \mathbb R}( {\color{blue} s_2})]  = [ L_{{\pi}, \mathbb R}( {\color{blue} s_2})] + 
a[ L_{{\pi}, \mathbb R}( {\color{red} s_1}{\color{blue} s_2})] + 
b [L_{{\pi}, \mathbb R}({\color{blue} s_2}{\color{red} s_1})] +
c[L_{{\pi}, \mathbb R}({\color{red} s_1}{\color{blue} s_2}{\color{red} s_1})]
\end{equation}
for some nonnegative integers $a, b,c$.  On the other hand, the sum formula ({\ref{sumformGR}})
gives in this case 
\begin{equation}
  \sum_{i>0} [ \Delta_{{\pi},\mathbb R}^i({\color{blue} s_2})] =
      [ \Delta_{{\pi},\mathbb R}( {\color{red} s_1}{\color{blue} s_2} )]  +
 [ \Delta_{{\pi},\mathbb R}({\color{blue} s_2} {\color{red} s_1} )]  = [ L_{{\pi},\mathbb R}( {\color{red} s_1}{\color{blue} s_2} )]  +
 [ L_{{\pi},\mathbb R}({\color{blue} s_2} {\color{red} s_1} )] +  2 [ L_{{\pi},\Bbbk}({\color{red} s_1}{\color{blue} s_2} {\color{red} s_1} )]
\end{equation}
where we for the second equality used ({\ref{onlypossi}}) and ({\ref{onlypossi2}}).
This information gives us that $ a= b = 1 $ but it does not determine the value of $ c$ since the filtration 
$ \{ \Delta_{{\pi},\mathbb R}^i({\color{blue} s_2}) \} $ 
may have one or two nonzero terms.
In fact it has two terms, as one sees via the restriction of the filtration to 
$ \Delta_{{\color{red} s_1 \color{blue} s_2 \color{red} s_1\color{blue}},\mathbb R }({\color{blue} s_2}) $. Indeed we have
$ \dim \Delta_{{\color{red} s_1 \color{blue} s_2 \color{red} s_1\color{blue}},\mathbb R }({\color{blue} s_2}) =1 $
and 
\begin{equation}
 \sum_{i >0} \dim \Delta^i_{{\color{red} s_1 \color{blue} s_2 \color{red} s_1\color{blue}},\mathbb R }({\color{blue} s_2}) =
\dim  \Delta_{{\color{red} s_1 \color{blue} s_2 \color{red} s_1\color{blue}},\mathbb R}( {\color{red} s_1}{\color{blue} s_2} ) +
\dim  \Delta_{{\color{red} s_1 \color{blue} s_2 \color{red} s_1\color{blue}},\mathbb R}
({\color{blue} s_2} {\color{red} s_1} )  = 2.
\end{equation}
and so there must be two terms in the filtration for
$ \Delta_{{\color{red} s_1 \color{blue} s_2 \color{red} s_1\color{blue}},\mathbb R }({\color{blue} s_2}) $, and
hence also in the filtration for $  \Delta_{{\pi},\mathbb R}({\color{blue} s_2})  $.
We then deduce that $ c=1$.

\medskip
We now consider $[ \Delta_{{\pi}, \mathbb R}( {\color{red} s_1})] $ which is 
treated almost the same way. 
We have that 
\begin{equation}
[ \Delta_{{\pi}, \mathbb R}( {\color{red} s_1})]  = [ L_{{\pi}, \mathbb R}( {\color{red} s_1})] + 
a_1 [ L_{{\pi}, \mathbb R}( {\color{red} s_1}{\color{blue} s_2})] + 
b_1 [L_{{\pi}, \mathbb R}({\color{blue} s_2}{\color{red} s_1})] +
c_1 [L_{{\pi}, \mathbb R}({\color{red} s_1}{\color{blue} s_2}{\color{red} s_1})]
\end{equation}
and the sum formula gives 
\begin{equation}
  \sum_{i>0} [ \Delta_{{\pi},\mathbb R}^i({\color{red} s_1})] =
      [ \Delta_{{\pi},\mathbb R}( {\color{blue} s_2}{\color{red} s_1} )]  +
 [ \Delta_{{\pi},\mathbb R}({\color{red} s_1} {\color{blue} s_2} )]  = [ L_{{\pi},\mathbb R}( {\color{blue} s_2}{\color{red} s_1} )]  +
 [ L_{{\pi},\mathbb R}({\color{red} s_1} {\color{blue} s_2} )] +  2 [ L_{{\pi},\Bbbk}({\color{blue} s_2}{\color{red} s_1} {\color{blue} s_2} )]
\end{equation}
and so $ a_1 = b_1 = 1$. In order to show that also $ c_1 =1 $ we must once again show that there are two
nonzero terms in the filtration $ \{  \Delta_{{\pi},\mathbb R}^i({\color{red} s_1}) \}$.
We here consider  
$ \Delta_{  {\color{red} s_1} {\color{blue} s_2} {\color{red} s_1},\mathbb R}({\color{red} s_1}) $
which is of dimension two, via the following light leaves basis
\begin{equation}
\begin{aligned}
  \begin{tikzpicture}[xscale=0.4, yscale=0.3]
\draw[red, very thick] (0+0,0)-- (0+0,2);
\draw[blue, very thick] (1+0,0)-- (1+0,2);
\node[circle,draw,blue,fill,minimum size=1mm, inner sep=0pt]  at (1+0,2) {};
\draw[red,very thick](0+0,2) to[out=90,in=180] (1+0, 3) to[out=0,in=90] (2+0,2);
\draw[red, very thick] (2+0,0)-- (2+0,2);
\draw[red, very thick] (1,3)-- (1,4);
\draw[red, very thick] (0+4,0)-- (0+4,2);
\node[circle,draw,red,fill,minimum size=1mm, inner sep=0pt]  at (0+4,2) {};
\draw[blue, very thick] (1+4,0)-- (1+4,2);
\node[circle,draw,blue,fill,minimum size=1mm, inner sep=0pt]  at (1+4,2) {};
\draw[red, very thick] (2+4,0)-- (2+4,4);
\end{tikzpicture} \, \, .  \end{aligned}
\end{equation}
With respect to this basis, 
the matrix of the bilinear form $  \langle \cdot, \cdot \rangle_{ \color{red}{s_1} } $ 
is
\begin{equation}
\begin{bmatrix}
-1  &   {\color{blue} \alpha_2}  \\
{ \color{blue} \alpha_2}    &  {\color{red} \alpha_1  \color{blue} \alpha_2}  \\
\end{bmatrix}
\sim
\begin{bmatrix}
-1  &   0  \\
0    &  ({\color{red} \alpha_1 + \color{blue} \alpha_2)\color{blue} \alpha_2}  \\
\end{bmatrix}   \stackrel{\varphi}{ \longrightarrow}
\begin{bmatrix}
-1  &   0  \\
0    & 2x^2  \\
\end{bmatrix} 
\end{equation}
where $ \sim  $ refers to the diagonalization procedure indicated in Theorem {\ref{sumformula}}.
The $ 2x^2 $ entry of this matrix tells us that there are two nonzero terms in the filtration
for $ \Delta_{  {\color{red} s_1} {\color{blue} s_2} {\color{red} s_1},\mathbb R}({\color{red} s_1}) $, and
so there are also two nonzero terms in the filtration for 
$     \Delta_{{\pi},\mathbb R} ( {\color{red} s_1})   $. We get from this 
that $c_1 = 1$ as claimed.

\medskip
Finally, using similar arguments we have that  
\begin{equation}
[ \Delta_{{\pi}, \mathbb R}( 1)]  =  [ L_{{\pi}, \mathbb R}( 1)] + 
 [ L_{{\pi}, \mathbb R}( {\color{red} s_1})] + 
 [ L_{{\pi}, \mathbb R}( {\color{blue} s_2})] + 
[ L_{{\pi}, \mathbb R}( {\color{red} s_1}{\color{blue} s_2})] + 
 [L_{{\pi}, \mathbb R}({\color{blue} s_2}{\color{red} s_1})] +
[L_{{\pi}, \mathbb R}({\color{red} s_1}{\color{blue} s_2}{\color{red} s_1})].
\end{equation}

\medskip
Using Theorem {\ref{firstmmain}}
one can also carry out these calculations in the characteristic $ p> 3 $ cases. In fact, one gets 
in these cases the same decomposition numbers.

\medskip
In the language of $p$-canonical bases, see \cite{JW}, we get that for $ A_2$ the $p$-canonical basis coincides with 
the Kazhdan-Lusztig basis. This is already proved in \cite{JW} using very different methods.

\medskip
In view of the results and conjectures of \cite{RW}, we note 
that 
it would be very interesting to generalize our results to the anti-spherical module of affine Weyl groups.

\sc Instituto de Matem\'atica y F\'isica, Universidad de Talca, Chile,
steen@inst-mat.utalca.cl.

\end{document}